\documentclass[11pt,reqno]{amsart}
\usepackage[margin=3cm]{geometry}
\usepackage{amsmath, amsfonts, amssymb, amscd, amsthm, amsbsy, amscd}
\usepackage{bbm, epsf, calc, enumerate, graphicx, verbatim, cite, setspace}
\usepackage[usenames,dvipsnames]{color}
\usepackage{mathrsfs, upgreek}
\usepackage{thmtools}
\usepackage{thm-restate}
\usepackage[margin=1cm]{caption}
\usepackage{mathtools}
\usepackage[colorlinks=true, linkcolor=blue, citecolor=red, pdfborder={0 0 0},pagebackref,hypertexnames=false]{hyperref}
\usepackage{cleveref}
\usepackage{accents}
\usepackage{wrapfig}
\usepackage{multirow}
\usepackage{float}
\usepackage{blkarray}
\usepackage{tikz}
\usepackage{yfonts}
\usepackage{bm}
\usepackage{stmaryrd}

\let\oldtocsubsection=\tocsubsection
\renewcommand{\tocsubsection}[2]{\hspace*{1.0cm}\oldtocsubsection{#1}{#2}}
\let\oldtocsubsubsection=\tocsubsubsection
\renewcommand{\tocsubsubsection}[2]{\hspace*{1.8cm}\oldtocsubsubsection{#1}{#2}}

\theoremstyle{plain}
 \newtheorem{theorem}{Theorem}[section]
 \newtheorem{lemma}[theorem]{Lemma}
 \newtheorem{prop}[theorem]{Proposition}
 \newtheorem{cor}[theorem]{Corollary}
 \newtheorem{conj}[theorem]{Conjecture}

 \newtheorem*{ack}{Acknowledgments}

\theoremstyle{definition}
 
 \newtheorem{assumption}[theorem]{Assumption}

 \newtheorem{claim}[theorem]{Claim}

\theoremstyle{remark}
 \newtheorem{remark}[theorem]{Remark}
 
\numberwithin{equation}{section}
\numberwithin{theorem}{section}

\newcommand\nc\newcommand
\nc\dmo\DeclareMathOperator
\nc\bs\boldsymbol
\nc{\labbr}[1]{{\sc\lowercase{#1}}}
\nc{\abbr}[1]{{\sc\uppercase{#1}}}
\nc{\ignore}[1]{}

\nc{\cf}{cf.\ }
\nc{\ie}{i.e.\ }
\nc{\eg}{e.g.\ }

\nc{\overeq}[1]{\stackrel{\text{\tiny #1}}{=}}
\nc{\overleq}[1]{\stackrel{\text{\tiny #1}}{\le}}

\nc{\filldi}{{\hfill$\Diamond$}}

\renewcommand*{\backref}[1]{\ifx#1\relax \else Page #1 \fi}
\renewcommand*{\backrefalt}[4]{
  \ifcase #1 \footnotesize{(Not cited.)}
  \or        \footnotesize{(Cited on page~#2.)}
  \else      \footnotesize{(Cited on pages~#2.)}
  \fi
}

\nc{\red}[1]{{\color{red} #1}}
\nc{\blue}[1]{{\color{blue} #1}}
\nc{\green}[1]{{\color{green} #1}}
\nc{\cyan}[1]{{\color{cyan} #1}}
\definecolor{purple}{rgb}{0.9,0,0.8}
\nc{\purple}[1]{{\color{purple} #1}}
\definecolor{gray}{rgb}{0.5,0.5,0.5}
\nc{\gray}[1]{{\color{gray} #1}}
\nc{\note}[1]{{\blue{\textup{[#1]}}}}
\nc{\nick}[1]{{\color{blue} #1}}
\nc{\nickC}[1]{{\bf \emph{\nick{[#1]}}}}
\nc{\comN}[1]{{\bf \emph{\nick{[#1]}}}}
\nc{\Nick}[1]{\emph{\nick{[#1]}}}

\nc{\al}{\alpha} \nc{\Al}{\Alpha}
\nc{\be}{\beta} \nc{\Be}{\Beta}
\nc{\gam}{\gamma} \nc{\Gam}{\Gamma}
\nc{\eps}{\varepsilon}
\nc{\epp}{\epsilon}
\nc{\vep}{\varepsilon}
\nc{\de}{\delta} \nc{\De}{\Delta}
\nc{\lam}{\lambda} \nc{\Lam}{\Lambda}
\nc{\ka}{\kappa} \nc{\Ka}{\Kappa}
\nc{\vpi}{\varpi} \nc{\vph}{\varphi}
\nc{\om}{\omega} \nc{\Om}{\Omega}
\nc{\si}{\sigma} \nc{\Si}{\Sigma}

\nc{\A}{\mathbb{A}}
\nc{\B}{\mathbb{B}}
\nc{\C}{\mathbb{C}}
\nc{\D}{\mathbb{D}}
\nc{\E}{\mathbb{E}}
\nc{\F}{\mathbb{F}}
\nc{\G}{\mathbb{G}}
\nc{\bI}{\mathbb{I}}
\nc{\J}{\mathbb{J}}
\nc{\K}{\mathbb{K}}

\nc{\M}{\mathbb{M}}
\nc{\N}{\mathbb{N}}
\nc{\bO}{\mathbb{O}}

	\renewcommand{\P}{\mathbb{P}}
\nc{\Q}{\mathbb{Q}}
\nc{\R}{\mathbb{R}}
\nc{\bS}{\mathbb{S}}
\nc{\T}{\mathbb{T}}
\nc{\U}{\mathbb{U}}
\nc{\V}{\mathbb{V}}
\nc{\W}{\mathbb{W}}
\nc{\X}{\mathbb{X}}
\nc{\Y}{\mathbb{Y}}
\nc{\Z}{\mathbb{Z}}

\nc{\cA}{\mathcal{A}}
\nc{\cB}{\mathcal{B}}
\nc{\cC}{\mathcal{C}}
\nc{\cD}{\mathcal{D}}
\nc{\cE}{\mathcal{E}}
\nc{\cF}{\mathcal{F}}
\nc{\cG}{\mathcal{G}}
\nc{\cH}{\mathcal{H}}
\nc{\cI}{\mathcal{I}}
\nc{\cJ}{\mathcal{J}}
\nc{\cK}{\mathcal{K}}
\nc{\cL}{\mathcal{L}}
\nc{\cM}{\mathcal{M}}
\nc{\cN}{\mathcal{N}}
\nc{\cO}{\mathcal{O}}
\nc{\cP}{\mathcal{P}}
\nc{\cQ}{\mathcal{Q}}
\nc{\cR}{\mathcal{R}}
\nc{\cS}{\mathcal{S}}
\nc{\cT}{\mathcal{T}}
\nc{\cU}{\mathcal{U}}
\nc{\cV}{\mathcal{V}}
\nc{\cW}{\mathcal{W}}
\nc{\cX}{\mathcal{X}}
\nc{\cY}{\mathcal{Y}}
\nc{\cZ}{\mathcal{Z}}
\nc{\co}{\scriptsize$\cO$\normalsize}

\nc\notni{\not\owns}
\def\({\left(}
\def\){\right)}
\nc{\ii}{\mathrm{i}}
\renewcommand{\d}{\,d}
\nc{\ls}{\lesssim}
\nc{\gs}{\gtrsim}
\def \lf {\lfloor}
\def \rf {\rfloor}

\dmo{\lls}{\,\ls\,}
\dmo{\ggs}{\,\gs\,}
\nc{\cond}{\,\middle|\,}
\dmo{\I}{{I}}
\dmo{\II}{{II}}
\dmo{\III}{III}
\dmo{\argmin}{arg\,min}
\dmo{\argmax}{arg\,max}
\dmo{\area}{area}
\dmo{\diag}{diag}
\dmo{\diam}{diam}
\dmo{\esssup}{ess\,sup}
\dmo{\essinf}{ess\,inf}
\dmo{\im}{im} 
\dmo{\loc}{loc}
\dmo{\sgn}{sgn}
\dmo{\supp}{supp}
\dmo{\vol}{vol}
\dmo{\dist}{{dist}}
\dmo\per{per}

\DeclareFontFamily{U}{mathx}{\hyphenchar\font45}
\DeclareFontShape{U}{mathx}{m}{n}{
      <5> <6> <7> <8> <9> <10>
      <10.95> <12> <14.4> <17.28> <20.74> <24.88>
      mathx10
      }{}
\DeclareSymbolFont{mathx}{U}{mathx}{m}{n}
\DeclareFontSubstitution{U}{mathx}{m}{n}
\DeclareMathAccent{\widecheck}{0}{mathx}{"71}
\DeclareMathAccent{\wideparen}{0}{mathx}{"75}
\nc{\ol}{\overline}

\nc{\ul}{\underline}
\nc{\wt}{\widetilde}
\nc{\wh}{\widehat}
\nc{\wch}{\widecheck}

\nc{\decto}{\downarrow}
\nc{\incto}{\uparrow}
\nc{\la}{\leftarrow}
\nc{\ra}{\rightarrow}
\nc{\La}{\Leftarrow}
\nc{\Ra}{\Rightarrow}
\nc{\lla}{\longleftarrow}
\nc{\lra}{\longrightarrow}
\nc{\Lla}{\Longleftarrow}
\nc{\Lra}{\Longrightarrow}
\nc{\iffa}{\Leftrightarrow}
\nc{\hra}{\hookrightarrow}
\nc{\inj}{\hookrightarrow}
\nc{\surj}{\twoheadrightarrow}
\nc{\bij}{\leftrightarrow}
\nc{\longbij}{\longleftrightarrow}
\nc{\dirlim}{\lim_{\lra}}
\nc{\invlim}{\lim_{\lla}}
\nc{\cra}{\stackrel{\sim}{\ra}}

\dmo{\adj}{adj}
\dmo{\Hom}{Hom}
\dmo{\proj}{proj}
\dmo{\tr}{tr}
\dmo{\Tr}{Tr}
\nc{\Span}{\operatorname{span}}

\dmo{\rank}{rank}

\dmo{\real}{Re}

\dmo{\Var}{Var}
\dmo{\Corr}{Corr}
\dmo{\Cov}{Cov}
\dmo{\1}{\mathbf{1}}
\dmo{\ind}{\mathbbm{1}}
\nc{\eqd}{\stackrel{\text{\tiny $d$}}{=}}
\nc{\dto}{\stackrel{d}{\to}}

\nc{\iid}{\abbr{iid}}
\dmo{\Ber}{Ber}
\dmo{\Bin}{Bin}
\dmo{\Poi}{Poi}
\nc{\ER}{Erd\H{o}s--R\'enyi }

\dmo{\Sparse}{Sparse}
\dmo{\Comp}{Comp}
\dmo{\Incomp}{Incomp}
\dmo{\row}{{row}}
\dmo{\col}{{col}}
\nc{\sphere}{\mathbb{S}}
\nc{\spheren}{{\sphere^{n-1}}}
\nc{\sphereN}{{\sphere^{N-1}}}
\nc{\ball}{\mathbb{B}}

\nc{\haotian}[1]{{\color{violet} #1}}
\nc{\comH}[1]{{\bf \emph{\haotian{[#1]}}}}
\nc{\mc}[1]{\mathcal{#1}}
\nc{\bone}{ {\mathbbm{1}} }
\nc{\floor}[1]{\lfloor #1\rfloor}
\nc{\ee}{\varepsilon}
\nc{\RoNum}[1]{\MakeUppercase{\romannumeral #1}}
\dmo\cue{\abbr{CUE}}
\nc\brw{\abbr{BRW}}
\dmo\upP{\mathbf{P}}
\nc\oop{o_{_{\upP}}\!}
\nc\OOP{O_{_{\upP}}\!}
\nc\speed{c_0}
\nc{\rot}{{\upsilon}}
\nc{\crot}{{\cU}}
\nc{\arr}{{\tau}}
\nc{\carr}{{\cT}}
\nc\hcarr{\widehat{\carr}}
\nc\Pcarr{\P^{\,\carr}\!}
\nc\Ecarr{\E^{\,\carr}\!}
\nc\Phcarr{\P^{\,\hcarr}\!}
\nc\Ehcarr{\E^{\,\hcarr}\!}
\nc\Pcars{\P^{\,\carr}_{*}}
\nc\Ecars{\E^{\,\carr}_{*}}
\nc\Phcars{\P^{\,\hcarr}_{*}}
\nc\Ehcars{\E^{\,\hcarr}_{*}}

\nc{\rev}{\overleftarrow}

\title[The maximum of a Poissonian log-correlated field]{On the maximum of Poissonian log-correlated fields} 

\author[N.\ A.\ Cook]{Nicholas A.\ Cook$^\star$}
\thanks{${}^{\star,\ddagger}$ Supported in part by NSF grant DMS-2154029.}
\address{$^\star$Department of Mathematics, Duke University, 120 Science Dr, Durham, NC 27710, USA}\email{nicholas.cook@duke.edu}
\author[H.\ Gu]{Haotian Gu$^{\ddagger}$}
\address{$^{\ddagger}$ Department of Mathematics, UCLA, Los Angeles, CA, USA.}\email{haotian.gu@math.ucla.edu}

\date{\today}

\begin{document}

\begin{abstract}
    We consider stationary log-correlated random fields of the form
    \[
    X_n(t):= \sum_{j\le n} f(e^{\arr_j} t+ \rot_{j})\,,\quad t\in \R
    \]
    for a centered 1-periodic analytic function $f:\R\to\R$, where $(\arr_j)_{j\ge1}$ are the arrival times of a rate-1 Poisson process on $\R_+$ and $(\rot_j)_{j\ge1}$ are iid uniform random variables in $[0,1)$. 
    We show that with high probability, 
    \[
         \max_{t\in[-1,1]}X_n(t) = c_1 n + c_2 \arr_n- (c_3+o(1))\log n
    \]
    where  $c_1,c_2,c_3$ are positive 
    constants depending on $f$. 
    In particular, the random frequency $\tau_n$ produces a central limit fluctuation of order $\sqrt{n}$, together with a logarithmic correction different from the ubiquitous Bramson correction for models with strong comparisons to Gaussian fields, such as homogeneous branching random walks and log-characteristic polynomials of C$\beta$E matrices.
    Our proof makes a connection with branching random walks in random time-environment, where coarse-scale fluctuations of the arrival times $\tau_j$ play the role of a quenched random environment, and fine-scale fluctuations are leveraged to establish a branching decorrelation structure.
    We conjecture analogous asymptotics for the characteristic polynomials of random permutation matrices, and expect the same mechanism to arise in other Poissonian fields such as characteristic polynomials of sparse random matrices.
\end{abstract}

\maketitle

\setcounter{tocdepth}{1}
\tableofcontents

\section{Introduction}

\subsection{Motivation}
\label{sec:motivation}

Large values of logarithmically correlated random fields (LCFs) have been intensely studied in recent years due to connections with areas such as random matrices, multiplicative chaos measures and Liouville quantum gravity, analytic number theory, branching processes and reaction-diffusion PDE. 
We refer to \cite{DRS2017LCGF, bovier2017gaussian, bailey2022maxima, ryzhik2023lausanne} and references therein for an introduction to this vast area. 

A fundamental example goes back to the pioneering work of Bramson \cite{bramson1978maximal,bramson1983convergence} on the maximum of binary branching Brownian motion (BBM), a particle system on the real line with standard Brownian displacements and binary branching at rate one. The maximal displacement $M_t$ satisfies 
\begin{equation} \label{BBM}
    M_t
    =\sqrt{2}\,\left(t-\frac{3}{4}\log t\right)+b_t
\end{equation}
for $b_t$ converging in law as $t\to\infty$.
This logarithmic correction and order-one fluctuation of the maximum appear, with suitable constants, in the extremal statistics of many other LCFs.

The model considered in the present work is motivated by results for characteristic polynomials of random matrices that we now recall, together with some connections to analytic number theory. 
A detailed review of the literature is postponed to $\mathsection$\ref{sec:literature}.

\subsubsection{The CUE field and the zeta function}

For $M_N$ drawn from the circular unitary ensemble (CUE), i.e.\ the Haar measure on the group of $N\times N$ unitary matrices, let
\begin{equation}    \label{LCF.CUE}
    P_N(t) := \det( e(t) - M_N)\,,\qquad t\in \R/\Z,
\end{equation}
where we use the common abbreviation $e(t):=e^{2\pi it}$.
It was conjectured in \cite{fyodorov2012freezing} and finally proved in \cite{paquette2022extremal} that
\begin{equation}    \label{FHK}
    \max_{t\in \R/\Z} \log|P_N(t)| = \log N - \frac34\log\log N + \xi_N
\end{equation}
for a sequence $\xi_N$ converging in law. 
The work \cite{paquette2022extremal} built on a sequence of works establishing \eqref{FHK} to precision  $o(\log N)$ \cite{ABB2017maximum}, $o(\log\log N)$ \cite{paquette2018maximum} and $\OOP(1)$ \cite{chhaibi2018maximum}.
We note that \cite{chhaibi2018maximum,paquette2022extremal} cover the more general class of circular beta ensembles. 

The conjecture \eqref{FHK} was motivated by the famous Montgomery--Dyson analogy between $P_N$ and the Riemann zeta function at a randomly chosen short interval at height of order $e^N$ on the critical line. In that setting the asymptotic was established in \cite{ABR1,ABR2}:
\begin{equation}\label{FHK2}
    \max_{|h|\le 1}\log\left|\zeta\left(\frac12+i(\tau+h)\right)\right| = \log\log T-\frac34\log\log\log T+\wt\xi_T\qquad\text{as }T\to\infty
\end{equation}
with $\tau\sim\mathrm{Unif}[T,2T]$ and a tight family $\wt\xi_T$ enjoying uniform tail estimates. 
Note with $T=e^N$, \eqref{FHK2} has the same asymptotic form as \eqref{FHK}.

The CUE model \eqref{FHK} enjoys strong comparisons with Gaussian fields. 
A CLT showing that $(2/\log N)^{1/2}\log |P_N(t)|$ converges to a standard Gaussian for fixed $t\in \R/\Z$ was established in \cite{KeSn}, with a multivariate CLT for values of the field at multiple points proved in \cite{HKO}.
A closely related, multivariate CLT for eigenvalue counts over fixed arcs is proved by Wieand \cite{wieand2002eigenvalue}.
In fact, the Gaussian tail holds out to the large deviation scale $\sim \log N$; see \cite[Lemma 2.1]{ABB2017maximum}.
The corresponding Gaussian behavior for the log-zeta function at one point is given by the Selberg's celebrated CLT; see \cite{radziwill2017selberg} for a modern proof and \cite{bourgade2010mesoscopic} for a multidimensional extension.

The logarithmic terms $-\frac34\log\log N,\, -\frac34\log\log\log T$ in \eqref{FHK} and \eqref{FHK2} show the familiar Bramson's correction for extremes of LCFs, with the same $3/4$ coefficient appearing in the normalization of \eqref{BBM}.

\subsubsection{A canonical model of Poissonian LCF}

The present work is motivated by a result in \cite{cook2020maximum} on the maximum for a field with \emph{Poissonian} features, where $M_N$ in \eqref{LCF.CUE} is instead drawn from the discrete subgroup of permutation matrices. That is, with $\pi$  a uniform random permutation on $N$ labels, let $A_\pi$ be the associated $N\times N$ permutation matrix, and set
\begin{equation}    \label{def:QN}
    Q_N(t):= \det(e(t)- A_\pi)\,,\qquad t\in \R/\Z.
\end{equation}
In contrast to $P_N$, the field $Q_N$ is not stationary under rotations of the torus, with the distribution of $Q_N(t)$ depending strongly on Diophantine approximation properties of $t$. Nonetheless, CLTs were proved in 
\cite{HKOS,DaZe} for the evaluation of $\log|Q_N|$ at points that are suitably generic in an arithmetic sense. 
It was shown in \cite{cook2020maximum} that
\begin{equation}    \label{CoZe}
    \max_{t\in\R/\Z} \log|Q_N(t)| = (\speed+\oop(1))\log N
\end{equation}
with a constant $\speed\approx 0.652$ determined by the equation $\lam^*(\speed)=1$, where $\lam^*$ is the Cram\'er large deviation rate function
of $\log|1-e(\rot)|$ for $\rot$ uniform in $[0,1)$. 
In \eqref{CoZe} and in the sequel we write $\oop(1)$ for a sequence converging to zero in probability.

The maximum of $\log|Q_N|$ differs from that of $\log|P_N|$ already to leading order.
That the leading order speed $\speed$ should be different from the speed 1 in \eqref{FHK} is already suggested by the fact that, while $\log|Q_N(t)|$ enjoys a CLT (for generic $t$), the tails are \emph{non-Gaussian} at the large deviation scale $\sim \log N$. 
In fact, the speed $\speed$ is determined by large deviations conditional on a typical distribution of cycle lengths, rather than by the unconditional pointwise tails; see \cite[$\mathsection2$]{cook2020maximum}.
This phenomenon also plays a key role in the present work.

Like the CUE field \eqref{LCF.CUE}, the permutation model \eqref{logQN} is also connected with an object from analytic number theory. 
In forthcoming work \cite{ABK},
Arguin, Bourgade and Koukoulopoulos consider the maximum of the twisted divisor function
\begin{equation}    \label{twisted-div}
    \tau(n,\xi):= \sum_{d|n} d^{i\xi}
\end{equation}
over $\xi\in[1,2]$, for $n$ drawn randomly from $[1,x]$ with $x\to\infty$. 
It is well known that the cycle factorization of a large random permutation shares many distributional properties with the prime factorization of a large random integer, so \eqref{twisted-div} can be expected to share features of \eqref{def:QN} dealt with in \cite{cook2020maximum} and the present work.

Indeed, the starting point in \cite{cook2020maximum} is to note that $Q_N$ is completely determined by the cycle structure of $\pi$: with $C_\ell(\pi)$ the number of cycles of length $\ell$, we have
\begin{equation}    \label{logQN}
    \log |Q_N(t)| = \sum_{\ell=1}^N C_\ell(\pi) f(\ell t)
\end{equation}
with $f(u)=\log|1-e(-u)|$.
A result of Arratia and Tavar\'e \cite{arratia_cycle_1992} allows, up to a small error, to replace the cycle counts $C_\ell(\pi)$ with a sequence of independent Poisson variables $Z_\ell$ of mean $\E Z_\ell=1/\ell$, leading to the model\footnote{The passage from $\log|Q_N|$ to $Y_N$ is not straightforward for proving the lower bound in \eqref{CoZe} due to the singularities of $f$; see \cite[$\mathsection$2]{cook2020maximum} for further discussion.}
\begin{equation}    \label{series.Poisson}
    Y_N(t) = \sum_{\ell=1}^N Z_\ell f(\ell t) = \sum_{j : L_j\le N} f(L_jt)\,,\qquad t\in \R/\Z
\end{equation}
where $1\le L_1\le L_2\le\cdots$ are the ordered points of a Poisson point process on $\Z_+$ with 1-point intensity function $\mu(\ell)=1/\ell$.
For the twisted divisor function \eqref{twisted-div}, its factorization for square-free $n$ gives
\begin{equation}\label{series.Twisted}
\log|\tau(n;\xi)|=\sum_{p\mid n}\log|1+e^{i\xi\log p}|\,,
\end{equation}
again a sum of periodic functions with frequencies $\log p$ selected by the prime factors of $n$.
Kubilius's classical independent-prime model, with selection probabilities $1/p$, suggests a similar sparse sampling of frequencies.

In this article, we introduce a canonical (stationary) Poissonian field inspired by \eqref{series.Poisson}, and prove a refined asymptotic for its maximum.
The related structure in \eqref{series.Twisted} suggests that this model is relevant beyond random permutation matrices. 

Let $\carr=(\arr_j)_{j\geq 1}$ be the arrival times of a rate-1 Poisson process on $\R_+$ and let $\crot=(\rot_j)_{j\ge1}$ be a sequence of iid random variables uniform in $[0,1)$, independent of $\carr$.
For a smooth 1-periodic function $f$ on $\R$ with $\E f(\rot_1)=0$,
we consider the stationary random field 
\begin{align}\label{def:field}
    X_n(t):=\sum_{j\le n} f(e^{\arr_j} t+\rot_j),
    \qquad t\in \R\,.
\end{align} 
(Note that $\sum_j \delta_{e^{\tau_j}}$ is a Poisson point process on $[1,\infty)$ of intensity $\mu(x)=1/x$.)
We make the following additional assumption on $f$:

\begin{assumption}  \label{assu:f}
    The distribution function of $f(\rot_1)$ is H\"older continuous.
\end{assumption}

For concreteness one may consider $f(x)=\cos(2\pi x)$, for which the distribution function of $f(\rot_1)$ is $\frac12$-H\"older continuous.

The field \eqref{def:field} is centered by our assumption $\E f(\rot_1)=0$, with covariance
\begin{align}   \label{Cov-sim0}
    \Cov(X_n(s),X_n(t)) 
    &\sim  \|f\|_2^2\,\min\Big\{ n\,,\,  \log\frac1{|s-t|}\Big\} +O(1) \,,
    \qquad s, t\in\R\,.
\end{align} 
See Appendix \ref{sec:covariance}, where we give a refined asymptotic in terms of the Fourier decomposition of $f$, with quantitative errors.
Thus the field is log-correlated down to the scale $e^{-n+o(n)}$, and is strongly correlated below that scale.

\subsection{Main result}
\label{sec:result}

Our main result for the maximum of $X_n$ involves some numerical constants depending on $f$ that we now define.
Denote the log-Laplace transform of $f(\rot_1)$ and its Legendre transform by 
\begin{align}\label{def:lambda}
    \lambda(\theta) := \log\E e^{\theta f(\rot_1)}
    \quad\text{ and }\quad
    \lambda^*(x) =\sup_{\theta \in \R}\{x\theta-\lambda(\theta)\}.
\end{align} 
Since $\lambda^*$ is strictly increasing and onto $\R_{+}$, there are unique $\speed,\theta_*>0$ such that
\begin{align}\label{def:speed}
    \lambda^*(\speed)=\speed\theta_*-\lambda(\theta_*)=1, 
    \qquad
    \lambda'(\theta_*)=\speed.
\end{align}
It is shown in \cite[Theorem 1.1]{MaMi:barrier} that there exists a symmetric, convex function $\gamma:\R\to\R^+$ with $\gamma(\beta)>\gamma(0)=1/2$ for all $\beta\ne 0$ such that for independent standard Brownian motions $B_t,W_t$, we have
\begin{equation}\label{def:gamma}
\lim_{t\to\infty}\frac{1}{\log t}\log\P(B_s+1\geq\beta W_s,\;\forall\,s\in[0,t]\mid W)= -\gamma(\beta)
\end{equation}
almost surely. We stress that $\gamma$ is universal, i.e.\ independent of $f$.\footnote{The existence of $\gamma$ is obtained by an application of Kingman's subadditive ergodic theorem, and no further properties of $\gamma$ are known at the time of this writing. The authors of \cite{MaMi:barrier} conjecture that $\gamma$ is strictly convex and grows quadratically at infinity.}
Define 
\begin{equation}
\label{def:phi*}
\phi_*:=\frac12+ 2\gamma\bigg(\frac1{\theta_*\sqrt{\lambda''(\theta_*)}}\bigg).
\end{equation}
Recall that $\oop(1)$ denotes a quantity that converges to zero in probability as $n\to \infty$.

\begin{theorem}
\label{thm:main}
We have
\begin{equation}
    \label{main-asymp}
    \max_{t\in [-1,1]} X_n(t) =  c_1n + c_2\arr_n -(c_3+\oop(1))\log n
\end{equation}
where $c_1=\lam(\theta_*)/\theta_*$, $c_2=1/\theta_*$ and $c_3=\phi_*/\theta_*$.
\end{theorem}

\begin{remark} 
We take the interval $[-1,1]$ for convenience;
the proof gives the same asymptotic for the maximum over any interval of fixed positive length. (Such an extension is already immediate from \Cref{thm:main} for intervals of bounded length at least 2, using stationarity of the field and the union bound, and can be deduced for shorter intervals by rescaling $t$ and discarding a small number of low-frequency terms from the series \eqref{def:field}.) 
\end{remark}

Since $\arr_n$ is a sum of $n$ iid exponential variables of rate 1, 
the leading order of the maximum is $\sim  (c_1+c_2)n=\speed n$, with fluctuations at scale  $\sqrt{n}$.
We thus obtain a central limit theorem:

\begin{cor}
    As $n$ goes to infinity, we have the following central limit theorem for the maximum:
    \[
    \theta_*\frac{\max_{t\in[-1,1]}X_n(t)-\speed n}{\sqrt{n}} \implies \mc{N}(0,1).
    \]
\end{cor}

However, \Cref{thm:main} gives a more detailed picture: the square-root fluctuations are entirely encapsulated by the ``global'' statistic $\arr_n$, the highest log-frequency in the sum \eqref{def:field} defining $X_n$, and after subtracting $c_1n+c_2\arr_n$ there is a logarithmic shift $-(c_3+\oop(1))\log n$ that is deterministic to leading order. 

A similar structure was established for the maximum of branching random walks in a time-inhomogeneous random environment by Mallein--Mi{\l}o{\'s} in \cite{MaMi:branching}, and indeed the main idea for the proof of \Cref{thm:main} is to make a connection with such branching models, by viewing a coarse-graining of the arrivals process $\carr$ as a ``random environment''. We expect that 1) CLT fluctuations from a global statistic like $\arr_n$, and 2) a log-correction involving the function $\gamma$ from \eqref{def:gamma} are generic features of extremes of Poissonian LCFs (see $\mathsection$\ref{sec:literature} for other related models).

However, proofs of lower bounds on maxima in \cite{MaMi:branching} and in the subsequent work \cite{Kriechbaum:time} of Kriechbaum establishing tightness
rely on the independence between descendant sub-family-trees conditional on the environment, whereas our field has no independent subtrees: each frequency and phase affects the field at every spatial point.
We overcome this difficulty by leveraging the fine-scale randomness of $\carr$ to prove approximate independence estimates for path events at well-separated points. 
See $\mathsection$\ref{sec:ideas} for more details.

\subsection{Conjectures for related Poissonian fields}

In the random field $X_n$ of \eqref{def:field}, the number $n$ of terms in the sum is a deterministic asymptotic parameter, while the largest frequency $e^{\arr_n}$ contributing to $X_n$ is random and enters into the asymptotic \eqref{main-asymp} for the maximum.
A natural variant of the model is to instead take a deterministic frequency cutoff $N$ for the asymptotic parameter, and consider the field    
\begin{equation}
    \label{def:field-freq}
    \cY_N(t):= \sum_{j\ge1\,:\, e^{\arr_j}\le N} f(e^{\arr_j}t+\rot_j)
    \,,\qquad t\in \R\,.
\end{equation}
Then the number of terms in the sum is the random variable
\[
Z_N = \sum_{j\ge1} \bone(e^{\arr_j}\le N)\,.
\]
We conjecture that the asymptotic \eqref{main-asymp} should carry over to this setting, with $c_1$ multiplying the (now random) number of terms in the sum, and $c_2$ multiplying the logarithm of the (now essentially deterministic) largest frequency in the sum. 

\begin{conj}    \label{conj:freq}
    We have
    \begin{equation}    \label{freq.conj}
        \max_{t\in[-1,1]} \log\cY_N(t) = c_1 Z_N  + c_2\log N - (c_3+\oop(1))\log\log N
    \end{equation}
    where $c_1,c_2,c_3>0$ are as in \Cref{thm:main}.
\end{conj}

Note that \Cref{thm:main} implies \eqref{freq.conj} holds along a subsequence $N_j\to\infty$. 

Note that $Z_N\sim \log N$ to leading order, with fluctuations of order $\sqrt{\log N}$. The asymptotic \eqref{freq.conj} thus has the same general form as \eqref{main-asymp} (with $\log N$ in place of $n$), but now it is the global statistic $Z_N$ -- the number of terms in the sum -- that has square-root fluctuations, while the largest log-frequency is $\log N + O_{\bf P}(1)$.

From \eqref{logQN}--\eqref{series.Poisson} we see that the field $\cY_N$ is closer to $Y_N$ in \eqref{series.Poisson} (hence also the permutation model $\log|Q_N|$ in \eqref{def:QN}), compared to $X_{\log N}$ in \eqref{def:field}.
Recall that $C_\ell(\pi)$ denotes the number of cycles of length $\ell$ in a permutation $\pi$. Let 
\[
C(\pi):=\sum_{\ell=1}^N C_\ell(\pi)
\]
denote the total number of cycles. Based on \Cref{conj:freq} we make the following:

\begin{conj}[Maximum for the permutation model] \label{conj:perm}
    For $Q_N$ as in \eqref{def:QN}, 
    \begin{equation}    \label{QN.conj}
        \max_{t\in\R/\Z} \log|Q_N(t)| = c_1C(\pi) + c_2\log N - (c_3+\oop(1))\log\log N
    \end{equation}
    where $c_1,c_2,c_3>0$ are as in \Cref{thm:main} for the case that $f(u)=\log|1-e(u)|$. 
    For this choice of $f$ we have $\lambda(x)=x\log 2-\frac12\log\pi+\log\Gamma((1+x)/2)-\log\Gamma(1+\frac x2)$, which, together with \eqref{def:speed} and \eqref{def:phi*} give the numerical values
    $c_1=\lambda(\theta_*)/\theta_*\approx 0.567, c_2=1/\theta_*\approx 0.085,$ and $c_3=(\frac12+2\gamma(c_4))/\theta_*>\frac{3}{2\theta_*}$ for $c_4
    =(\theta_*\sqrt{\lambda''(\theta_*)})^{-1}
    \approx 1.478$ (recall that the function $\gamma$ is unknown).
\end{conj}
It is classical that $C(\pi) = (1+\oop(1))\log N$, and one verifies the constant $c_0$ in \eqref{CoZe} is equal to $c_1+c_2$, so \eqref{QN.conj} is consistent with \eqref{CoZe}.
A recent preprint \cite{lu2026spectral} claims to establish this conjecture up to the correct $O_{\P}(\log\log N)$ order, without the sharp constant $c_3$.

The differences between $Y_N$ and $\cY_N$ are as follows, and we expect none of them to affect the asymptotic behavior of the maximum, apart from appropriate modification of constants based on the Laplace transform of $f$ (see \eqref{def:lambda}--\eqref{def:phi*}). 
\begin{enumerate}
    \item Instead of $f=\log|1-e(-u)|$, which has singularities, $\cY_N$ takes $f$ to be a bounded smooth 1-periodic function on $\R$ (satisfying \Cref{assu:f}).
    \item $\cY_N$ is stationary due to perturbing the argument of $f$ for each summand by independent uniform random shifts $\rot_j$.
    \item $\cY_N$ replaces the discrete Poisson point process of frequencies $L_j\in\Z_+$ in $Y_N$ by a Poisson point process on the continuum $[1,+\infty)$ of intensity $\mu(x)=1/x$. 
\end{enumerate}
Indeed, the analysis in \cite{cook2020maximum} showed that with high probability the maximum for \eqref{series.Poisson} is attained at ``minor arc'' points $t$ that are not well approximated by rationals with small denominator. For such points, $L_j t$ rapidly equidistributes in $\R/\Z$ and does not come too close to the logarithmic singularities of $f$, so $f(L_jt)$ can be modeled by a bounded function evaluated at a uniform point in $\R/\Z$. (The joint distribution of $(f(L_js),f(L_jt))$ further depends on arithmetic properties of the difference $s-t$.)

\subsection{Related literature}
\label{sec:literature}

Apart from results for the CUE and permutation models discussed in $\mathsection$\ref{sec:motivation}, there is a large literature establishing asymptotics to varying levels of precision for the characteristic polynomial of random matrices. 
These include Hermitian random matrices \cite{fyodorov2016GUE,fyodorov2016moments,claeys2021much,LaPa,BLZ}, iid matrices \cite{lambert2020maximum, CiLa}, and 2-dimensional Coulomb gases \cite{LLZ:Coulomb,Peilen:Coulomb}.
There is also work on the related problem of showing that the modulus of the characteristic polynomial converges in some sense to a Gaussian multiplicative chaos measure \cite{webb2015characteristic,berestycki2018hermitian,nikula2020multiplicative,Forkel2021compact,junnila2024multiplicative,lambert2024subcritical,kivimae2024gaussian,Chhaibi2025circle}.
Parallel to work on the Fyodorov--Hiary--Keating conjecture \eqref{FHK}, corresponding asymptotics have been established for the Riemann zeta function, or random models of the zeta function \cite{arguin2017maxima,arguin2019maximum,ABR1,ABR2,harper2019partition,najnudel2018extreme}. 
See also the surveys \cite{Harper:Bourbaki,bailey2022maxima}.

Another line of work has considered the asymptotic distribution of the characteristic polynomial of various matrix ensembles outside the limiting spectral support, where it converges locally uniformly to the exponential of a random analytic function; consequently there are a.a.s.\ no outlier eigenvalues. This beautiful argument, inspired by \cite{BaZe}, was first carried out for matrices with iid entries of fixed standardized entry distribution in \cite{BCG}, and has since been applied to a range of models \cite{Coste23,CLZ,HaLo,FrGZ,CGZX}.
The argument was extended to sparse random matrices, first in \cite{Coste23} for the case of 
iid Bernoulli($d/n$) entries with $d>0$ fixed 
(the adjacency matrix of a sparse directed Erd\H os--R\'enyi graph), and then in \cite{CLZ} for the sum of $d$ iid uniform permutation matrices (the adjacency matrix of a sparse directed $d$-regular graph). Remarkably, the limiting analytic function for characteristic polynomials is the same in both cases.
Fran\c{c}ois~\cite{francois2026ewens} recently obtained a related extension for permutation matrices sampled from generalized Ewens measures, where the limiting function is again the exponential of a Poisson series, recovering the uniform permutation case and producing a weighted family of Poissonian limits.
In \cite{Coste23} Coste defines the \emph{Poissonian holomorphic chaos} (PHC) as a boundary limit of this function (see \cite[Definition 2.9]{Coste23}). 
The PHC is a random distribution, and a Poissonian analogue of the Gaussian holomorphic multiplicative chaos studied in \cite{najnudel2023secular,soundararajan2022model,najnudel2025fourier,atherfold2025fourier}. (Non-Gaussian holomorphic multiplicative chaoses of a different form have been studied in \cite{gu_universality_2025}.)
As for the Gaussian case, it is expected that regularity properties of the PMC are closely related to the distribution of large values of the characteristic polynomials on the circle at the edge of the limiting spectral support; see \cite[$\mathsection$8.4]{Coste23}. 
A further connection between the Poissonian and Gaussian settings was established in \cite{franccois2026characteristic}, who showed that the limiting random analytic functions associated with characteristic polynomials of self-normalized random matrices interpolate between Poissonian and Gaussian multiplicative chaos as the tail index of the entries before normalization varies from $0$ to $2$.

As we explain in $\mathsection$\ref{sec:ideas}, extremes for the Poissonian field \eqref{def:field} are well modeled by those of a branching random walk in time-inhomogeneous random environment (BRWtRE), and our proof makes use of ideas and results from recent works \cite{MaMi:branching,Kriechbaum:time} on this model. We defer further discussion of the  literature on BRWtRE to $\mathsection$\ref{sec:ideas}. 
A time-inhomogenous branching random walk also arose in the work \cite{CiLa} establishing the leading order of the maximum of the log-characteristic polynomial of a matrix with real iid entries. 
There the time-inhomogeneity is different from the quenched random inhomogeneity of the earlier work \cite{cook2020maximum} on random permutation matrices and of the present work: the sequence of branching laws and step sizes is deterministic and piecewise constant in time, with a single change at a macroscopic time.  The maximum for such a branching random walk was studied in \cite{FaZe}, and generalized to sequences of branching laws piecewise-constant in multiple macroscopic time windows in \cite{Mallein:piecewise}.

\subsection{Organization}

For \Cref{thm:main}, the high-probability upper bound on the maximum of $X_n$ is proved in $\mathsection$\ref{sec:upper}, and the matching lower bound is proved in $\mathsection$\ref{sec:lower}. 
The latter section is more involved, with several technical lemmas deferred to later sections. 
In $\mathsection$\ref{sec:decouple} we establish a crucial pathwise decoupling property for the field, and $\mathsection$\ref{sec:barrier} gathers various technical estimates for random walks bounded by a quenched random barrier, making use of results from \cite{MaMi:branching,Kriechbaum:time}.
In Appendix \ref{sec:covariance} we verify the logarithmic correlation structure \eqref{Cov-sim0} for $X_n$.

\subsection{Notation}

We write $C,C'$ etc.\ for finite positive constants depending only on the function $f$. Their values may change from line to line.
We write $a\ls b$ to mean $a=O(b)$.  
Dependence of the implicit constant in $O(\cdot)$ on a parameter $q$ is indicated with a subscript, e.g.  $a\ls_q b$ means $|a|\le C(q)b $.

The asymptotic notation $f=o(g)$ is with respect to the limit $n\to\infty$. Dependence of the rate of convergence of $f/g$ to zero on a parameter $q$ is indicated by writing $f=o_q(g)$. 
We will often implicitly assume that $n$ is larger than any fixed constant (possibly depending on $f$). 
We say that an event depending on $n$ holds \emph{asymptotically almost-surely (a.a.s.)} if it holds with probability tending to 1 as $n\to\infty$.

For a sequence $(z_j)$ we often abbreviate
\begin{equation}    \label{seq-diff}
    z_{i,j}:= z_j-z_i
\end{equation}
for $i\le j$. 

We let $(B_j)_{j\ge1}$ be a discrete-time random walk with independent standard Gaussian increments $B_{j+1}-B_j$. Thus, $(B_j)_{j\ge1}$ has the distribution of a Brownian motion evaluated at integer times. 

\subsubsection*{Tilted measures}

Define a titled probability measure $\P_*$ and expectation $\E_*$ by
\begin{equation}    \label{def:tiltE}
    \E_*(\,\cdot\,) = \E e^{\theta_* X_{2n}(0)-\lam(\theta_*)\cdot 2n}(\,\cdot\,)\,.
\end{equation}
Note the change of measure depends only on the sequence $\crot$.
(The choice of index $2n$ is arbitrary; we only want to ensure it is larger than any index we consider in the paper.)
We abbreviate
\begin{align}
    \E^\carr(\,\cdot\,) &:= \E( \,\cdot\,|\carr)\,, \label{def:ET}\\
    \E^\carr_*(\,\cdot\,) &:= \E_*(\,\cdot\,|\carr) = \E \big(  e^{\theta_* X_{2n}(0)-\lam(\theta_*)\cdot 2n}(\,\cdot\,) \,\mid\, \carr\big)    \label{def:tiltET}
\end{align}
and $\P^\carr(\,\cdot\,) = \E^\carr \ind(\,\cdot\,)$,
$\P^\carr_*(\,\cdot\,) = \E^\carr_* \ind(\,\cdot\,)$.
We define $\E^{\hcarr}$, $\E^{\hcarr}_*$ in the same way with $\hcarr$ in place of $\carr$, for the rounded arrivals process $\hcarr$ defined in \eqref{def:hcarr}.

\subsection{Statement on AI usage}
Artificial intelligence tools were used to proofread the prose and check for notational consistency. They played no role in the mathematical content of this work.

\section{Some ideas of the proof}
\label{sec:ideas}

Our proof adapts the barrier method for branching random walks in a time-inhomogeneous random environment to a field with no independent subtrees.
Two ideas are crucial: conditioning on rounded log-frequencies $(\tau_j)_{j\ge 1}$ provides the randomness needed for approximate independence, and retaining an environment-dependent correction in the centering of the maximum allows accurate comparisons of quenched barrier probabilities before obtaining the logarithmic shift.

As with many works on extreme values of LCFs, the above-mentioned works \cite{ABB2017maximum,paquette2018maximum,chhaibi2018maximum,paquette2022extremal,cook2020maximum} on Haar unitary and permutation matrices follow a paradigm introduced by Bramson in his seminal work on the extremal particle for branching Brownian motion \cite{bramson1978maximal}, where accurate upper and lower bounds are obtained from sharp estimates for the first and second moments of the number of particles attaining an extreme value whose path stays below a barrier. Insertion of the barrier is crucial to capture the correct constant for the logarithmic correction (the constant $\frac34$ in \eqref{FHK})  and to get concentration of particle counts. In the application to random fields $(F(t))_{t\in S}$, ``particles'' are associated to points $t$ in a discretization $S'$ for $S$, and the ``particle path'' for a point $t\in S'$ is a partial sum $(\sum_{i\le j}f_i(t))_{j=1}^n$, where $F(t)=\sum_{j=1}^n f_j(t)$ is a hierarchical decomposition of $F$. The index $j$ is thus viewed as a time variable, and levels $f_j$ in the decomposition are analogous to the displacements made by particles in generation $j$ of a branching random walk. 

In the present work, the fluctuations of the frequencies $\arr_j$ in \eqref{def:field} requires us to insert a \emph{random} barrier, which is a random walk with a drift. The square-root fluctuations of the barrier lead to the presence of the endpoint $\arr_n$ in the asymptotic \eqref{main-asymp} for the maximum, as well as the appearance of the function $\gamma$ of \eqref{def:gamma} in the definition of the constant $c_3$ for the logarithmic correction.  

For an analogy to branching processes, similarly to previous works we can view the sum \eqref{def:field} as a hierarchical decomposition of $X_n$, with the level-$j$ contribution $f(e^{\arr_j}t+\rot_j)$ corresponding to the displacements made by the $j$th generation. In the analogy, 
$e^{\arr_j}$ is the effective size of the population at generation $j$. Hence, the population typically grows like $\exp(n)$, but fluctuations in $\arr_j$ correspond to fluctuations in the reproduction rate of the population. We are thus led to consider a branching random walk model in which the law for the number of offspring is independently resampled at each generation. For a simple example, take a sequence $(\xi_j)_{j=1}^n$ of iid Bernoulli variables, and conditional on this sequence, for each $j$, everyone has $\ind(\xi_j=0) + 3\ind(\xi_j=1)$ children. 
Then, conditional on the random genealogical tree thus obtained, in each generation the distances of each child from its parent are iid standard Gaussian. 

The process just described is (a special case of) the branching random walk in time-inhomogeneous random environment (BRWtRE), which first appeared in \cite{biggins2004measure}. It is important to stress that the ``environment'' in this model is the random genealogical tree (which in the above examples were encoded by $(\xi_j)_{j=1}^n$). There is also work on branching random walks in \emph{space}-inhomogeneous random environment, where displacements depend on the spatial locations of the particles; see \cite{Kriechbaum:space} and references therein. 
We note that the BRWtRE is a natural model for a spreading population for which random occurrences such as long periods of drought affect the reproduction rate of an entire generation.

The leading order speed of the extremal particle for BRWtRE was established in \cite{huang2014branching}, and a refined asymptotic of the form \eqref{main-asymp} was proved in \cite{MaMi:branching}, with $\carr=(\arr_j)_{j\ge1}$ a random walk determined by the random sequence of branching laws. (We note that these works consider a more general class of BRWtRE than described above, in which the time environment is a random sequence of laws on point processes giving the displacements of children from their parents.) In  \cite{Kriechbaum:time} it was shown that conditional on the environment $\carr$, the law of the extremal particle is tight after appropriate recentering, with high probability in $\carr$.

For the field $X_n(s) = \sum_{j\le n} f(e^{\arr_j}s + \rot_j)$, 
the above analogy suggests viewing the arrivals process $(\arr_j)$ of frequencies as the quenched randomness analogous to the random genealogical tree, while the shifts $\rot_j$ provide additional randomness for jump displacements of the particles.
This is indeed the view taken in the proof of the upper bound on $X_n$ in $\mathsection$\ref{sec:upper}. 
However, note that while in a branching process, the steps taken by each individual in a fixed generation are independent, for the field $X_n$ there is only one random variable $\rot_j$ per generation $j$. This is not enough randomness to implement the second moment method for the lower bound in $\mathsection$\ref{sec:lower}. 

Instead, for the lower bound we take the  quenched randomness to be the \emph{rounded arrivals process} $(\wh\arr_j)$, where $\wh\arr_j= \rho\lf \arr_j/\rho\rf$ is obtained by rounding $\arr_j$ down to the nearest integer multiple of a parameter $\rho=\rho_n$ that tends to 0 slowly with $n$ ($\rho_n=\exp( - \sqrt{\log n})$ works). It turns out that the residual randomness $\arr_j-\wh\arr_j$ together with the rotation $\rot_j$ is enough to establish an approximate branching decorrelation structure for the quenched process. 

Our arguments make use of several estimates established in \cite{MaMi:barrier,Kriechbaum:time} for the probability that a random walk stays under a quenched random barrier, giving pointwise control on the field $X_n$. 
As is typically the case in the literature on LCFs, the argument to  control the maximum by the second moment method is relatively clean for branching models, where particle paths switch from perfect correlation to perfect independence at the branching point. 
Establishing an approximate version of this decoupling is a significant technical hurdle for the proof of \Cref{thm:main}. In $\mathsection$\ref{sec:decouple} we establish approximate decoupling for quenched barrier events for a pair of paths 
\[
(X_j(s)-X_m(s),X_j(t)-X_m(t))_{j=m}^n
\]
traced out by the evolving field at two separated points $s,t\in\R$, as soon as $|s-t|$ is large compared with the reciprocal $e^{-\arr_{m+1}}$ of the base frequency in the sum $X_j-X_m = \sum_{m<i\le j} f(e^{\arr_i}\cdot + \rot_i)$.
The decoupling bound for barrier events (see \Cref{lem:decouple}) introduces perturbations to the barriers that must then be removed with the help of fine stability estimates for barrier probabilities (see \Cref{lem:stab}).

For the barrier estimates, we combine the comparison arguments in \cite{Kriechbaum:time} with the asymptotic estimates in \cite{MaMi:barrier,MaMi:branching}, adapting them to our setting.
A key point is to retain the logarithm of a reference quenched barrier probability, denoted by $\beta_n$ in $\mathsection$\ref{sec:lower}, in the environment-dependent centering of the maximum throughout the moment estimates.
We compare the relevant barrier probabilities with this reference probability, and only after obtaining the quenched bounds do we use the logarithmic asymptotic of $\beta_n$ to identify the constant $c_3$.

\section{Upper bound}
\label{sec:upper}

In this section we prove the asymptotic in \Cref{thm:main} holds as an upper bound, i.e.\ for any fixed $\eps>0$, 
\begin{equation}
    \label{main-asymp-upper}
    \theta_*\max_{t\in [-1,1]} X_n(t) \le  K_n -(\phi_*-\eps)\log n
\end{equation}
with probability $1-o_\eps(1)$,
where here and in the sequel we denote
\begin{equation}
    \label{def:Kj}
    K_j:= \lam(\theta_*)j+\arr_j\,.
\end{equation}
For compactness of notation, in this section we set
\begin{equation}
    Y_j(s):= \theta_* X_j(s) - K_j\,.
\end{equation}
\eqref{main-asymp-upper} is an immediate consequence of the following:

\begin{prop}
    [Quenched upper bound]
    \label{prop:upper}
    For any fixed $\eps>0$,
    \[
    \Pcarr\bigg( \sup_{|s|\le 1} Y_n(s) \ge -(\phi_*-\eps)\log n)\bigg) \le n^{-\eps/8}
    \]
    with probability $1-o_\eps(1)$.
\end{prop}

The following allows us to control the maximum of $X_n$ using finite meshes of the domain.

\begin{lemma}
    [Continuity]
    \label{lem:cont}
    Let 
    \[
    H_j(\delta):=\sup_{s,t\in\R: |s-t|\le \delta}\max_{i\le j}|X_i(s)-X_i(t)|\,.
    \]
    There is a finite constant $C$ such that, with
    \[
    \cG_{\ref{lem:cont}} :=\bigcap_{\substack{i,j\in [n]:\\i-j\ge C\log n}}\{ \arr_i-\arr_j \ge (i-j)/2 \}\,,
    \]
    we have 
    \begin{equation}    \label{G1-Pbound}
        \P(\cG_{\ref{lem:cont}}) \ge 1-n^{-10}
    \end{equation}
    and
    \begin{equation}    \label{G1-Hbound}
        \ind(\cG_{\ref{lem:cont}}) \max_{j\le n} H_j(e^{-\arr_j}/\log n) \le C \qquad a.s.
    \end{equation}
\end{lemma}

\begin{proof}
    For any $i,j\in\N$ with $j-i\ge k$, since $\arr_{i+k}-\arr_i$ has the distribution of a sum of independent rate-1 exponential variables, 
    \[
    \P( \arr_j-\arr_i \le k/2) \le \P( \arr_{i+k}-\arr_i\le k/2) \le \exp( - ck)
    \]
    for a constant $c>0$. By the union bound, 
    \[
    \P((\cG_{\ref{lem:cont}})^c)\le 
    \P( \exists i,j\in[n]: j-i\ge k, \arr_j-\arr_i\le k/2) \le n^2 \exp( - ck).
    \]
    With $k=C\log n$ for a sufficiently large constant $C$, the right hand side is at most $n^{-10}$, giving \eqref{G1-Pbound}.

    On $\cG_{\ref{lem:cont}}$, for any $i\le n$ and $\delta>0$, 
    \[
    |X_i(s)-X_i(t)| \le \sum_{k \le i} |f(e^{\arr_k} s+\rot_k) - f(e^{\arr_k}t+\rot_k)| \ls \delta \sum_{k\le i}e^{\arr_k}
    \]
    so for any $j\le n$, 
    \begin{equation}    \label{bd:Hj}
        H_j(\delta) \ls \delta E_j
    \end{equation}
    where $E_j:= \sum_{k\le j}e^{\arr_k}$.
    On $\cG_{\ref{lem:cont}}$ we have
    \begin{align*}
        e^{-\arr_j}E_j 
        &=\sum_{k\le j} e^{-(\arr_j-\arr_k)}
        \ls \log n + \sum_{k\le j-C\log n} e^{-(j-k)/2} \ls \log n
    \end{align*}
    and combining with \eqref{bd:Hj} for $\delta=e^{-\arr_j}/\log n$ gives \eqref{G1-Hbound}.
\end{proof}

Recall the tilted probability measure $\P_*$ and expectation $\E_*$ in \eqref{def:tiltE}
\begin{align*}
    \E_*(\,\cdot\,) = \E e^{\theta_*X_{2n}(0)- \lam(\theta_*)\cdot 2n}(\,\cdot\,)\,.
\end{align*}
Note the change of measure depends only on the sequence $\crot$.
We have
\begin{equation}    \label{E*V*}
    \E_* X_j(0)=\speed j\,, \qquad \Var_*(X_j(0))=\lambda''(\theta_*)j\, \qquad \forall j\le 2n.
\end{equation}
Recall from \eqref{def:tiltET} the notations $\E^\carr(\,\cdot\,) := \E( \,\cdot\,|\carr)$ and
\begin{align*}
    \E^\carr_*(\,\cdot\,) &:= \E_*(\,\cdot\,|\carr) = \E \big(  e^{\theta_* X_{2n}(0) - \lam(\theta_*)\cdot 2n}(\,\cdot\,) \,\mid\, \carr\big)  \,.
\end{align*}

The following is an adaptation of \cite[Theorem 3.3, Equation (3.4)]{MaMi:branching} in our setting.

\begin{lemma}
    [Quenched barrier probability upper bound]
    \label{lem:quenched.upper}
    For any fixed $\eps>0$, the event $\cG_{\ref{lem:quenched.upper}}(\eps)$ that 
    \[
    \sup_{0\le x,y\le \log^3n} \Pcars\bigg(\max_{j\le n} Y_j(0)\le x\,,\; Y_n(0) \ge x-y\bigg) \le n^{-\phi_*+\eps}
    \]
    holds with probability $1-o_\eps(1)$.
\end{lemma}

\begin{proof}
    While \cite[Theorem 3.3, Equation (3.4)]{MaMi:branching} is phrased in terms of a spinal walk of a BRWtRE, their argument applies to a more general setting of a centered random walk staying below another quenched centered random walk and ending in a slowly-growing interval, under certain moment assumptions.
    In our setting the two walks are $(\theta_*(X_j(0)-\speed j))$ and $(K_j-\theta_*\speed j)=(\arr_j-j)$. Their increments are centered and have variances $\theta_*^2\lambda''(\theta_*)$ and $1$ under $\P_*$, respectively. Thus \Cref{lem:quenched.upper} corresponds to \cite[Theorem 3.3, Equation (3.4)]{MaMi:branching} with the choices $x_n=(\log n)^3$, $k=k'=0$, and $T_j=Y_j(0)$ and exponent $\phi_*$ (in place of their $\lambda$) as in \eqref{def:phi*}.

    We follow their argument towards \cite[Theorem 3.3, Equation (3.4)]{MaMi:branching}. That is, decomposing the barrier probability into thirds at times $\lfloor n/3\rfloor$ and $(n-\lfloor n/3\rfloor)$ by the Markov property under conditioning, and estimate them separately. 
    The decomposition and the estimates of the first and last third follow the same argument towards \cite[Theorem 3.3, Equation (3.4)]{MaMi:branching}, and we remark below the modifications to the proof of \cite[Lemma 3.6, Equation (3.10)]{MaMi:branching} in our setting, which is needed in estimating the (conditional) probability of the middle third.
    
    In our setting the quenched centered random walk $(K_j-\theta_*\speed j)$ under $\P_*$ has a weaker moment property than their global assumptions \cite[(1.7) and (1.8)]{MaMi:branching}.
    However, the assumption \cite[(1.7)]{MaMi:branching} is not used in the proof of \cite[Theorem 3.3]{MaMi:branching}.
    Moreover, they use the assumption \cite[(1.8)]{MaMi:branching} to establish an almost sure Berry--Esseen error of order $O(n^{-1/2})$ in \cite[Equation (3.12)]{MaMi:branching}, whereas an error of order $n^{-1/2+o_{\P}(1)}$ suffices in our \Cref{lem:quenched.upper}.
    
    Let $E_j:=\arr_j-\arr_{j-1}, j\ge 1$ be the iid Exp$(1)$ arrival gaps, so $Y_j(0)-Y_{j-1}(0)=\theta_*(f(\rot_j)-\speed)-(E_j-1)$.
    Let $\cG=\{\max_{j\le n}|E_j-1|^3\le C\log^3 n\}$.
    By a union bound and the exponential tail of $E_1$, there is a large constant $C>0$ such that $\lim_{n\to\infty}\P(\cG)=1$.
    On this event, we have
    \[
       \bone(\cG)\max_{j\le n}\frac{\E_{*}^{\carr}[|Y_j(0)-Y_{j-1}(0)|^3]}{\mathrm{Var}_*^{\carr}(Y_j(0)-Y_{j-1}(0))} = \bone(\cG)\max_{j\le n}\frac{\E_{*}^{\carr}[|Y_j(0)-Y_{j-1}(0)|^3]}{\theta_*^2\lambda''(\theta_*)} \lesssim_f \log^3 n=n^{o(1)}.
    \]
    So, on the event $\cG$, we can replace the $O(n^{-1/2})$ upper bound in \cite[Equation (3.12)]{MaMi:branching} by $O(\frac{\log^3 n}{\sqrt{n}})$. With the choice $x_n=\log^3 n$ in their proof, the same argument still gives \cite[Lemma 3.6, Equation (3.10)]{MaMi:branching}.
    The rest of the proof of \cite[Theorem 3.3, Equation (3.4)]{MaMi:branching} applies without any further modifications.
\end{proof}

\begin{proof}
    [Proof of \Cref{prop:upper}]
    For each $j$ let $\Sigma_j= [-1,1]\cap ( e^{-\arr_j}/\log n)\Z$.
    From \Cref{lem:cont}, the union bound and stationarity we have
    \begin{align}
        &\Pcarr \Big( \sup_{|s|\le 1,j\le n} Y_j(s) > \log^2n\Big)  \notag\\
        &\qquad\le n^{-10} + \ind(\cG_{\ref{lem:cont}})\Pcarr \Big( \sup_{|s|\le 1,j\le n} Y_j(s) > \log^2n\Big)  \notag\\
        &\qquad\le n^{-10}+\sum_{j\le n} |\Sigma_j|\Pcarr(Y_j(0)>\tfrac12\log^2n)  \notag\\
        &\qquad = n^{-10}+\sum_{j\le n} |\Sigma_j|\Ecars e^{\lam(\theta_*)j-\theta_*X_j(0)}\ind(Y_j(0)>\tfrac12\log^2n)  \notag\\
        &\qquad\ls n^{-10}+ \sum_{j\le n} e^{\arr_j}(\log n ) e^{\lam(\theta_*)j-K_j - \frac12\log^2 n} \Pcars(Y_j(0)>\tfrac12\log^2n)\\
        &\qquad\ls n^{-10}+n(\log n) e^{-\frac12\log^2n} = o(1) \label{upper-bprob}
    \end{align}
    where in the final line we bounded $\Pcars(Y_j(0)>\tfrac12\log^2n)\le1$.

    Now set
    \begin{align}
        \cE_s &:= \Big\{ \max_{j\le n} Y_j(s) \le 2\log^2n\,,\; Y_n(s) \ge -(\phi_*-\tfrac\eps 2)\log n\Big\}\,,\\
        \cG^2 &:= \Big\{\, \Pcars(\cE_0) \le n^{-\phi_*+\eps/4}\Big\}.
    \end{align}
    From \eqref{upper-bprob} and Lemmas \ref{lem:cont}, \ref{lem:quenched.upper}, to complete the proof it suffices to show
    \begin{equation}
        \label{upper-goal}
        \ind(\cG_{\ref{lem:cont}}\cap\cG^2) \Pcarr\bigg( \bigcup_{s\in \Sigma_n} \cE_s\bigg) \le n^{-\eps/8}
    \end{equation}
    for all $n$ sufficiently large. 
    For this we apply the union bound, stationarity and a change of measure:
    \begin{align*}
        \ind(\cG_{\ref{lem:cont}}\cap\cG^2) \Pcarr\bigg( \bigcup_{s\in \Sigma_n} \cE_s\bigg) 
        &\le \ind(\cG^2) \sum_{s\in \Sigma_n} \Pcarr(\cE_s) \\
        &= \ind(\cG^2) |\Sigma_n| \Pcarr(\cE_0)\\
        &=\ind(\cG^2) |\Sigma_n| \Ecars e^{\lam(\theta_*)n-\theta_*X_n(0)}\ind(\cE_0)\\        
        &\ls \ind(\cG^2) (\log n)e^{\arr_n+\lam(\theta_*)n-K_n+(\phi_*-\frac\eps2)\log n}\Pcars(\cE_0)\\
        &\le (\log n) n^{\phi_*-\eps/2} n^{-\phi_*+\eps/4} 
        =(\log n)n^{-\eps/4}
    \end{align*}
    giving \eqref{upper-goal} to complete the proof.
\end{proof}

\section{Lower bound}
\label{sec:lower}

In this section we prove the asymptotic in \Cref{thm:main} holds as a lower bound, 
that is, that for any fixed $\eps>0$, 
\begin{equation}
    \label{lower-goal0}
    \P\,\Big(\theta_* \sup_{|s|\le 1} X_n(s) \ge K_n - (\phi_*+\eps) \log n\Big) = 1-o_\eps(1).
\end{equation}
recalling that $K_n:=\lam(\theta_*)n+\arr_n$. 

As in the proof of the upper bound, we  deduce \eqref{lower-goal0} from a quenched lower bound; however, rather than quenching the precise arrival times $\arr_j$, in this case we need a little more randomness left over in order to establish a branching correlation structure. We take the random environment to be a coarse-graining of the arrivals sequence $\carr=(\arr_j)$. 

To that end, let $\rho=\rho_n$ satisfy
\begin{equation}    \label{assu:rho}
    \omega(\log\log n) \le \log (1/\rho) \le o(\log n).
\end{equation}
Define rounded arrival times $\hcarr=(\wh\arr_j)_{j\ge1}$ by 
\begin{equation}    \label{def:hcarr}
    \wh\arr_j:= \rho \lf \arr_j/\rho\rf. 
\end{equation}
Conditional on $\hcarr$, form a new sequence $\wt\carr=(\wt\arr_j)_{j\ge1}$ by independently and uniformly permuting the labels of the arrival times $\carr_k:=\{\arr_j: \wh\arr_j=\rho k\}$ for each $k\ge0$.
(Note that the size of each $\carr_k$ is Poisson-distributed of rate $\rho=o(1)$, so the randomization of labels is trivial for most $k$.)
Conditional on $\hcarr$, the variables
\begin{equation}
    \label{tarr-harr}
    \{\wt\arr_j-\wh\arr_j\}_{j\ge1}
\end{equation}
are independent and uniform in $[0,\rho]$. 
We then apply the same permutations (within each bin $\{j:\wh\arr_j=\rho k\}$) to the phases to form a new sequence $(\wt\rot_j)_{j\ge 1}$.
Let
\begin{equation}    \label{def:tX}
    \wt X_j(s):= \sum_{i\le j} f(e^{\wt\arr_i}s + \wt\rot_i)\,.
\end{equation}
We reorder the arrival times in this way so that we can apply a decoupling lemma later (see \Cref{lem:decouple}).
Note that for all $j$ and $s$,
\begin{equation}
    \label{tX-at0}
    \wt X_j(s)\,\overset{law}{=}\, X_j(s)\,.
\end{equation}
However, $(X_j(s))_{j\ge 1}$ and $(\wt X_j(s))_{j\ge 1}$ need not be equal pathwise.
The following shows the difference is negligible for our purpose:

\begin{lemma}
    \label{lem:XtX}
    With probability $1-n^{-\omega(1)}$, 
    \[
    \max_{j\le n} \|X_j-\wt X_j\|_\infty \ls \frac{\log n}{\log(1/\rho)} = o(\log n/\log\log n)\,.
    \]
\end{lemma}

\begin{proof}
    The second bound follows from \eqref{assu:rho}. For the first bound, 
    recall $N_k=|\carr\cap[k\rho,(k+1)\rho)|\sim\text{Poisson}(\rho)$ independently.
        Let $M_N=\max_{0\le k\le N}N_k$ for any $N\ge 0$.
        Since for any integer $j$ of the form $j=\sum_{k\le \ell}N_k$ for some $\ell$,
        \[
        \sum_{i\le j}f\left(e^{\arr_i}t+\rot_i\right) = \sum_{i\le j}f\left(e^{\wt{\arr}_i}t+\wt\rot_i\right)
        \]
        we have
        \begin{align*}
            &\P\left(\max_{j\le n}\|X_j-\wt X_j\|_{\infty}\ge C\log n/\log(1/\rho)\right)\\
            &\le \P\left(M_{\floor{\arr_n/\rho}}\ge C\log n/\log(1/\rho)\right) \\
            &\le \P(\arr_n\ge 2n)+\P(M_{\floor{2n/\rho}}\ge C\log n/\log(1/\rho))\\
            &\le \P(\arr_n\ge 2n)+\frac{2n}{\rho}\P\left(N_1\ge C\log n/\log(1/\rho)\right).
        \end{align*}
        Since $\arr_n$ is the sum of $n$ iid Exp$(1)$ random variables and $N_1$ is a Poisson$(\rho)$ random variable, standard tail estimates give the desired probability upper bound after taking $C>2$ to be a large constant.
\end{proof}

For a sequence $(w_j)_{j\ge1}$ denote 
\begin{equation}\label{beta-w}
      \beta_n(w) := \P\Big( \, \max_{j\le n} B_j -w_j/\sigma_* \le 10\,,\; B_n-w_n/\sigma_*\in [-1,0]\Big) 
\end{equation}
where we recall $\sigma_*=\theta_*\sqrt{\lam''(\theta_*)}$.
(Note that $(\sigma_*B_j)_{j\ge1}$ has the same covariance as the random walk $(\sigma_*X_j(0))_{j\ge1}$.)
In particular, if $W$ is a random process, then $\beta_n(W)$ denotes the conditional probability obtained by evaluating this functional at the realized path of $W$.
Set
\begin{equation}    \label{def:n0}
    n_0 := \lf C_0 \log(1/\rho)\rf
\end{equation}
for a constant $C_0$ to be taken sufficiently large (depending only on $f$). 
Recalling the notation \eqref{seq-diff}, denote
\begin{align}
    V_j(s) &:= \wt X_{n_0,n_0+j}(s)\\
    \wh K_j &:= \wh\arr_{n_0,n_0+j} + \lam(\theta_*)j\\
    \wh W_j &:= \wh\arr_{n_0,n_0+j} - j\,.
\end{align}

\begin{prop}
    [Quenched lower bound]
    \label{prop:lower}
    Assuming the constant $C_0$ in \eqref{def:n0} is sufficiently large, there is a constant $c>0$ such that 
    \[
    \Phcarr \Big( \theta_*\sup_{|s|\le 1} V_n(s) \ge \wh K_n + \log \beta_n(\wh W)\Big) \ge 1-\rho^c
    \]
    with probability $1-o(1)$.
\end{prop}

Before proceeding to the proof of \Cref{prop:lower} we combine it with \Cref{lem:XtX} 
and the following
to establish \eqref{lower-goal0}.

\begin{lemma}
    [Quenched Gaussian barrier probability]
    \label{lem:gauss-lower}
    For any fixed $\eps>0$, the event 
    \[
    \cG_{\ref{lem:gauss-lower}}(\eps)=\Big\{\beta_n(\wh W)\ge n^{-\phi_*-\eps}\Big\}
    \]
    holds with probability $1-o_\eps(1)$.
\end{lemma}

\begin{proof}
    This is a consequence of \Cref{prop:LB-wiggle}. Indeed, define $W_j=\tau_{n_0,n_0+j}-j$, then $W$ is a centered random walk satisfying the moment assumptions in \Cref{prop:LB-wiggle}, and $\max_{j\le n}|\wh W_j-W_j|\le \rho<\frac1{40}$ a.s.~for sufficiently large $n$. Since the residuals $\{(\tau_j-\wh\tau_j):j\ge 1\}$ are independent of the barrier event in \eqref{beta-w} when conditioned on $\hcarr$, we have 
    \[
    \beta_n(\wh W) = \P\Big(\max_{j\le n} B_j-\wh W_j/\sigma_* \le 10,\, B_n-\wh W_n/\sigma_* \in[-1,0]\,\Big|\, W\Big).
    \]
    Then the lower bound follows from \Cref{prop:LB-wiggle}.
\end{proof}

\begin{proof}[Proof of \eqref{lower-goal0}]
Since $n+n_0=(1+o(1))n$, it suffices to show \eqref{lower-goal0} with $X_{n_0+n}, K_{n_0+n}$ in place of $X_n,K_n$. 
From \Cref{lem:XtX} we can replace $X_{n_0+n}$ with $\wt X_{n_0+n}$, and since
\[
\|\wt X_{n_0}\|_\infty \ls n_0 \ls \log(1/\rho) = o(\log n)\qquad a.s.
\]
we may further replace $\wt X_{n_0+n}$ with $V_n$. 
Finally, since $\E \arr_{n_0} = n_0$, from Markov's inequality we get 
\[
K_{n_0} = \arr_{n_0} + \lam(\theta_*)n_0=O(n_0)= o(\log n)\qquad a.a.s.
\]
so
\[
K_{n_0+n} = \wh K_n + K_{n_0} + O(\rho) = \wh K_n + o(\log n) \qquad a.a.s.
\]
Hence, it suffices to show
\begin{equation}
    \label{lower-goal1}
    \P\Big(\theta_* \sup_{|s|\le 1} V_n(s) \ge \wh K_n - (\phi_*+\eps) \log n\Big) \ge 1-o_\eps(1).
\end{equation} 
But this is immediate
from \Cref{prop:lower} and \Cref{lem:gauss-lower}.
\end{proof}

In the proof of \Cref{prop:lower} we will condition on $\cG_{\ref{lem:gauss-lower}}$ and several other ``good'' events, denoted by $\cG$ with various subscripts, all of which will be $\sigma(\hcarr)$-measurable and hold a.a.s.

The next such event ensures that path segments  $(V_j(s))_{j=\ell+1}^n, (V_j(t))_{j=\ell+1}^n$ are approximately independent once $|s-t|$ is large compared to $\exp( -\wh\arr_{n_0+\ell})$, which is roughly the reciprocal of the largest frequency in the sum for $V_\ell=\wt X_{n_0,n_0+\ell}$ (see \eqref{def:tX}). 
This path-wise decoupling is precisely the form of branching de-correlation structure that we need for the field $V_n$.

\begin{lemma}
    [Path decoupling]
    \label{lem:decouple}
    There exists a constant $C<\infty$ depending only on $f$ such that the following holds for all $n\ge C$. 
    Let 
    \[
    \cG_{\ref{lem:decouple}}:=\bigg\{ \bigwedge_{1\le i<j\le n} \Big( j-i\ge \log^2n\;\Rightarrow\; \wh\arr_j-\wh\arr_i \ge (j-i)/3 \Big) \bigg\}.
    \]
    Note this event holds a.a.s.\ by \Cref{lem:cont} and the fact that $\wh\arr_i=\arr_i+O(\rho)$.
    Let $\ell\ge 0$ and let $s,t\in \R$ with $|s-t|\ge \rho^{-C}e^{-\wh\arr_{n_0+\ell+1}}$. 
    For each $m\in [\ell+1,n]$ let $J_{m,1}, J_{m,2}$ be arbitrary compact intervals. 
    Then on $\cG_{\ref{lem:decouple}}$ we have
    \begin{align*}
        &\Phcarr \bigg( \bigwedge_{\ell<m\le n} V_{\ell,m}(s)\in J_{m,1}\,,\; V_{\ell, m}(t) \in J_{m,2}\bigg) \\
        &\qquad\le (1+\rho^{10}) 
        \Phcarr \bigg( \bigwedge_{\ell<m\le n} V_{\ell,m}(0)\in J_{m,1}^\rho\bigg) 
         \Phcarr\bigg( \bigwedge_{\ell<m\le n} V_{\ell,m}(0)\in J_{m,2}^\rho\bigg) 
    \end{align*}
    where $J_{m,i}^\rho:=J_{m,i}+[-\rho,\rho]$ denotes the closed $\rho$-neighborhood of interval $J_{m,i}$.
\end{lemma}

The proof of \Cref{lem:decouple} is deferred to $\mathsection$\ref{sec:decouple}.

\begin{proof}[Proof of \Cref{prop:lower}]

For brevity, let
\begin{equation}    \label{def:an}
    a_n:=\log(1/\rho)\,.
\end{equation}
Then
$a_n=\Theta(n_0)$, $\omega(\log\log n) \le a_n \le o(\log n)$.
Set
\begin{align}
    \wh L_j &:= \min(j, n-j)^{1/6} - \frac jn\log   \beta_n(\wh W)\,, \quad j\le n\\
    Z_j(s) &:= \theta_* V_j(s) - \wh K_j + \wh L_j \\
    \cE_s(y_-,y_+) &:= \Big\{ \max_{j\le n} Z_j(s) \le y_+ \,,\; Z_n(s) \in [y_-,y_+]\Big\} \label{def:Ey-pm}\\
    \cE_s(y) &:= \cE_s(-y,a_n+y)\\
    \cE_s &:= \cE_s(0) 
    = \Big\{ \max_{j\le n} Z_j(s) \le a_n \,,\; Z_n(s) \in [0,a_n]\Big\}\,.
\end{align}

\begin{lemma}[Gaussian comparison]
    \label{lem:compare}
    For any $\delta>0$, with probability at least $1-\delta$ the following holds for all sufficiently large $n$: 
    \begin{equation}
        \Phcars(\cE_0(y_-,y_+)) = e^{O_\delta(\log\log n)}\beta_n(\wh W) 
    \end{equation}
    uniformly for $y_+\in[a_n/2,\log^2n]$, $y_-\in y_+-[a_n/2,\log n]$.
    In particular, by \eqref{assu:rho}, the event $\cG_{\ref{lem:compare}}$ that for all $y\in[a_n/2,\log^2n]$ and all $z\in[a_n/2,\log n]$,
    \begin{equation}
        \Phcars(\cE_0(y-z,y)) = e^{O(\log(1/\rho))}\beta_n(\wh W) 
    \end{equation}
    holds a.a.s.
\end{lemma}

\begin{proof}
    This follows from \Cref{lem:Kriech3} (with $r_n=a_n/4$) and \Cref{lem:gauss-lower}. 
\end{proof}

\begin{lemma}
    [Stability]
    \label{lem:stab}
    There exists $c>0$ such that the event  $\cG_{\ref{lem:stab}}$ that $\Phcarr(\cE_0(\rho))\le (1+\rho^c)\Phcarr(\cE_0)$ holds a.a.s.
\end{lemma}

\begin{proof}
    This follows from \Cref{lem:stab-gen} (with $x=y=r_n=a_n$, so that $\eps=\rho$) and \Cref{lem:gauss-lower}.
\end{proof}

Let
\begin{align}
    \Sigma_n &:= [-1,1] \cap e^{-\wh\arr_{n_0+n}}\Z 
    \,,\qquad 
    \cN_n :=\sum_{s\in \Sigma_n}\ind_{\cE_s}\,.
\end{align}
For the conditional first moment, by stationarity we have
\begin{align}   \label{1mom-on3}
    \E^{\hcarr}\cN_n
    &= 
    |\Sigma_n| \Phcarr(\cE_0) \asymp  e^{\wh\arr_{n_0+n}}\Phcarr(\cE_0) 
\end{align}
For general $y\ge0$ we note that
\begin{align*}
    \E_*^{\hcarr}e^{\lam(\theta_*)n-\theta_*V_n(0)}\ind(\cE_0(y))
    &= \E^{\hcarr} e^{\theta_*(X_{2n}(0)-V_n(0))-\lam(\theta_*)n}\ind(\cE_0(y))\\
    &= \E^{\hcarr} e^{\theta_*(X_{n_0}(0) + X_{n_0+n,2n}(0)) - \lam(\theta_*)n}\ind(\cE_0(y))\\
    &= \P^{\hcarr}(\cE_0(y))\E^{\hcarr} e^{\theta_*(X_{n_0}(0) + X_{n_0+n,2n}(0)) - \lam(\theta_*)n} \\
    &=\P^{\hcarr}(\cE_0(y))
\end{align*}
where in the third line we used that the event $\cE_0$ is independent of $(\rot_j)_{j\notin[n_0+1,n_0+n]}$, and for the final line:
\begin{align*}
    \E^{\hcarr} e^{\theta_*(X_{n_0}(0) + X_{n_0+n,2n}(0)) - \lam(\theta_*)n} 
    &= \E \exp\bigg( -\lam(\theta_*)n +\sum_{j\in [n_0]\cup[n_0+n+1,2n]}\theta_* f(\rot_j) \bigg) = 1\,.
\end{align*}
Hence,
\begin{align}
    \Phcarr(\cE_0(y))
    &= 
    \E_*^{\hcarr}e^{\lam(\theta_*)n-\theta_*V_n(0)}\ind_{\cE_0(y)} \notag\\
    &=
    e^{\lam(\theta_*)n-\wh K_n+\wh L_n +O(a_n+y)}\P^{\hcarr}_*(\cE_0(y)) \notag\\
    &= 
    e^{-\wh\arr_{n_0,n_0+n}+O(a_n+y)} \P^{\hcarr}_*(\cE_0(y))/\beta_n(\wh W) \,. \notag
\end{align}
From \Cref{lem:compare} we hence have
\begin{align}
    \ind_{\cG_{\ref{lem:compare}}} 
    \Phcarr(\cE_0(y))
    &= 
    \ind_{\cG_{\ref{lem:compare}}} e^{-\wh\arr_{n_0,n_0+n}+O(y+ \log(1/\rho))} \,.\label{probE00}
\end{align}
From this and \eqref{1mom-on3} we get the first moment estimate:
\begin{equation}    \label{lower:1mom}
    \ind_{\cG_{\ref{lem:compare}}}\E^{\hcarr}\cN_n =
    \ind_{\cG_{\ref{lem:compare}}}
    e^{\wh\arr_{n_0}+ O(\log(1/\rho))}.
\end{equation}
(Note this gives both a lower and an upper bound.)

For the conditional second moment, we decompose the sum as
\begin{align}   \label{2mom-123}
    \Ehcarr \cN_n^2 = \I + \sum_{\ell=1}^n \II_\ell + \III 
\end{align}
where
\begin{align*}
    \I &= \sum_{\substack{s,t\in\Sigma_n:\\ |s-t|<\rho^{-C}e^{-\wh\arr_{n_0+n}}}} \Phcarr(\cE_s\cap\cE_t)\\
    \II_\ell &= \sum_{\substack{s,t\in\Sigma_n:\\ \rho^C|s-t|\in[e^{-\wh\arr_{n_0+\ell}},e^{-\wh\arr_{n_0+\ell-1}})}} \Phcarr(\cE_s\cap\cE_t)\\
    \III &= \sum_{\substack{s,t\in\Sigma_n:\\ |s-t|\ge\rho^{-C}e^{-\wh\arr_{n_0}}}} \Phcarr(\cE_s\cap\cE_t)\,.
\end{align*}
We expect the main contribution to $\Ehcarr\cN_n^2$ to come from the sum $\III$ over pairs of well-separated points, which in the branching process analogy correspond to pairs of particles that split early on. 

\begin{claim}[Late-branching pairs]
    \label{claim:I}
    On $\cG_{\ref{lem:compare}}$, $\I \le e^{\wh\arr_{n_0}+O(\log (1/\rho))}$.
\end{claim}

\begin{claim}[Early-branching pairs]
    \label{claim:III}
    On $\cG_{\ref{lem:decouple}}\cap\cG_{\ref{lem:stab}}$, 
    $\III \le (1+O(\rho^c))(\E \cN_n)^2$.
\end{claim}

\begin{claim}[Intermediate-branching pairs]
    \label{claim:II}
    $\sum_{\ell=1}^n\II_\ell \le e^{\wh\arr_{n_0}+O(\log (1/\rho))}$ a.a.s.
\end{claim}

We now conclude the proof of \Cref{prop:lower} on the claims. 
We may restrict to the event that the bounds in Claims \ref{claim:I}, \ref{claim:III}, \ref{claim:II} hold.
(Note this event is $\sigma(\hcarr)$-measurable.)
From these bounds and \eqref{2mom-123} we have
\begin{equation}    \label{2mom-1}
    \Ehcarr \cN_n^2
    \le (1+O(\rho^c)) (\Ehcarr \cN_n)^2 + e^{\wh\arr_{n_0} + O(\log(1/\rho))}.
\end{equation}
Squaring \eqref{lower:1mom} gives
\[
e^{\wh\arr_{n_0}} \le e^{-\wh\arr_{n_0}+O(\log(1/\rho))} (\Ehcarr \cN_n)^2
\]
which, combined with \eqref{2mom-1}, yields
\begin{equation}    \label{2mom-2}
    \Ehcarr \cN_n^2
    \le \Big(1+O(\rho^c) + e^{-\wh\arr_{n_0} + O(\log(1/\rho))}\Big) (\Ehcarr \cN_n)^2.
\end{equation}
Since $\arr_{n_0}=\wh\arr_{n_0}+O(\rho)$ is a sum of $n_0$ iid rate-1 exponential variables, we have
\[
\wh\arr_{n_0}  \ge n_0/2 = \frac12\lf C_0\log(1/\rho)\rf \qquad a.a.s.
\]
and restricting to this event we have
\[
e^{-\wh\arr_{n_0} + O(\log(1/\rho))} \le e^{-\frac12C_0\log(1/\rho) + O(\log(1/\rho))} \le e^{-\frac14 C_0\log (1/\rho)} = \rho^{C_0/4}
\]
if $C_0$ is sufficiently large. 
Returning to \eqref{2mom-2}, we obtain
\begin{equation}    \label{2mom-3}
    \Ehcarr \cN_n^2
    \le \big(1+O(\rho^c) \big) (\Ehcarr \cN_n)^2.
\end{equation}
By Cauchy--Schwarz,
\[
\Phcarr(\cN_n\ge 1) \ge \frac{(\Ehcarr \cN_n)^2}{\Ehcarr\cN_n^2} \ge 1-O(\rho^c).
\]
It only remains to note that on the event that $\cN_n\ge1$, there is a point $s\in\Sigma_n$ such that
\[
0 \le Z_n(s) = \theta_*V_n(s) - \wh K_n + \wh L_n
= \theta_*V_n(s) - \wh K_n - \log\beta_n(\wh W)
\]
which concludes the proof of \Cref{prop:lower}. 
\end{proof}

In the remaining subsections we prove Claims \ref{claim:I}, \ref{claim:III} and \ref{claim:II}.

\subsection{Late-branching pairs (Proof of \Cref{claim:I})}

For the contribution of pairs $(s,t)$ of nearby points, we drop the intersection with $\cE_t$ to bound:
\begin{align}
    \I &\le \sum_{\substack{s,t\in\Sigma_n:\\ |s-t|<\rho^{-C}e^{-\wh\arr_{n_0+n}}}} \Phcarr(\cE_s)
    \asymp \rho^{-C} e^{-\wh \arr_{n_0+n}}|\Sigma_n|^2 \Phcarr(\cE_0)
    \asymp \rho^{-C} \Ehcarr \cN_n\,.
\end{align}
The claim follows from this and \eqref{lower:1mom}.

\subsection{Early-branching pairs (Proof of \Cref{claim:III})}

Turning to bound $\III$, consider $\Phcarr(\cE_s\cap\cE_t)$ for a fixed pair $s,t\in\R$ with $|s-t|\ge\rho^{-C}e^{-\wh\arr_{n_0}}$.
By stationarity of the field we may take $t=0$ and $s$ positive.
Note that since $f$ is bounded, the event $\cE_s$ takes the form $\wedge_{n_0<m\le n_0+n}V_m(s)\in J_m$ for compact intervals $J_m$. We may hence apply \Cref{lem:decouple},   taking $C$ sufficiently large, to bound
\begin{align}   \label{III-decouple}
    \Phcarr(\cE_0\cap\cE_s) \ind_{\cG_{\ref{lem:decouple}}}
    &\le  (1+\rho^{10})\Phcarr(\cE_0(\rho))^2 \ind_{\cG_{\ref{lem:decouple}}}\,.
\end{align}
Combining \Cref{lem:stab} with \eqref{III-decouple} we conclude
\begin{align}   \label{III-Pbd}
\Phcarr(\cE_0\cap\cE_s)\ind_{\cG_{\ref{lem:decouple}}\cap\cG_{\ref{lem:stab}}}
    &\le  (1+O(\rho^c))\Phcarr(\cE_0)^2 \ind_{\cG_{\ref{lem:decouple}}\cap\cG_{\ref{lem:stab}}}
\end{align}
and hence
\begin{align}
    \label{III-bd}
\III \ind_{\cG_{\ref{lem:decouple}}\cap\cG_{\ref{lem:stab}}}&\le 
     (1+O(\rho^c)) |\Sigma_n|^2 \Phcarr(\cE_0)^2 \ind_{\cG_{\ref{lem:decouple}}\cap\cG_{\ref{lem:stab}}}
     \,.
\end{align}
The claim follows from this and \eqref{1mom-on3}. 

\subsection{Intermediate-branching pairs (Proof of \Cref{claim:II})}

Turning to the terms $\II_\ell$, fix arbitrary $\ell\in[n]$ and $s\in [\rho^{-C} e^{-\wh\arr_{n_0+\ell}}, \rho^{-C}e^{-\wh\arr_{n_0+\ell-1}})$.
The idea (by now standard in the literature on maxima of LCFs) is to exploit the branching correlation structure (as provided by \Cref{lem:decouple}) by splitting $\cE_0\cap\cE_s$ into independent events depending on increments of $(Z_j(0),Z_j(s))_{j=1}^n$ before and after $j=\ell$.  
We decompose
\begin{align}   \label{Es-split-q}
    \cE_s 
    &= \bigcup_{q\le a_n} \cE_s\cap\{ Z_\ell(s)\in (q-1,q]\}
\end{align}
with the union running over integer $q$. (By the boundedness of increments of the field these events are nonempty for only finitely many $q$.)
With
\begin{equation}\label{def:Eq-le}
    \cE_{s,q}^\le=\cE_{s,q}^{\le \ell} := \Big\{ \max_{j\le \ell} Z_j(s) \le a_n, \; Z_\ell(s)\in (q-1,q]\Big\}
\end{equation}
(we will usually omit $\ell$ from the superscript to lighten notation)
we have
\begin{align*}
    \cE_s\cap\{Z_\ell(s)\in (q-1,q]\}
    &= \cE_{s,q}^\le \cap \Big\{Z_\ell(s)\in (q-1,q]\,,\; \max_{\ell<j\le n} Z_j(s) \le a_n , \;Z_n(s)\in [0,a_n]\Big\}\,.
\end{align*}
On the second event on the right hand side above, recalling the notation \eqref{seq-diff}, for $j>\ell$ we have
\[
Z_j(s) = Z_\ell(s) + Z_{\ell,j}(s) \in Z_{\ell,j}(s) + (q-1,q]
\]
so
\begin{align*}
    \cE_s\cap\{Z_\ell(s)\in (q-1,q]\} 
    \subseteq \cE_{s,q}^\le \cap \cE_{s,q}^>(1)
\end{align*}
where
\begin{equation}
    \label{def:Eq->}
    \cE_{s,q}^>(y)= \cE_{s,q}^{>\ell}(y):= \Big\{ \max_{\ell<j\le n} Z_{\ell,j}(s) \le a_n+y -q\,,\; Z_{\ell, n}(s) \in -q+ [-y, a_n+y] \Big\}\,.
\end{equation}
Combining with \eqref{Es-split-q} we have shown
\begin{equation}\label{Es-split-q2}
    \cE_s \subseteq \bigcup_{q\le a_n} 
    \cE_{s,q}^\le \cap \cE_{s,q}^>(1)\,.
\end{equation}
The point is that after fixing $q$, the events $\cE_{s,q}^\le$,  $\cE_{s,q}^>(1)$ are conditionally independent given $\hcarr$.
Thus, with the summation running over integers $q, r$,
\begin{align*}
    \Phcarr(\cE_0\cap\cE_s) 
     &\le \sum_{q,r\le a_n} \Phcarr\big(\cE_{0,q}^\le\cap \cE_{s,r}^\le \cap \cE_{0,q}^>(1)\cap \cE_{s,r}^>(1)\big) \\
     &=\sum_{q,r\le a_n} \Phcarr\big(\cE_{0,q}^\le\cap\cE_{s,r}^\le\big) \Phcarr\big(\cE_{0,q}^>(1)\cap \cE_{s,r}^>(1)\big).
\end{align*}
Now since $s\ge \rho^{-C} e^{-\wh\arr_{n_0+\ell}}$, from \Cref{lem:decouple} and stationarity, we have that on $\cG_{\ref{lem:decouple}}$,
\[
\Phcarr\big(\cE_{0,q}^>(1)\cap \cE_{s,r}^>(1)\big) 
\le 2 \Phcarr\big(\cE_{0,q}^>(2)\big)\Phcarr\big(\cE_{0,r}^>(2)\big)\,,
\]
say. (Here we do not need as much precision as in \eqref{III-decouple}.)
Apply the inequality $ab\le\frac12(a^2+b^2)$, we have that on $\cG_{\ref{lem:decouple}}$,
\begin{align*}
    \Phcarr(\cE_0\cap\cE_s)
    &\le 2\sum_{q,r\le a_n}\Phcarr\big(\cE_{0,q}^\le\cap\cE_{s,r}^\le\big)\Phcarr\big(\cE_{0,q}^>(2)\big)\Phcarr\big(\cE_{0,r}^>(2)\big)\\
    &\le \sum_{q,r\le a_n}\Phcarr\big(\cE_{0,q}^\le\cap\cE_{s,r}^\le\big)\left(\Phcarr\big(\cE_{0,q}^>(2)\big)^2+\Phcarr\big(\cE_{0,r}^>(2)\big)^2\right).
\end{align*}
Since the events $\cE_{s,r}^\le$ are disjoint as $r$ varies, we have
\[
    \sum_{r\le a_n}\Phcarr\big(\cE_{0,q}^\le\cap\cE_{s,r}^\le\big)\le\Phcarr\big(\cE_{0,q}^\le\big).
\]
Similarly, summing over $q$ gives an upper bound of $\Phcarr(\cE_{s,r}^\le)$.
By stationarity, $\Phcarr(\cE_{s,r}^\le)=\Phcarr(\cE_{0,r}^\le)$.
Thus, on $\cG_{\ref{lem:decouple}}$,
\begin{align}
    \Phcarr(\cE_0\cap\cE_s) 
    &\le 2\sum_{q\le a_n} 
    \Phcarr\big(\cE_{0,q}^\le\big)\Phcarr\big(\cE_{0,q}^>(2)\big)^2\,.
    \label{E0Es-bd2}
\end{align}
Again by conditional independence, we can recombine:
\[
\Phcarr\big(\cE_{0,q}^\le\big)\Phcarr\big(\cE_{0,q}^>(2)\big)
= \Phcarr\big(\cE_{0,q}^\le\cap\cE_{0,q}^>(2)\big)
\]
and note that
\[
\cE_{0,q}^\le\cap\cE_{0,q}^>(2)
\subseteq \big\{ Z_\ell(0)\in (q-1,q]\big\} \cap \cE_0(3).
\]
Returning to \eqref{E0Es-bd2}, we have that on $\cG_{\ref{lem:decouple}}$,
\begin{equation}    \label{E0Es-bd3}
    \Phcarr(\cE_0\cap\cE_s) 
    \le 2\sum_{q\le a_n} 
    \Phcarr\big( \cE_0(3)\,,\; Z_\ell(0)\in (q-1,q]\big) \Phcarr(\cE_{0,q}^>(2))\,.
\end{equation}

We turn to estimate $\Phcarr(\cE_{0,q}^>(2))$. 
We consider $\cE_{0,q}^>(y)$ for general $y\ge0$. 
Recalling the change of measure defined under \eqref{def:tiltET}, we have
\begin{align*}
    \Phcarr(\cE_{0,q}^>(y))
    &= \E^{\hcarr}_* e^{\lam(\theta_*)(n-\ell) - \theta_*V_{\ell, n}(0)}\ind( \cE_{0,q}^>(y))\,.
\end{align*}
Recall that on $\cE_{0,q}^>(y)$ we have
\[
Z_{\ell,n}(0) 
= \theta_*V_{\ell, n}(0) - \wh K_{\ell, n} + \wh L_{\ell, n}
\in -q+[-y,a_n+y]
\]
where
\begin{align*}
    \wh K_{\ell, n} & = \wh\arr_{n_0+\ell, n_0+n} + \lam(\theta_*)(n-\ell) \\
    \wh L_{\ell, n}&= -\min(\ell, n-\ell)^{1/6} - (1-\frac\ell n)\log \beta_n(\wh W)\,.
\end{align*}
Thus, on $\cE_{0,q}^>(y)$ we have
\[
\theta_*V_{\ell, n}(0) \in \wh\arr_{n_0+\ell, n_0+n} + \lam(\theta_*)(n-\ell) + \min(\ell, n-\ell)^{1/6} +(1-\frac\ell n)\log \beta_n(\wh W) -q + [-y,a_n+y]
\]
so
\begin{equation*}
    \Phcarr(\cE_{0,q}^>(y))
    = e^{O(a_n+y)} \exp\Big( q - \wh\arr_{n_0+\ell, n_0+n} - 
    \min( \ell, n-\ell)^{1/6} - (1-\frac\ell n) \log\beta_n(\wh W)
    \Big) \Phcars( \cE_{0,q}^>(y))\,.
\end{equation*}
Towards having a clean cancellation when we multiply by the cardinality of pairs of points $(s,t)$ in the sum defining $\II_\ell$, it will be convenient to replace $\wh\arr_{n_0+\ell, n_0+n}$ with $\wh\arr_{n_0+\ell-1, n_0+n}$.
Since 
$\arr_{n_0+\ell}- \arr_{n_0+\ell-1}$ is exponential with rate 1, and is a.s.\ within $O(\rho)$ of $\wh\arr_{n_0+\ell}- \wh\arr_{n_0+\ell-1}$, one easily sees that the event $\cG_{\ref{arr-steps}}$ that
\begin{equation}
    \label{arr-steps}
    \wh\arr_{n_0+\ell}- \wh\arr_{n_0+\ell-1}\le \max\Big(\log\log n, \frac1{10}\min(\ell,n-\ell)^{1/6}\Big) \qquad \forall 1\le \ell \le n
\end{equation}
holds a.a.s.
Thus, on $\cG_{\ref{arr-steps}}$,
\begin{align}
    \label{PE0q-bd1}
    \Phcarr(\cE_{0,q}^>(y))
    = e^{O(a_n+y+\log\log n)} \exp\Big( q - \wh\arr_{n_0+\ell-1, n_0+n} - 
    R(\ell,n)
    \Big) \Phcars( \cE_{0,q}^>(y))
\end{align}
where
\begin{equation}
    \label{def:Rln}
    R(\ell, n):= \frac12\min( \ell, n-\ell)^{1/6} + (1-\frac\ell n) \log\beta_n(\wh W)\,.
\end{equation}

We combine \eqref{PE0q-bd1} with \eqref{E0Es-bd3} to estimate $\II_\ell$. Since the number of pairs $s,t\in\Sigma_n$ with $\rho^{C}|s-t|\in [e^{-\wh\arr_{n_0+\ell}}, e^{-\wh\arr_{n_0+\ell-1}})$ is at most 
\[
|\Sigma_n|^2 e^{-\wh\arr_{n_0+\ell-1}}\asymp\rho^{-C} \exp( 2\wh\arr_{n_0+n} - \wh\arr_{n_0+\ell-1})
\]
by stationarity of the field we have that on $\cG_{\ref{lem:decouple}}\cap\cG_{\ref{arr-steps}}$,
\begin{align*}
    \II_\ell 
    &\ls |\Sigma_n|^2 e^{-\wh\arr_{n_0+\ell-1}} 
    \sum_{q\le a_n} 
    \Phcarr\Big( \cE_0(3)\,,\; Z_\ell(0)\in (q-1,q]\Big) \Phcarr\big(\cE_{0,q}^>(2)\big)\\
    &\le \exp\Big( O(\log(1/\rho))+ 2\wh\arr_{n_0+n} - \wh\arr_{n_0+\ell-1} -\wh\arr_{n_0+\ell-1, n_0+n} - R(\ell,n)\Big) \\
    &\qquad\qquad\qquad\times \sum_{q\le a_n} \Phcarr\Big( \cE_0(3)\,,\; Z_\ell(0)\in (q-1,q]\Big)
    \,e^q\,\Phcars\big( \cE_{0,q}^>(2)\big)\\
    &= \exp\Big( O(\log(1/\rho))+ \wh\arr_{n_0+n}  - R(\ell,n)\Big) \\
    &\qquad\qquad\qquad\times \sum_{q\le a_n} \Phcarr\Big( \cE_0(3)\,,\; Z_\ell(0)\in (q-1,q]\Big)
    \,e^q\,\Phcars\big( \cE_{0,q}^>(2)\big)\,.
\end{align*}
(Note the error $O(y+\log\log n)$ in \eqref{PE0q-bd1} is negligible (for $y=2$) compared to $\log(1/\rho)$ by \eqref{assu:rho}, and $a_n=\log(1/\rho)$.)
Now we bound
\begin{align*}
    &\sum_{q\le a_n} \Phcarr\big( \cE_0(3)\,,\; Z_\ell(0)\in (q-1,q]\big)
    e^q\,\Phcars( \cE_{0,q}^>(y))\\
    &\quad \le 
    \sum_{q\le a_n} \Phcarr\big( \cE_0(3)\,,\; Z_\ell(0)\in (q-1,q]\big)
    \sup_{z\in [-a_n,+\infty)}e^{-z}\,\Phcars(\cE_0^>(2,z))\\
    &\quad = \Phcarr\big( \cE_0(3)\big)\sup_{z\in [-a_n,+\infty)}e^{-z}\,\Phcars(\cE_0^>(2,z))
\end{align*}
where
\begin{equation}
    \cE_0^>(y,z) := \Big\{ \max_{\ell<j\le n} Z_{\ell, j}(0) \le a_n+z+y\,,\; Z_{\ell, n}(0) \in z+ [-y,a_n+y]\Big\}\,.
\end{equation}
Hence, on $\cG_{\ref{lem:decouple}}\cap\cG_{\ref{arr-steps}}$,
\begin{align}
    \II_\ell &\le 
    e^{ \wh\arr_{n_0+n} + O(\log(1/\rho))  }
    \Phcarr\big( \cE_0(3)\big)
    e^{- R(\ell,n)}
    \sup_{z\in [-a_n,+\infty)}e^{-z}\,\Phcars(\cE_0^>(2,z))  \,.\label{II-bd2}
\end{align}
We will bound the right hand side differently for two ranges of $\ell$, namely 
\begin{align*}
    &\text{large $\ell$:}\qquad \log^7n< \ell\le n\\
    &\text{small $\ell$:} \qquad 1\le \ell \le \log^{7}n\,.
\end{align*}

\subsubsection{Late-intermediate pairs}

For large $\ell$, on the event $\cG_{\ref{lem:gauss-lower}}(1)$ we can bound
\[
R(\ell, n) \ge \frac12\min( \ell, n-\ell)^{1/6} - \frac{n-\ell}{n}(\phi_*+1)\log n
\ge \frac14 \min( \ell, n-\ell)^{1/6} .
\]

Taking the trivial bound
\[
\sup_{z\in [-a_n,+\infty)}e^{-z}\,\Phcars(\cE_0^>(2,z))\le e^{a_n}=e^{\log(1/\rho)}
\]
we conclude 
\begin{align*}
    \sum_{\log^7n\le \ell\le n }\II_\ell 
    &\le e^{ \wh\arr_{n_0+n} + O(\log(1/\rho))  }\Phcarr\big( \cE_0(3)\big) \sum_{\log^7n\le \ell\le n} e^{
    - \frac14\min(\ell, n-\ell)^{1/6}}\\
    &\le e^{ \wh\arr_{n_0+n} + O(\log(1/\rho))  }\Phcarr\big( \cE_0(3)\big)\,.
\end{align*}
Now recall from \eqref{probE00} that
\[
\ind_{\cG_{\ref{lem:compare}}}
\Phcarr\big( \cE_0(3)\big)\asymp
\ind_{\cG_{\ref{lem:compare}}} e^{-\wh\arr_{n_0,n_0+n}+O(\log(1/\rho))} 
\]
and hence, on $\cG_{\ref{lem:decouple}}\cap\cG_{\ref{lem:compare}}\cap\cG_{\ref{lem:gauss-lower}}(1)\cap\cG_{\ref{arr-steps}}$, 
\begin{equation}    \label{large-ell}
    \sum_{\log^7n\le \ell\le n }\II_\ell 
    \le e^{\wh\arr_{n_0} + O(\log(1/\rho)) }.
\end{equation}

\subsubsection{Early-intermediate pairs}

Turning to bound the contribution of small $\ell$, in this case we bound the final expression in \eqref{II-bd2} as follows:
\begin{equation}
    \sup_{z\in [-a_n,+\infty)}e^{-z}\,\Phcars(\cE_0^>(2,z))\le e^{-\log^{3/2}n} + e^{a_n}
    \sup_{z\in [-a_n,\log^{3/2}n]} \Phcars(\cE_0^>(2,z))\,.
\end{equation}
Consider the event $\cE_0^>(2,z)$ for an arbitrary $z\in [-a_n,\log^{3/2}n]$. 
Let 
\begin{equation}\label{def:Eq-le-long-ending}
    \cE_{s}^{\le}(L)=\cE_{s}^{\le\ell}(L):=\Big\{\max_{j\le\ell}Z_j(s)\le a_n,\, Z_{\ell}(s)-a_n\in[-L,0]\Big\}.
\end{equation}
Note that $\cE_0^>(2,z)$ is conditionally independent of $\cE_{0}^\le(2a_n)$ given $\hcarr$. 
Thus,
\begin{equation}
    \label{mid-combine-E}
    \Phcars(\cE_{0}^\le(2a_n))\Phcars ( \cE_0^>(2,z)) = \Phcars(\cE_{0}^\le(2a_n) \cap \cE_0^>(2,z))\,.
\end{equation}
Now for any $y\ge0$ we have
\begin{align*}
    \cE_{0}^\le(2a_n) \cap \cE_0^>(y,z)
    &= \Big\{ \max_{j\le \ell} Z_j(0)\le a_n\,,\; 
    Z_\ell(0)\in [-a_n,a_n]\,,\;\\
    &\qquad\qquad 
    \max_{\ell<j\le n} Z_{\ell,j}(0)\le a_n+y+z\,,\;
    Z_{\ell,n}(0)\in z+[-y-a_n,a_n+y]\Big\}\\
    &\subseteq\Big\{ \max_{j\le n} Z_j(0)\le 2a_n+y+z\,,\;
    Z_n(0)\in z+ [-y-2a_n,2a_n+y]\Big\}\\
    &= \cE_0(z-y-2a_n,z+y+2a_n)
\end{align*}
recalling the notation \eqref{def:Ey-pm}.
From \Cref{lem:compare}, we have
\begin{equation*}
    \Phcars(\cE_0(z-y-2a_n,z+y+2a_n)) \ind_{\cG_{\ref{lem:compare}}} = e^{O(\log(1/\rho))} \beta_n(\wh W) \ind_{\cG_{\ref{lem:compare}}} 
\end{equation*}
uniformly for $y\le \frac13\log n$ and $-a_n\le z\le \log^{3/2}n$. 
Thus,
\begin{equation}
    \label{mid-combine-E2}
    \Phcars(\cE_{0}^\le(2a_n))\Phcars ( \cE_0^>(2,z)) \ind_{\cG_{\ref{lem:compare}}} \le
    e^{O(\log(1/\rho))}\beta_n(\wh W)
    \ind_{\cG_{\ref{lem:compare}}}
\end{equation}
for all $z\in [-a_n,\log^{3/2}n]$.

We first assume the following technical lemma to prove \Cref{claim:II}, and postpone its proof afterwards.

\begin{lemma}
    \label{lem:Kr-small}
    There is a constant $C>0$ such that the event
    \[
    \Phcars(\cE_{0}^{\le\ell}(2a_n)) \ge (\log n)^{-C}
    \]
    for all $1\le \ell\le \log^7 n$ holds a.a.s.
\end{lemma}

For small $\ell$ we also have
\begin{equation}    \label{small-l-R}
    e^{-R(\ell,n)} = \bigg(\frac1{\beta_n(\wh W)}\bigg)^{1-\frac\ell n} e^{-\frac12\ell^{1/6}}
    \le \frac1{\beta_n(\wh W)} e^{-\frac12\ell^{1/6}}\,.
\end{equation}

Hence, inserting the bounds \eqref{mid-combine-E2} and \eqref{small-l-R} in \eqref{II-bd2} and summing over $\ell$, on $\cG_{\ref{lem:decouple}}\cap\cG_{\ref{lem:compare}}\cap\cG_{\ref{arr-steps}}\cap\cG_{\ref{lem:Kr-small}}$ we have
\begin{align*}
    &\sum_{\ell \le \log^7 n} \II_\ell\\
    &\le e^{ \wh\arr_{n_0+n} + O(\log(1/\rho))  }
    \Phcarr\big( \cE_0(3)\big)
    \sum_{\ell\le \log^7n}e^{-R(\ell,n)}\sup_{z\in [-a_n,+\infty)}e^{-z}\,\Phcars(\cE_0^>(2,z))  \\
    &\le  e^{ \wh\arr_{n_0+n} + O(\log(1/\rho))  }
    \Phcarr\big( \cE_0(3)\big) 
    \sum_{\ell\le \log^7n}e^{-\frac12\ell^{1/6}}
    \bigg( e^{-\log^{3/2}n} + \frac1{\Phcars(\cE_{0}^{\le\ell}(2a_n))}\bigg)\\
    &\le  e^{ \wh\arr_{n_0+n} + O(\log(1/\rho))  }
    \Phcarr\big( \cE_0(3)\big) 
    \bigg(1+ 
    (\log n)^C\sum_{\ell\le \log^7n}e^{-\frac12\ell^{1/6}}\bigg)\,.
\end{align*}

Since $\sum_{\ell\ge 1}e^{-\frac12\ell^{1/6}}<\infty$, on $\cG_{\ref{lem:decouple}}\cap\cG_{\ref{lem:compare}}\cap\cG_{\ref{arr-steps}}\cap\cG_{\ref{lem:Kr-small}}$, 
\begin{equation}
    \sum_{\ell \le \log^7 n} \II_\ell
    \le e^{ \wh\arr_{n_0+n} + O(\log(1/\rho))  }
    \Phcarr\big( \cE_0(3)\big) 
    =e^{ \wh\arr_{n_0} + O(\log(1/\rho))  }
\end{equation}
where we used \eqref{probE00} in the last equality.
Together with \eqref{large-ell} this yields \Cref{claim:II}.    \qed\\

We turn to prove \Cref{lem:Kr-small}, which is a version of \cite[Lemma 11.6]{Kriechbaum:time} for our setting. \cite[Lemma 11.6]{Kriechbaum:time} does not directly apply, as the proof there assumes that increments of the quenched barrier random walk are sub-Gaussian, whereas the increments in our model are only sub-exponential. To deal with this we use the KMT coupling, recalled in the following lemma, to replace the random walk and random barrier with Gaussian random walks, at which point the argument of \cite{Kriechbaum:time} applies.

\begin{lemma}[KMT coupling, \cite{komlos1975approximation,komlos1976approximation}]
\label{lem:KMT}
    Let $(S_j)_{j\ge 1}$ be a centered random walk with iid increments such that $\E[e^{\eta |S_1|}]<\infty$ for some $\eta>0$. Denote $\sigma^2:=\Var(S_1)$.
    Then there exist constants $c_0,c_1,c_2>0$, depending only on the law of $S_1$, such that for every integer $N\ge 1$, one can couple a copy of $(S_j)_{j\le N}$ with a standard Brownian motion $(B_t)_{t\ge 0}$, so that for $G_t:=\sigma B_t$ and for any $\lambda\ge0$,
    \begin{equation}\label{eqn:KMT}
        \P\left(\max_{j\le N}|S_j-G_j|>c_0\log N+\lambda\right)\le c_1e^{-c_2 \lambda}.
    \end{equation}
\end{lemma}

\begin{proof}[Proof of \Cref{lem:Kr-small}]
    Set $W_j:=\tau_{n_0,n_0+j}-j$.
    Recall $\sup_{j\ge 1}|\wh W_j-W_j|\le \rho$. 
    On $\cG_{\ref{lem:gauss-lower}}$, for $j\le\log^7 n$, we have
    \[
    \wh L_j = j^{1/6}-\frac{j}{n}\log\beta_n(\wh W) = j^{1/6}+o(n^{-0.9}).
    \]
    Let $S_j:=\theta_*(V_j(0)-\speed j)$, for all $j\ge 1$. By definition $S_j$ has iid centered increments under $\P_*$, and $Z_j(0)=S_j-\wh W_j+\wh L_j$.
    Recall that under $\P_*$, the process $S$ is independent of $\carr$, and $\hcarr$ is measurable with respect to $\carr$.
    Thus, for an event depending only on $S$ and $\wh W$, its conditional probability given $\carr$ equals its conditional probability given $\hcarr$.
    Therefore, on $\cG_{\ref{lem:gauss-lower}}$ we have
    \begin{align}
        \Phcars(\cE_{0}^{\le\ell}(2a_n)) &= \Phcars\left(\max_{j\le\ell}S_j-\wh W_j+\wh L_j\le a_n,\,S_{\ell}-\wh W_\ell+\wh L_{\ell}\in[-a_n,a_n]\right) \notag\\
        &\ge \Phcars\left(\max_{j\le\ell}S_j-\wh W_j+j^{1/6}\le a_n/2,\,S_{\ell}-\wh W_\ell+\ell^{1/6}\in[-a_n/2,a_n/2]\right) \notag\\
        &= \Pcars\left(\max_{j\le\ell}S_j-\wh W_j+j^{1/6}\le a_n/2,\,S_{\ell}-\wh W_\ell+\ell^{1/6}\in[-a_n/2,a_n/2]\right) \notag\\
        &\ge \Pcars\left(\max_{j\le\ell}S_j-W_j+j^{1/6}\le a_n/3,\,S_{\ell}- W_\ell+\ell^{1/6}\in[-a_n/3,a_n/3]\right). \label{eqn:4.11-preparation}
    \end{align}

    Let $\sigma_S:=\Var_*(S_1)^{1/2}$ and note $\Var_*(W_1) = \Var(W_1) = 1$.
    With $N=\lceil\log^7 n\rceil$ and $\lambda=C_{\mathrm{KMT}}\log N$ in \Cref{lem:KMT}, there are two-dimensional standard Brownian motion $(B'_t, B_t)_{t\ge 0}$ and constants $c_i, i\in\{0,1,2\}$ (depending on laws of $S_1$ and $W_1$ under $\P_*$), such that $(S,B')$ is independent of $(\carr, B)$ under $\P_*$, and for $G_t:=\sigma_S B'_t$ we have
    \[
    \P_*\left(\max_{j\le\log^7 n}|S_j-G_j|\vee \max_{j\le\log^7 n}|W_j-B_j| \le (c_0+C_{\mathrm{KMT}})\log N\right) \ge 1-c_1N^{-c_2C_{\mathrm{KMT}}}.
    \]
    Let $\P_*^B(\cdot):=\P_*(\cdot\mid B)$.
    On the event $\cG_{W}:=\{\max_{j\le\log^7 n}|W_j-B_j| \le (c_0+C_{\mathrm{KMT}})\log N\}$, since $\log N=O(\log\log n)=o(a_n)$, for sufficiently large $n$ we have
    \begin{align}
        &\Pcars\left(\max_{j\le\ell}S_j-W_j+j^{1/6}\le a_n/3,\,S_{\ell}- W_\ell+\ell^{1/6}\in[-a_n/3,a_n/3]\right) \notag\\
        \ge &\P_*^B\left(\max_{j\le\ell}S_j-B_j+j^{1/6}\le a_n/4,\,S_{\ell}- B_\ell+\ell^{1/6}\in[-a_n/4,a_n/4]\right) \notag\\
        \ge &\P_*^B\left(\max_{j\le\ell}G_j-B_j+j^{1/6}\le a_n/5,\,G_{\ell}- B_\ell+\ell^{1/6}\in[-a_n/5,a_n/5]\right) - c_1 N^{-c_2C_{\mathrm{KMT}}} \notag\\
        \ge &\P_*^B\left(\max_{j\le\ell}G_j-B_j+j^{1/6}\le \sigma_S,\,G_{\ell}- B_\ell+\ell^{1/6}\in[0,\sigma_S]\right) - c_1 N^{-c_2C_{\mathrm{KMT}}}. \label{eqn:4.11-afterKMT}
    \end{align}

    Let $c_{\ell}:=\ell^{1/6}/\ell=\ell^{-5/6}$, and $h_{\ell}(j):=j^{1/6}-c_{\ell}j$ for $0\le j\le\ell$ and $\ell\ge 1$, and interpolate linearly in between integer times to define $h_{\ell}$ on $[0,\ell]$. Note $h_\ell(0)=h_\ell(\ell)=0$, and the family $(h_{\ell})_{\ell\ge 1}$ is nice in the sense of \cite[Definition 11.5]{Kriechbaum:time}.
    Define
    \[
    \frac{d \P_{\ell}}{d \P_*} := \exp\left(-\frac{c_{\ell}}{\sigma_S^2} G_\ell-\frac{c_{\ell}^2}{2\sigma_S^2}\ell\right),\qquad \ell\ge 1.
    \]
    Then $(G_t^{(\ell)}:=G_t+c_{\ell}t)_{t\in[0,\ell]}$ is a Brownian motion with variance $\sigma_S^2$ under $\P_{\ell}$. Hence we have
    \begin{align}
        &\P_*^B\left(\max_{j\le\ell}G_j-B_j+j^{1/6}\le \sigma_S,\,G_{\ell}- B_\ell+\ell^{1/6}\in[0,\sigma_S]\right) \notag\\
        &\ge e^{-\sigma_S^{-2}\left(\frac{|B_{\ell}|}{\ell^{5/6}}+\frac{\ell^{-2/3}}{2}\right)}\P_{\ell}^B\left(\max_{j\le\ell}G_j^{(\ell)}-B_j+h_{\ell}(j)\le \sigma_S,\,G^{(\ell)}_{\ell}- B_\ell\in[0,\sigma_S]\right). \label{eqn:4.11-tilted}
    \end{align}

    Under $\P_\ell$, the process $G^{(\ell)}/\sigma_S$ is a standard Brownian motion independent of $B$, and $(B_j)_{j\ge 1}$ has iid increments with sub-Gaussian tails.
    Apply \cite[Lemma 11.6]{Kriechbaum:time} to the linearly interpolated environment $(B_j/\sigma_S)_{j\ge0}$, the fixed family $(h_\ell/\sigma_S)_{\ell\ge1}$, and signed starting and endpoint parameters $y=-1$ and $y_0=0$.
    Since the continuous-time barrier event is contained in the discrete-time event above, there is a deterministic constant $C_0>0$, depending only on $\sigma_S$, such that almost surely the latter probability is at least $\ell^{-C_0}$ for all $\ell\ge M$, where $M\in\N$ is a $\sigma(B_j:j\ge0)$-measurable random variable.
    For each remaining integer $\ell$, this probability is strictly positive.
    Consequently, there is a $\sigma(B_j:j\ge0)$-measurable random variable $D>0$ a.s.~such that
    \begin{equation}\label{eqn:4.11-asLB}
        \P_{\ell}^B\left(\max_{j\le\ell}G_j^{(\ell)}-B_j+h_{\ell}(j)\le \sigma_S,\,G^{(\ell)}_{\ell}-B_\ell\in[0,\sigma_S]\right) \ge D\ell^{-C_0},\quad \ell\ge 1, \quad a.s.
    \end{equation}
    Here $D$ can be chosen as a fixed measurable function of the Gaussian environment path, so its distribution does not depend on $n$ or on the choice of the KMT coupling parameters.

    Therefore, on the event
    \[
    \cG:=\cG_{\ref{lem:gauss-lower}}\cap\cG_{W}\cap\left\{\sup_{\ell\ge 1}\frac{|B_{\ell}|}{\ell^{5/6}}\le\log N\right\}\cap\{D>N^{-1}\},
    \]
    combining \eqref{eqn:4.11-preparation}, \eqref{eqn:4.11-afterKMT}, \eqref{eqn:4.11-tilted}, and \eqref{eqn:4.11-asLB} gives, when $n$ is sufficiently large,
    \[
    \Phcars(\cE_{0}^{\le\ell}(2a_n)) \ge N^{-C_1}-c_1 N^{-c_2C_{\mathrm{KMT}}}
    \]
    for some $C_1<\infty$ depending on $\sigma_S$. 
    Note $\sup_{\ell\ge 1}\frac{|B_{\ell}|}{\ell^{5/6}}<\infty$ a.s.~by the law of iterated logarithm, we have
    \[
    \lim_{n\to\infty}\P\left(\sup_{\ell\ge 1}\frac{|B_{\ell}|}{\ell^{5/6}}\le\log N,\; D>N^{-1}\right) = 1.
    \]
    Take $C_{\mathrm{KMT}}$ sufficiently large finishes the proof.
\end{proof}

The proof of \Cref{prop:lower} is complete.

\section{Path decoupling}
\label{sec:decouple}

In this section we prove \Cref{lem:decouple}. 
We deduce it from a more general result:

\begin{prop}
    [Path decoupling]
    \label{prop:decouple}
    For any $A,\lam\ge1$, and $K>0$ there exist $C,\alpha<\infty$ depending only on $A,\lam,K$ and $f$ such that the following holds.
    Let $d\in \N$ and $\delta\in (0,c)$ for $c$ sufficiently small depending only of $f$, 
    and let $[a_k,b_k)$, $k=1,\dots, d$ be a  sequence of (not necessarily disjoint) intervals of length at least $\delta$, such that
    \begin{equation}    \label{decouple-growth}
        a_k - a_1\ge \alpha \log k\qquad \forall\; k\ge \delta^{-K}.
    \end{equation}
    Let $T_1,\dots,T_d$, $U_1,\dots, U_d$ be jointly independent random variables with $T_k$ uniform in $[a_k,b_k)$ and $U_k$ uniform in $[0,1)$, and for $u\in \R, m\in[d]$ set 
    \[
    F^u_m:= \sum_{k=1}^m f( e^{T_k}u + U_k)\,.
    \]
    For any fixed $s,t\in \R$ with $|s-t|\ge \delta^{-C}e^{-a_1}$ and $2d$ intervals $J_k^s,J_k^t$, $k=1,\dots, d$, we have
    \begin{align}
        &\P \bigg( \bigwedge_{1\le m\le d} F_m^s\in J_m^s\,,\;
        F_m^t \in J_m^t 
        \bigg) \notag\\
        & \le \big(1+ O(\delta)^{A\lam})\big)\,
        \P \bigg(\bigwedge_{1\le m\le d} F_m^s\in J_m^s(\delta') 
        \bigg)\,
        \P \bigg( \bigwedge_{1\le m\le d} F_m^t \in J_m^t(\delta') 
        \bigg) \notag
    \end{align}
    for $\delta'=O(\delta^\lam)$, 
    where $J(\delta')$ is the $\delta'$-neighborhood of $J$.
\end{prop}

\begin{proof}[Proof of \Cref{lem:decouple}]
     It is easily verified that on  the event $\cG_{\ref{lem:decouple}}$
     the condition \eqref{decouple-growth} 
     holds for the sequence $a_k = \wh\arr_{n_0+m+k}$, with $\delta=\rho$ and $K=1$, say. Recalling the observation \eqref{tarr-harr} that $T_k:=\wt\arr_{n_0+m+k}-\wh\arr_{n_0+m+k}$ are iid and uniform in $[0,\rho)$,   the claim then follows from \Cref{prop:decouple}. 
\end{proof}

We begin with a lemma for asymptotic independence of $(f(e^{\arr} s+\rot),f(e^{\arr}t+\rot))$ with $\arr,\rot$ independent and uniformly distributed, and $\{s,t\}$ well separated.

\begin{lemma}
    \label{lem:block-indep}
    For any $A>0$ there exists $B>0$ depending only on $A$ and $f$ such that the following holds.
    Let $0<a<b$, let $\arr$ be uniform in $[a,b)$ and let $\rot$ be uniform in $[0,1)$ and independent of $\arr$. Let $\eps\in(0,\frac12]$, let 
    $I_1,I_2$ be intervals contained in the range of $f$, with $|I_1|,|I_2|\ge \eps$,
    and let $s,t\in\R$ with $|s-t|\ge \eps^{-B}(b-a)^{-1}e^{-a}$. Then
    \[
    \P\big( f( e^\arr s+\rot)\in I_1, \,f(e^\arr t+\rot)\in I_2\big) 
    \le (1+ O(\eps^A)) \P( f(\rot)\in I_1) \P( f(\rot)\in I_2).
    \]
\end{lemma}

\begin{proof}
    Since $f$ is 1-periodic, without loss of generality we may take $s=0$ and $t\ge \eps^{-B}(b-a)^{-1}e^{-a}$.
    We abbreviate $(X,Y):= (f(\rot),f(e^\arr t+\rot))$ and let $\mu_t$ be the joint distribution of $(X,Y)$. 
    Then $X,Y\eqd f(\rot)$ are identically distributed, and we denote their common law
    by $\nu$.

    Let $\delta\in(0,\eps)$ to be chosen later, and let $h_1,h_2$ be smooth functions satisfying 
    \[
    1_{I_k}\le h_k\le 1_{I_k^\delta} \; \text{ pointwise}\,,\quad \|h_k'\|_{L^\infty}\le 2/\delta
    \]
    where $I_k^\delta:= I_k+(-\delta,\delta)$.
    We have
    \begin{align}
        \mu_t(I_1\times I_2)
        &\le \mu_t( h_1\otimes h_2) \notag\\
        &\le \nu(h_1)\nu(h_2) + |\mu_t(h_1\otimes h_2)-\nu(h_1)\nu(h_2)|    \notag\\
        &\le \nu(I_1)\nu(I_2) + \nu(I_1) (\nu(h_2)-\nu(I_2)) +\nu(h_2)(\nu(h_1)-\nu(I_1))   \notag\\
        &\qquad+ |\mu_t(h_1\otimes h_2)-\nu(h_1)\nu(h_2)|\,.    \label{indep1}
    \end{align}
    Here we abbreviate $\mu(h):=\int hd\mu$. 
    By \Cref{assu:f} we have
    \begin{equation}
       |I|/C\le \nu(I) \le C|I|^{1/C}\label{eqn:small-ball-estimate}
    \end{equation}
    for any interval $I$ contained in the range of $f$, and some $C>0$ depending only on $f$. It follows that
    \[
    \nu(h_k)-\nu(I_k) \le \nu(I_k^\delta\setminus I_k)
    \lesssim \delta^{1/C} \lesssim \delta^{1/C} \eps^{-1}\nu(I_k)
    \]
    for $k=1,2$. Hence, assuming $\delta<\eps^C$,
    \begin{align}
        \mu_t(I_1\times I_2)
        &\le \nu(I_1)\nu(I_2) ( 1+ O(\delta^{1/C}/\eps)) 
        + |\mu_t(h_1\otimes h_2)-\nu(h_1)\nu(h_2)|\,.   \label{indep2}
    \end{align}
    Turning to bound the final term on the right hand side above, denote the Fourier transforms of $h_1,h_2$ by
    \[
    \cF h_k(\xi):=\int_\R f(x) e^{-i\xi x}dx
    \]
    and the (joint) characteristic functions of $X$ and $(X,Y)$ by 
    \begin{align*}
        \phi_X(\xi)&:=\E e^{i\xi X} = \int_\R e^{i\xi x}d\nu(x)\\
        \phi_{X,Y}(\xi,\eta) &:= \E e^{i\xi X+i\eta Y} = \int_\R e^{i\xi x + i\eta y}d\mu_t(x,y).
    \end{align*}
    By Fourier inversion,
    \begin{align}
        \int_{\R^2} h_1\otimes h_2 d(\mu_t -\nu\otimes\nu)
        &= \frac1{4\pi^2} \int_{\R^2} \cF h_1(\xi) \cF h_2(\eta) \big(\phi_{X,Y}(\xi,\eta) - \phi_X(\xi)\phi_X(\eta) \big) d\xi d\eta \notag\\
        &\lesssim \|\cF h_1\|_{L^1} \|\cF h_2\|_{L^1} \|\phi_{X,Y}-\phi_X\otimes \phi_X\|_{L^\infty}\,.  \label{indep3}
    \end{align}
    From integration by parts and our assumptions on $h_1,h_2$ we have a pointwise bound
    \[
    |\cF h_k(\xi)| \le |\xi|^{-2} \|h_k''\|_{L^1} \lesssim |\xi|^{-2} /\delta\,.
    \]
    Combining with the trivial bound $\|\cF h_k\|_{L^\infty} \le \|h_k\|_{L^1}\lesssim \eps$ and integrating over $\xi$, we find
    \begin{equation}
        \label{indep4}
            \|\cF h_k\|_{L^1} \lesssim \sqrt{\eps/\delta}.
    \end{equation}
    It only remains to get a uniform estimate on the difference of characteristic functions. Fix arbitrary $\xi,\eta\in\R$ and denote the 1-periodic function
    \[
    g:\R\to\C\,,\quad g(x):= \E_\rot\exp( i\xi f(\rot) + i\eta f(x+\rot)).
    \]
    Note that for $\rot'\eqd \rot$ and independent from $\rot$, we have $\E g(\rot') = \phi_X(\xi)\phi_X(\eta)$. 
    Then with $g_0:= g-\E g(\rot')$ we have
    \begin{align}
        \phi_{X,Y}(\xi,\eta) - \phi_X(\xi)\phi_X(\eta)
        &= \E_\arr  g_0(e^\arr  t).
    \end{align}
    Since $g_0$ is smooth it has an absolutely convergent Fourier series expansion
    \[
    g_0(x) = \sum_{k\in\Z\setminus\{0\}} \hat g(k) e(kx)
    \]
    where $\hat g_0(0)=0$ since $g_0$ is centered.
    Substituting this in the previous line, we have
    \begin{align}
        \phi_{X,Y}(\xi,\eta) - \phi_X(\xi)\phi_X(\eta)
        &= \sum_{k\in\Z\setminus\{0\}} \hat g(k) \E_\arr  e(k e^\arr t).
    \end{align}
    From change of variable and integration by parts, for any $v\ne 0$ we have
    \begin{align*}
        \int_a^b e(ve^x) dx &= \int_{e^a}^{e^b} e(vy)\frac{dy}y
        = \frac1{2\pi i v} \bigg[ \frac{e(ve^b)}{e^b}- \frac{e(ve^a)}{e^a} + \int_{e^a}^{e^b} e(vy) \frac{dy}{y^2}  \bigg] \ls e^{-a}/|v|
    \end{align*}
    and hence
    \begin{equation*}
        \E_\arr  e(k e^\arr t) = \frac1{b-a} \int_a^b e(kte^x) dx \lesssim \frac{e^{-a}}{|k|t(b-a)}\,.
    \end{equation*}
    Combining with Cauchy--Schwarz and the Plancherel theorem,
    \begin{align*}
        |\phi_{X,Y}(\xi,\eta) - \phi_X(\xi)\phi_X(\eta)|
        &\le \|g_0\|_{L^2([0,1))} \bigg( \sum_{k\ne 0} |\E_\arr  e(k e^\arr t)|^2\bigg)^{1/2}
        \lesssim \frac{e^{-a}}{t(b-a)} \lesssim \eps^B
    \end{align*}
    by our assumption on $t$. Substituting this bound in \eqref{indep3} and combining with \eqref{indep1}, \eqref{indep2} and \eqref{indep4}, we obtain
    \[
    \mu_t(I_1\times I_2) \le (1+ O(\delta^{1/C}/\eps) + O(\eps^{B+1}/\delta))\nu(I_1)\nu(I_2).
    \]
    The claim now follows by optimizing $\delta$ and taking $B$ sufficiently large. 
\end{proof}

\begin{proof}
    [Proof of \Cref{prop:decouple}]
    For each $k\in[d]$ fix a partition of the range $\cR$ of $f$ into intervals $\{I_q^k\}_q$ of length 
    \[
    \Delta_k \in [1,2]\cdot \min (k^{-2}, \delta^{2\lam})\,.
    \]
    For $q=(q_1,\dots, q_d)$ denote the box $B_q:= \prod_{k=1}^d I_{q_k}^k\subset \R^d$.
    For $u=s,t$ set 
    \begin{align*}
        \cX_u&:= \bigg\{ x\in \cR^d: \sum_{k=1}^m x_k \in J_m^u \quad \forall\,m\in [d]\bigg\}\\
        Q_u&:= \big\{ q: B_q\cap \cX_u \ne \emptyset\big\}\,.
    \end{align*}
    Thus $\cX_u\subset \bigcup_{q\in Q_u}B_q$. 
    We abbreviate the random variables and random vector
    \[
    f_k^u:= f(e^{T_k}u + U_k)\, , \; f^u=(f_k^u)_{k=1}^d
    \]
    so that $F_m^u=\sum_{k=1}^m f_k^u$. 
    By independence of $(f_k^s,f_k^t)$ across $k\in[d]$, we have
    \begin{align}
        &\P \big( F_m^s\in J_m^s, F_m^t\in J_m^t \quad\forall m\in [d]\big) \notag\\
        &\quad= \P( f^s\in \cX_s \,, \; f^t\in \cX_t)   \notag\\
        &\quad\le \sum_{q^s\in Q_s, q^t\in Q_t} \P( f^s\in B_{q^s}\,,\; f^t\in B_{q^t}) \notag\\
        &\quad = \sum_{q^s\in Q_s, q^t\in Q_t} \prod_{k=1}^d \P\big( f^s_k\in I^k_{q^s_k}\,,\; f^t_k\in I^k_{q^t_k} \big)\,.    \label{Pbox1}
    \end{align}
    Now we use \Cref{lem:block-indep} to show the events $f^s_k\in I^k_{q^s_k}$, $f^t_k\in I^k_{q^t_k}$ are approximately independent. 
    For this we must verify that
    \begin{equation*}
        |s-t| \ge \Delta_k^{-B} \delta^{-1} e^{-a_k}
    \end{equation*}
    for each $k\in[d]$, where $B=B(A)$ is the constant from \Cref{lem:block-indep}. 
    By our assumption $|s-t|\ge \delta^{-C} e^{-a_1}$ it suffices to verify 
    \begin{equation}\label{decouple-verify}
        \Delta_k^Be^{a_k-a_1} \ge \delta^{C-1}
    \end{equation}
    for $C,\alpha$ sufficiently large. 
    For $k\le \delta^{-K}$ we can lower bound the left hand side above by
    \[
    \Delta_k^B \ge \min( k^{-2}, \delta^{2\lam})^B \ge \min(\delta^{2K}, \delta^{2\lam})^B
    \]
    giving \eqref{decouple-verify} upon taking $C:= 1+ 2B(\lam+K)$. 
    For $k>\delta^{-K}$ we instead bound the left hand side of \eqref{decouple-verify} by
    $
    \min( k^2,\delta^{2\lam})^B k^\alpha \ge \delta^{2B\lam} 
    $
    for $\al\ge 2B$, and \eqref{decouple-verify} follows from our choice of $C$.
    From \Cref{lem:block-indep} we hence have
    \[
    \P\big( f^s_k\in I^k_{q^s_k}\,,\; f^t_k\in I^k_{q^t_k} \big)
    \le (1+ O(\Delta_k^A)) 
    \P\big( f^s_k\in I^k_{q^s_k}\big)\,\P\big( f^t_k\in I^k_{q^t_k} \big)
    \]
    for each $k\in[d]$. 
    Inserting this in \eqref{Pbox1} we obtain
    \begin{align}
        &\P \big( F_m^s\in J_m^s, F_m^t\in J_m^t \quad\forall m\in [d]\big) \notag\\
        &\qquad \le \prod_{k=1}^d (1+ O(\Delta_k^A)) 
        \sum_{q^s\in Q_s, q^t\in Q_t} \prod_{k=1}^d \P\big( f^s_k\in I^k_{q^s_k}\big)\,\P\big( f^t_k\in I^k_{q^t_k} \big) \notag\\
        &\qquad \le \prod_{k=1}^d (1+ O(\Delta_k^A)) 
        \bigg( \sum_{q^s\in Q_s}\prod_{k=1}^d \P\big( f^s_k\in I^k_{q^s_k}\big) \bigg) \bigg( \sum_{q^t\in Q_t}\prod_{k=1}^d \P\big( f^t_k\in I^k_{q^t_k} \big) \bigg)\notag\\
        &\qquad= \prod_{k=1}^d (1+ O(\Delta_k^A)) 
        \P \bigg( f^s \in \bigcup_{q^s\in Q_s} B_{q^s}\bigg) \P \bigg( f^t \in \bigcup_{q^t\in Q_t} B_{q^t}\bigg)\,.\notag
    \end{align}
    One easily verifies
    \begin{equation}    \label{sum-DeltakA}
        \sum_{k=1}^d \Delta_k^A \ls 2^A \delta^{\lam(2A-1)}
    \end{equation}
    so that
    \begin{align*}
         \prod_{k=1}^d (1+ O(\Delta_k^A))
         &\le  \exp \big( O( 2^A\delta^{\lam(2A-1)})\big) = 1+ O(\delta)^{\lam A} \,.
    \end{align*}
    It only remains to show that
    \begin{equation}    \label{fu-goal}
         f^u \in \bigcup_{q\in Q_u} B_{q^u}
         \quad\Longrightarrow \quad
         F_m^u \in J_m^u(O(\delta^\lam)) \quad\forall m\in[d]\,.
    \end{equation}
    For this, note that if $x\in B_q$ for some $q\in Q_u$ then by definition there is a point $y\in B_q\cap \cX_u$, for which we have
    \[
    \sum_{k=1}^m y_k \in J_m^u
    \]
    for each $m\in[d]$. 
    Since $B_q$ has side lengths $\Delta_k$, we have
    \[
    \bigg|\sum_{k=1}^m x_k - y_k \bigg|\le 
    \sum_{k=1}^m |x_k -y_k| \le \sum_{k=1}^d \Delta_k \ls \delta^\lam
    \]
    where we applied \eqref{sum-DeltakA} with $A=1$. It follows that $\sum_{k=1}^m x_k\in J_m^u(O(\delta^\lam))$ for all $m$, which gives  \eqref{fu-goal} to complete the proof. 
    \end{proof}

\section{Barrier probability estimates}
\label{sec:barrier}

Let $(B_t)_{t\ge0}$ be a standard Brownian motion. Its values at integer times form the standard Gaussian random walk. For a process $X$ with $X_0=0$, deterministic paths $w$, integers $0\le a<b$, and an interval $J\subseteq\R$, define
\begin{align}
    p_{a,b;X}^{x;J}(w):=\P\Big(&\max_{0\le j\le b-a}\{X_j-(w_{a+j}-w_a)\}\le x,\,X_{b-a}-(w_b-w_a)-x\in J\Big).
    \label{eqn:defn-ballot}
\end{align}
The connection with the notation \eqref{beta-w} of Section~\ref{sec:lower} is
\begin{equation}\label{def:beta}
    \beta_n(w)=p_{n,B}^{10;10,1}(w/\sigma_*).
\end{equation}
For continuous $X,w$, define
\begin{align}
    \overline p_{a,b;X}^{x;J}(w):=\P\Big(&\sup_{s\in [0,b-a]}\{X_s-(w_{a+s}-w_a)\}\le x,\, X_{b-a}-(w_b-w_a)-x\in J\Big). \label{eqn:defn-continuous-ballot}
\end{align}
Unless a continuous path is specified, $w$ is linearly interpolated between integer times. On the other hand, a Brownian motion is not replaced by its linear interpolation. 
We use the abbreviations
\[
    p_{n,X}^{x;J}:=p_{0,n;X}^{x;J},\qquad p_{a,b;X}^{x;y,z}:=p_{a,b;X}^{x;[-y-z,-y]},
\]
and the same abbreviations for $\overline p$. The argument $w$ is unchanged in these abbreviations. When the paths are random, these expressions mean evaluation at their realized paths, with probability taken only over $X$. If the paths are measurable with respect to an environment $\sigma$-field independent of $X$, this is also the conditional probability given that environment.

\begin{prop}\label{prop:LB-wiggle}
    Let $\sigma>0$ and let $W=(W_j)_{j\ge1}$ be a centered random walk with iid steps of unit variance satisfying $\E[e^{\ee |W_1|}]<\infty$ for some $\ee>0$. Let
    $\wt W$ be any process such that $|W_j-\wt W_j|\le \frac{\sigma}{40}$ a.s.\ for all $1\le j\le n$. For any fixed $\eta>0$, 
    \[
        \inf_{\substack{x\in [5,\log^5 n] \\ y\in[5,\log^5 n]\\ z\ge \frac1{10}}} p_{n,B}^{x;y,z}(\wt W/\sigma) \ge n^{-\phi-\eta} \qquad a.a.s.
    \]
    where $\phi = \frac12+2\gamma(1/\sigma)$.
    The same bound holds with $p_{n,B}^{x;y,z}(\wt W/\sigma)$ replaced by its conditional expectation given $W$.
    Note that
    \[
        \E\big[ p_{n,B}^{x;y,z}(\wt W/\sigma) \,\big|\, W\big]= \P\Big(\, \max_{j\le n} B_j -\wt W_j/\sigma \le x, \; x-y-z \le  B_n - \wt W_n/\sigma \le x-y \,\Big|\, W\,\Big) \,.
    \]
\end{prop}

To prove \Cref{prop:LB-wiggle} we use the following two lemmas. 
The first is an adaptation of \cite[Theorem 3.3, Equation (3.5)]{MaMi:branching} to our setting.

\begin{lemma} \label{lem:MM3.3}
    Let $W=(W_j)_{j\ge1}$ be a centered random walk with iid steps of unit variance satisfying $\E[e^{\eta |W_1|}]<\infty$ for some $\eta>0$. 
    For any fixed $\sigma,\ee>0$, 
    \[
    \inf_{\log^4 n\le x,z\le \log^5n}
    p_{n,B}^{x;0,z}(W/\sigma)
    \ge n^{-\phi-\ee}\qquad a.a.s.
    \]
    where $\phi = \frac12+2\gamma(1/\sigma)$.
\end{lemma}

\begin{proof}
    We follow the proof of \cite[Theorem~3.3, Equation (3.5)]{MaMi:branching}, with the same reduction as in \Cref{lem:quenched.upper}. 
    Here we take $T_j=\sigma B_j-W_j$ conditional on $W$, so the barrier height and endpoint width in \cite[Theorem~3.3, Equation (3.5)]{MaMi:branching} are $\sigma x$ and $\sigma z$, respectively. With $\alpha=0$, the event in \cite[Theorem~3.3, Equation (3.5)]{MaMi:branching} is contained in the event defining $p_{n,B}^{x;0,z}(W/\sigma)$. 
    Similarly as in \Cref{lem:quenched.upper}, it remains to verify the local estimate \cite[Lemma 3.6, Equation (3.11) and Lemma 3.7]{MaMi:branching} for an endpoint interval of width $\sigma(\log n)^4$.

    Under the exponential moment assumption on $W_1$, we have
    \[
    \lim_{n\to\infty}\P\!\left(\max_{1\le j\le n}|W_j-W_{j-1}|\ge \frac{2}{\eta}\log n\right)= 0.
    \]
    Hence
    \[
    \max_{1\le j\le n}\frac{\E[|\sigma (B_j-B_{j-1})-(W_j-W_{j-1})|^3\mid W]}{\Var\left(\sigma (B_j-B_{j-1})-(W_j-W_{j-1})\middle| W\right)} = O_{\eta,\sigma}(\log^3 n)\qquad a.a.s.
    \]
    Similarly as in Lemma~\ref{lem:quenched.upper}, this gives a Berry--Esseen error of order $O_{\P}((\log n)^3/\sqrt n)$ in place of the almost sure $O(n^{-1/2})$ error in \cite[Equation (3.12)]{MaMi:branching}.
    In the proof of the lower bound \cite[Theorem~3.3, Equation (3.5)]{MaMi:branching}, we take $a_n=\sigma(\log n)^4$ and $x_n=\sigma(\log n)^5$.
    Then, in the proofs of \cite[Lemma~3.6, equation (3.11), and Lemma~3.7]{MaMi:branching}, the Gaussian mass of an interval of width $a_n$ within $A\sqrt{n}$ of the conditional mean for each fixed $A>0$ is of order $C_{A,\sigma}\frac{a_n}{\sqrt{n}}$. Since the Berry--Esseen error is of a smaller order compared to it, the same arguments towards \cite[Lemma~3.6, equation (3.11), and Lemma~3.7]{MaMi:branching} still apply.
    The rest of the proof of \cite[Theorem~3.3, Equation (3.5)]{MaMi:branching} remains unchanged and gives the claimed bound.
\end{proof}

\begin{lemma}[Gaussian barrier comparison]\label{lem:Kriech1}
    Let $W$ be a centered random walk with iid increments, $W_0=0$, $\E[W_1^2]<\infty$, and $\E[e^{\eta|W_1|}]<\infty$ for some $\eta>0$, linearly interpolated between integer times.
    Let $B$ be a standard Brownian motion independent of $W$.
    Fix $D>0$, $C_1, C_2\ge0$, and $0<c_1\le c_2<\infty$.
    For every $\delta>0$, there exist $C_\delta\ge1$ and $N_\delta\in\N$ with the following property.
    For every pair of deterministic integers $0\le a<b$ with $N:=b-a\ge N_\delta$, there is a $\sigma(W_{a+s}-W_a:s\le N)$-measurable event of probability at least $1-\delta$ on which
    \begin{equation}\label{eqn:Gaussian-comparison}
        (xy)^{-C_\delta} \le \frac{\overline p_{a,b;B}^{x;y,z}(W-h)}{p_{a,b;B}^{10;10,1}(W)} \le \frac{p_{a,b;B}^{x;y,z}(W-h)}{p_{a,b;B}^{10;10,1}(W)} \le (xy)^{C_\delta},
    \end{equation}
    uniformly for $x,y\in[4,(\log N)^D]$, $z\in[c_1,c_2]$, and all continuous curves $h$ on $[a,b]$ satisfying
    \begin{align}
        c&:=\frac{h_b-h_a}{N},\qquad |c|\le \frac{C_1}{\sqrt{N}},\notag\\
        |h_{a+s}-h_a-cs|&\le C_2\left(((1+s)\wedge(1+N-s))^{1/6}-1\right),\qquad s\in[0,N].\label{eqn:Gaussian-comparison-curve}
    \end{align}
\end{lemma}
This is a version of the comparisons in \cite[Section 12]{Kriechbaum:time} for the moment assumptions used here.
The event is for each fixed block $[a,b]$, not for all pairs with the same gap simultaneously.
We prove the lemma after \Cref{lem:Gaussian-crude-LB}, which provides the polynomial lower bounds needed to absorb the errors in the comparison argument.

\begin{proof}[Proof of \Cref{prop:LB-wiggle} assuming \Cref{lem:Kriech1}]
    For $z\ge 0.1$ we have $p_{n,B}^{x;y,z}(\wt W/\sigma)\ge p_{n,B}^{x;y,0.1}(\wt W/\sigma)$, so for a lower bound we only need to consider $z=0.1$.
    
    By Lemma \ref{lem:MM3.3} with $x=y=\log^4 n$,
    there is an event $\cG_n$ of probability $1-o_\eta(1)$ on which we have
    \[
    p_{n,B}^{\log^4 n+4;4,\log^4 n}(W/\sigma) \ge p_{n,B}^{\log^4 n;0,\log^4 n}(W/\sigma) \ge n^{-\phi-\eta/2}
    \]
    where $\phi = \frac12+2\gamma(1/\sigma)$.
    By the union bound, for any $c>0$ we have
    \[
    p_{n,B}^{\log^4 n+4;4,\log^4 n}(W/\sigma) \le \sum_{0\le \ell\le \log^4 n}p_{n,B}^{\log^4 n+4;\ell+4,1}(W/\sigma).
    \]
    Fixing $\delta>0$ arbitrary, by the upper and lower bounds in \Cref{lem:Kriech1} there is an event $\cG_n'(\delta)$ of probability at least $1-\delta-o_{\delta}(1)$ on which we have
    \[
    p_{n,B}^{\log^4 n+4;\ell+4,1}(W/\sigma) \le (\log^{O_\delta(1)}n )p_{n,B}^{x;y,z}(W/\sigma)
    \]
    for all $0\le \ell\le 10 \log^4 n$, $x, y\in[4,\log^5 n]$, and $z\in[0.01, 100]$.
    Hence on $\cG_n'(\delta)$, for any $x, y\in[4,\log^5 n]$ and $z\in[0.01, 100]$ we have
    \[
    p_{n,B}^{\log^4 n+4;4,\log^4 n}(W/\sigma)  \le (\log^{O_\delta(1)}n )p_{n,B}^{x;y,z}(W/\sigma).
    \]
    Thus, for all $n$ sufficiently large depending on $\eta,\delta$, 
    \[
    \P\left( \inf_{x, y\in [4,\log^5 n], z\in [0.01, 100]}p_{n,B}^{x;y,z}(W/\sigma) \ge n^{-\phi-\eta}\right) \ge \P( \cG_n\cap\cG_n'(\delta)) \ge 1- \delta - o_{\delta,\eta}(1),
    \]
    which implies
    \[
    \liminf_{n\to\infty}\P\left(\inf_{x, y\in [4,\log^5 n], z\in [0.01, 100]}p_{n,B}^{x;y,z}(W/\sigma) \ge n^{-\phi-\eta}\right) \ge 1-\delta.
    \]
    
    Since $|\frac{W_j-\wt{W}_j}{\sigma}|\le \ee_0=\frac1{40}$, 
    for any $1\le j\le n$, we obtain by monotonicity that
    \[
    p_{n,B}^{x-\ee_0;y,0.1-2\ee_0}(W/\sigma) \le p_{n,B}^{x;y,0.1}(\wt W/\sigma)
    \]
    for any $x,y>0$.
    This inclusion is pathwise and also holds after taking conditional expectation given $W$.
    Therefore, we have
    \[
    \liminf_{n\to\infty}\P\left( \inf_{x, y\in [5,\log^5n], z\ge 0.1}p_{n,B}^{x;y,z}(\wt W/\sigma)\ge n^{-\phi-\eta}\right) \ge 1-\delta.
    \]
    Taking $\delta\downarrow0$, the claim follows.
\end{proof}

We repeatedly use the following crude lower bound estimate of the Gaussian quenched barrier probability. 
A similar crude lower bound can be found in \cite[Lemma 11.6]{Kriechbaum:time}.
\begin{lemma}[Crude Gaussian lower bound]\label{lem:Gaussian-crude-LB}
    Let $W$ be a centered random walk with iid increments and $\E[W_1^2]<\infty$, linearly interpolated between integer times. Let $B$ be a standard Brownian motion independent of $W$. Fix $0<\alpha<1/2$, positive $C_1,C_2<\infty$, and $0<c_0\le1$.
    For every $\delta>0$, there exist $C_\delta<\infty$ and $N_\delta\in\N$ with the following property. 
    For every pair of deterministic integers $0\le a<b$ with $N:=b-a\ge N_\delta$, there is a $\sigma(W_{a+s}-W_a:0\le s\le N)$-measurable event of probability at least $1-\delta$ on which
    \begin{equation}\label{eqn:Gaussian-crude-LB}
        \overline p_{a,b;B}^{x;y,z}(W-h)\ge N^{-C_\delta}
    \end{equation}
    holds uniformly for $1\le x\le\sqrt N$, $0\le y\le\sqrt N$, $z\ge c_0$, and all continuous curves $h$ on $[a,b]$ satisfying
    \begin{align}
        c&:=\frac{h_b-h_a}{N},\qquad |c|\le C_1\sqrt{\frac{\log N}{N}},\notag\\
        f(s)&:=h_{a+s}-h_a-cs\le C_2\left(((1+s)\wedge(1+N-s))^\alpha-1\right),\qquad 0\le s\le N.\label{eqn:Gaussian-crude-curve}
    \end{align}
\end{lemma}
\begin{proof}
    By stationarity, it suffices to take $a=0$ and $b=N$.
    We bound the quenched probability from below by \cite[Theorem 10.2]{Kriechbaum:time}. To match its assumptions, we make the following reductions.
    By monotonicity, we may replace $z$ with $c_0$. Recall that $f(s)=h_s-h_0-cs$ and $f(0)=f(N)=0$.
    By Girsanov's formula, with $\cA$ denoting the barrier event defining $\overline p_{N,B}^{x;y,c_0}(W-f)$,
    \[
        \overline p_{N,B}^{x;y,c_0}(W-h)=\E\left[\ind_{\cA}\exp\left\{cB_N-\frac{c^2N}{2}\right\}\,\middle|\,W\right].
    \]
    On this event, $B_N\in x+W_N+[-y-c_0,-y]$, so
    \[
        \left|cB_N-\frac{c^2N}{2}\right|\le |c|(x+y+c_0+|W_N|)+\frac{Nc^2}{2}.
    \]
    Set $f_*(s):=C_2\big(((1+s)\wedge(1+N-s))^\alpha-1\big)$; since $f\le f_*$ and both vanish at the endpoints, monotonicity gives
    \begin{equation}\label{eqn:Gaussian-crude-reduction}
        \overline p_{N,B}^{x;y,z}(W-h)\ge e^{-|c|(x+y+c_0+|W_N|)-Nc^2/2}\overline p_{N,B}^{x;y,c_0}(W-f_*).
    \end{equation}

    Before applying \cite[Theorem 10.2]{Kriechbaum:time}, we remark that the proof of this theorem also applies to an endpoint interval of fixed length $c_0$, with constants depending on $c_0$. Next, we introduce the following notation. Let $k_1:=\lfloor\log_2(N/3)\rfloor$ and define
    \begin{align}
        \Xi_j(W)&:=\min\left\{\inf_{0\le s\le2^{j+1}}W_s,\inf_{N-2^{j+1}\le s\le N}(W_s-W_N)\right\},\qquad 0\le j<k_1,\notag\\
        \Xi_{k_1}(W)&:=\min\left\{\inf_{0\le s\le N}W_s,\inf_{0\le s\le N}(W_s-W_N),0\right\}.\label{eqn:defn-Xi}
    \end{align}
    Let $\Xi_0(-f_*):=-C_2\alpha$, and for $1\le j\le k_1$ define $\Xi_j(-f_*)$ by replacing $W$ with $-f_*$ above.
    The quantities in \cite[Definition 10.1]{Kriechbaum:time} are
    \[
        t_1:=2^{k_1},\qquad t_2:=N-t_1,\qquad \Xi_0:=\Xi_0(W)+\Xi_0(-f_*),\qquad \Delta_j:=\Delta_j(W)+\Delta_j(-f_*),
    \]
    where, for $0\le j<k_1$,
    \begin{equation}
        \Delta_j(W):=\frac{\Xi_{j+1}(W)-\Xi_j(W)}{2^j},\qquad \Delta_j(-f_*):=\frac{\Xi_{j+1}(-f_*)-\Xi_j(-f_*)}{2^j}.\label{eqn:defn-Delta}
    \end{equation}
    Since $f_*(N-s)=f_*(s)\le C_2\alpha s$ for $0<s\le1$, our choice satisfies the required inequality
    \[
        \Xi_0(-f_*)\le\inf_{0<s\le1}\min\left\{\frac{-f_*(s)}{s},\frac{-f_*(N-s)}{s}\right\}.
    \]

    Since $-\sqrt N\le-x,-y\le0$ and $x\ge1$, \cite[Theorem 10.2]{Kriechbaum:time} gives deterministic constants $C,\gamma_0>0$ such that, uniformly over the stated ranges of $x,y$,
    \begin{equation}\label{eqn:Gaussian-crude-pathwise}
        \overline p_{N,B}^{x;y,c_0}(W-f_*)\ge N^{-\gamma_0}e^{-CR_N(W,-f_*)},\qquad a.s.
    \end{equation}
    where
    \begin{equation}\label{eqn:Gaussian-crude-R}
        R_N(W,-f_*):=1+\sum_{j=0}^{k_1-1}\left(2^j\Delta_j^2+2^{j/2}|\Delta_j|\right)+\Xi_0(\Xi_0-4)+\frac{W_N^2+(\sqrt{t_1}+2\sqrt N)|W_N|}{t_2-t_1}.
    \end{equation}
    Here we used $f_*(N)=0$ and $x+y\le2\sqrt N$. Since $\Xi_0\le0$, $R_N(W,-f_*)\ge0$.
    We claim that there is a deterministic constant $D_\delta<\infty$ such that for sufficiently large $N$,
    \begin{equation}\label{eqn:Gaussian-crude-event}
        \cG_{N}(\delta):=\{R_N(W,-f_*)\le D_\delta\log N\},\qquad \P\big(\cG_{N}(\delta)\big)\ge1-\delta.
    \end{equation}
    On this event, $W_N^2/N\le R_N(W,-f_*)$ and $Nc^2\le C_1^2\log N$ give, there is a $C_{\delta}<\infty$ such that
    \[
        |c|(x+y+c_0+|W_N|)+Nc^2/2\le 3C_1\sqrt{\log N}+C_1\sqrt{D_\delta}\log N+(C_1^2/2)\log N\le C_\delta\log N.
    \]
    Combining \eqref{eqn:Gaussian-crude-reduction} and \eqref{eqn:Gaussian-crude-pathwise} therefore proves \eqref{eqn:Gaussian-crude-LB}, after increasing $C_\delta$.

    It remains to prove \eqref{eqn:Gaussian-crude-event}. We assume $N$ is sufficiently large in the rest of the proof, by enlarging it when necessary. Split the exponent as in the proof of \cite[Lemma 11.6]{Kriechbaum:time} up to (11.15):
    \begin{align*}
        R_N^{(1)}&:=1+\sum_{j=0}^{k_1-1}\left(2^j\Delta_j(-f_*)^2+2^{j/2}|\Delta_j(-f_*)|\right)+\frac{W_N^2+(\sqrt{t_1}+2\sqrt N)|W_N|}{t_2-t_1}\\
        &\quad+2^{k_1-1}\Delta_{k_1-1}(W)^2+2^{(k_1-1)/2}|\Delta_{k_1-1}(W)|,\\
        R_N^{(2)}&:=\Xi_0(\Xi_0-4)+\sum_{j=0}^{k_1-2}\left(2^j\Delta_j(W)^2+2^{j/2}|\Delta_j(W)|\right).
    \end{align*}
    Then $R_N(W,-f_*)\le2(R_N^{(1)}+R_N^{(2)})$.
    
    The proof of \cite[Lemma 11.6]{Kriechbaum:time} up to (11.15) directly implies that $\P(R_N^{(1)}\le C_\delta^{(1)}\log N)\ge1-\delta/2$ for a sufficiently large $C_\delta^{(1)}$.
    For $0\le j\le k_1-2$, define
    \[
        a_{j,N,1}:=\sup_{2^{j+1}\le s\le2^{j+2}}|W_s-W_{2^{j+1}}|,\qquad a_{j,N,2}:=\sup_{N-2^{j+2}\le s\le N-2^{j+1}}|W_s-W_{N-2^{j+1}}|.
    \]
    By the definitions of $\Xi_j(W)$ and $\Delta_j(W)$,
    \[
        2^j|\Delta_j(W)|=|\Xi_{j+1}(W)-\Xi_j(W)|\le\max\{a_{j,N,1},a_{j,N,2}\}.
    \]
    Since $W$ has centered iid increments with finite variance, Doob's $L^2$ maximal inequality and Cauchy--Schwarz, applied also to the time-reversed increments, give a deterministic $C<\infty$ such that, uniformly over $0\le j\le k_1-2$ and $i\in\{1,2\}$,
    \[
        \E[a_{j,N,i}^2]\le C2^j,\qquad \E[a_{j,N,i}]\le C2^{j/2}.
    \]
    To control the $\Xi_0$-term, the intervals defining $\Xi_0(W)$ have length two, so Doob's inequality and stationarity give
    \[
        \E[\Xi_0(W)^2]\le\E\left[\sup_{0\le s\le2}|W_s|^2\right]+\E\left[\sup_{0\le s\le2}|W_N-W_{N-s}|^2\right]\le8\E[W_2^2].
    \]
    which implies $\E[\Xi_0(\Xi_0-4)]\le C$.
    Consequently (constants may vary from line to line),
    \begin{align*}
        \E[R_N^{(2)}]&=\E[\Xi_0(\Xi_0-4)]+\sum_{j=0}^{k_1-2}\E\left[2^j\Delta_j(W)^2+2^{j/2}|\Delta_j(W)|\right]\\
        &\le C+C\sum_{j=0}^{k_1-2}\E\left[2^{-j}(a_{j,N,1}^2+a_{j,N,2}^2)+2^{-j/2}(a_{j,N,1}+a_{j,N,2})\right]\\
        &\le Ck_1\le C\log N.
    \end{align*}
    Therefore, by Markov's inequality, there is a deterministic constant $C_\delta^{(2)}<\infty$ such that
    \[
        \P\left(R_N^{(2)}>C_\delta^{(2)}\log N\right)\le\delta/2.
    \]
    Combining the two probability bounds proves \eqref{eqn:Gaussian-crude-event} with $D_\delta=2(C_\delta^{(1)}+C_\delta^{(2)})$.
\end{proof}

\begin{proof}[Proof of \Cref{lem:Kriech1}]
    By stationarity, it suffices to consider $a=0, b=N$.
    Let $f_s:=h_s-h_0-cs$ and $g_s:=\left(((1+s)\wedge(1+N-s))^{1/6}-1\right),  \forall s\in[0,N]$. Condition \eqref{eqn:Gaussian-comparison-curve} then says $|f|\le C_2 g$ on $[0,N]$.

    Recall $|c|\le C_1 N^{-1/2}$. Similarly as in the Girsanov change-of-measure argument leading to \eqref{eqn:Gaussian-crude-reduction}, on the event $\cG_1:=\{|W_N|\le C_{\delta}\sqrt{N}\}$ with $C_{\delta}<\infty$ chosen large enough such that $\P(\cG_1)\ge 1-\delta/2$, uniformly for $x, y, z$ in the stated parameter ranges and all $h$ satisfying the stated condition, we have
    \begin{equation}\label{eqn:Gaussian-comparison-reduction}
        \overline p_{N,B}^{x;y,z}(W-h) = e^{O_{\delta}(1)}\overline p_{N,B}^{x;y,z}(W-f),\qquad p_{N,B}^{x;y,z}(W-h) = e^{O_{\delta}(1)} p_{N,B}^{x;y,z}(W-f).
    \end{equation}
    Therefore, it suffices to show \eqref{eqn:Gaussian-comparison} for curves $f$ satisfying $f_0=f_N=0$ and $|f|\le C_2 g$. We first prove it assuming the following claims, which are analogues of the comparison lemmas in \cite[Section 12]{Kriechbaum:time}, with the moment assumptions, curves, and parameter ranges needed here.
    We follow the same proof strategy and give the necessary modifications afterwards.

    \begin{claim}[Move starting and ending points]\label{claim:start-end-comparison}
        For every $\delta>0$, there exist $C_\delta\ge1$ and $N_\delta\in\N$ such that, for every integer $N\ge N_\delta$, there is a $\sigma(W_s:0\le s\le N)$-measurable event of probability at least $1-\delta$ on which
        \[
        (xy)^{-C_\delta}\le \frac{\overline p_{N,B}^{x;y,z}(W-f)}{\overline p_{N,B}^{10;10,1}(W-f)}\le (xy)^{C_\delta}
        \]
        holds uniformly for $x, y, z$ in the stated ranges, and all continuous curves $f$ on $[0,N]$ satisfying $f_0=f_N=0$ and $|f|\le(C_2+1)g$.
    \end{claim}
    
    \begin{claim}[Remove the additional correction]\label{claim:remove-deterministic}
        For every $\delta>0$, there exist $C_\delta\ge1$ and $N_\delta\in\N$ such that, for every integer $N\ge N_\delta$, there is a $\sigma(W_s:0\le s\le N)$-measurable event of probability at least $1-\delta$ on which
        \[
        C_\delta^{-1}\le \frac{\overline p_{N,B}^{10;10,1}(W-f)}{\overline p_{N,B}^{10;10,1}(W)}\le C_\delta
        \]
        holds uniformly over all continuous curves $f$ on $[0,N]$ satisfying $f_0=f_N=0$ and $|f|\le(C_2+1)g$.
    \end{claim}

    \begin{claim}[Pass to continuous time]\label{claim:discrete-cts-comparison}
        There is a deterministic constant $C\ge1$ such that, for every integer $N\ge2$, every linearly interpolated path $W$, $x,y\ge4$, $z>0$, and every sequence $(f_j)_{j=0}^N$ satisfying $f_0=f_N=0$ and $f_j\ge-C_2g_j$ for $1\le j\le N-1$, we have
        \[
        p_{N,B}^{x;y,z}(W-f)\le C\,\overline p_{N,B}^{x;y,z}(W+(C_2+1)g).
        \]
    \end{claim}

    By monotonicity, we have
    \begin{equation}\label{eqn:comparison-monotonicity}
    \overline p_{N,B}^{x;y,z}(W-C_2 g) \le \overline p_{N,B}^{x;y,z}(W-f) \le p_{N,B}^{x;y,z}(W-f) \le p_{N,B}^{x;y,z}(W+C_2 g).
    \end{equation}
    Apply \Cref{claim:start-end-comparison} and \Cref{claim:remove-deterministic} with error probability $\delta/4$, and let $\cG_2$ be the intersection of their events.
    Then $\P(\cG_2)\ge1-\delta/2$, and the following inequalities hold uniformly over the stated parameter ranges on $\cG_2$ for any sufficiently large $N$ (constants may vary from line to line):
    \begin{align}
        p_{N,B}^{x;y,z}(W+C_2 g) \overset{\text{\Cref{claim:discrete-cts-comparison}}}{\le} &C\overline p_{N,B}^{x;y,z}(W+(C_2+1)g) \notag\\
        \overset{\text{\Cref{claim:start-end-comparison}}}{\le} &C_{\delta}(xy)^{C_{\delta}}\overline p_{N,B}^{10;10,1}(W+(C_2+1)g) \notag\\
        \overset{\text{\Cref{claim:remove-deterministic}}}{\le} &C_{\delta}(xy)^{C_{\delta}}\overline p_{N,B}^{10;10,1}(W) \le C_{\delta}(xy)^{C_{\delta}} p_{N,B}^{10;10,1}(W), \label{eqn:comparison-UBs}
    \end{align}
    and
    \begin{align}
    \overline p_{N,B}^{x;y,z}(W-C_2g) &\overset{\text{\Cref{claim:start-end-comparison}}}{\ge} C_\delta^{-1}(xy)^{-C_\delta}\overline p_{N,B}^{10;10,1}(W-C_2g) \notag\\
    &\overset{\text{\Cref{claim:remove-deterministic}}}{\ge} C_\delta^{-1}(xy)^{-C_\delta}\overline p_{N,B}^{10;10,1}(W+(C_2+1)g) \notag\\
    &\overset{\text{\Cref{claim:discrete-cts-comparison}}}{\ge} C_\delta^{-1}(xy)^{-C_\delta}p_{N,B}^{10;10,1}(W), \label{eqn:comparison-LBs}
    \end{align}
    where the second inequality follows by applying \Cref{claim:remove-deterministic} to both curves $C_2 g$ and $-(C_2+1)g$.
    Combining \eqref{eqn:comparison-monotonicity}, \eqref{eqn:comparison-UBs}, \eqref{eqn:comparison-LBs}, and \eqref{eqn:Gaussian-comparison-reduction} proves the desired comparison on $\cG_1\cap\cG_2$, whose probability is at least $1-\delta$.
    This concludes the proof of \Cref{lem:Kriech1} assuming Claims \ref{claim:start-end-comparison}. \ref{claim:discrete-cts-comparison} and \ref{claim:remove-deterministic}.
    \end{proof}
    
    \begin{proof}[Proof of \Cref{claim:start-end-comparison}]
    The proof follows the arguments of \cite[Lemmas 12.4 and 12.8]{Kriechbaum:time}, using \Cref{lem:Gaussian-crude-LB} in place of Lemma 11.6 there.
    The arguments of \cite[Lemmas 11.1--11.4]{Kriechbaum:time} give tight bounds on the random constants appearing in those proofs, uniformly over $f$ in our stated range.
    Indeed, comparing $f$ with $\pm(C_2+1)g$ and monotonicity of probability give common bounds on the quenched barrier probabilities defining $\gamma$, and the remaining constants are controlled by the formulas in (8.5)--(8.7).
    Note that $\pm(C_2+1)g$ differ from their $h^{\frown}$ and $h^{\smile}$ only by fixed multiplicative constants.
    \end{proof}
    
    \begin{proof}[Proof of \Cref{claim:discrete-cts-comparison}]
    By monotonicity, for $f$ satisfying the stated conditions, we have
    \[
    p_{N,B}^{x;y,z}(W-f) \le p_{N,B}^{x;y,z}(W+C_2 g).
    \]
    It suffices to show that $p_{N,B}^{x;y,z}(W+C_2 g) \le C\,\overline p_{N,B}^{x;y,z}(W+(C_2+1)g)$ for some $C\ge 1$. This follows from the Brownian bridge argument leading to \cite[Lemma 12.2 or equation (11.2)]{Kriechbaum:time}.
    Indeed, in place of their $c''((1+j)^{1/6}-1)$ in (11.12) we have $\min\{g(j-1), g(j)\}$ in our setting for $j\in\{2,\dots, N-1\}$.
    For the first and last unit time-interval, the respective replacement is $g(1)=2^{1/6}-1<4\le x, y$.
    Then we use the standard Brownian bridge estimate $\P(\max_{t\in[0,1]}U_t>a) = e^{-2a^2}, \forall a>0$ (where $(U_t)_{t\in[0,1]}$ is a standard Brownian bridge) and the fact that
    \[
    \sum_{j\ge 2}\exp\Big(-2(j^{1/6}-1)^2\Big) < \infty
    \]
    to conclude the proof.
    \end{proof}    

    \begin{proof}[Proof of \Cref{claim:remove-deterministic}]
            For every admissible $f$, monotonicity gives
    \[
    \overline p_{N,B}^{10;10,1}(W-(C_2+1)g)\le\overline p_{N,B}^{10;10,1}(W-f)\le\overline p_{N,B}^{10;10,1}(W+(C_2+1)g).
    \]
    Thus, we only need to establish an upper bound for $-(C_2+1)g$ and a lower bound for $(C_2+1)g$.
    Their proofs follow the arguments of \cite[Lemmas 12.6 and 12.9]{Kriechbaum:time}, respectively. We describe only the necessary modifications.
    Just like in the proof of \Cref{claim:start-end-comparison}, we use \Cref{lem:Gaussian-crude-LB} in place of \cite[Lemma 11.6]{Kriechbaum:time} in these proofs.
    Meanwhile,
    the following claim replaces
    \cite[Assumptions (9.4, 9.5), Lemma 11.8]{Kriechbaum:time} 
    (which, similarly to \Cref{lem:Kr-small}, does not directly apply in our setting due to a sub-Gaussian hypothesis for the random barrier increments).

    \begin{claim}\label{claim:short-time-uniform-LB}
        Recall $(W_s)_{s\ge 0}$ is defined as linearly interpolating $(W_j)_{j\in\N}$ between integer times.
        For every $\delta>0$, there is $C_\delta<\infty$ such that, for every sufficiently large integer $N$, there is a $\sigma(W_s:0\le s\le N)$-measurable event of probability at least $1-\delta$ on which
        \begin{equation}\label{eqn:simultaneous-short}
        \overline p_{r,B}^{1;1,1}(W)\ge\exp\{-C_\delta\log^2(r+2)\},\qquad 0<r\le N,
        \end{equation}
        and the same estimate holds with $W$ replaced by $\rev W_s:=W_{N-s}-W_N$.
    \end{claim}
    
    Finally, we remark that in \cite[Equations (9.11, 9.12)]{Kriechbaum:time}, $c_2^{-1}(l+1)^{C_4}$ and $\hat{c}_2^{-1}\,(l+1)^{\hat{C}_4}$ should be replaced by $\exp\{C_\delta\log^2(l+2)\}$ using the above claim.
    The choice of $\eta$ remains possible since
    \[
    \sum_{l\ge 1}e^{-cl^{1/3}+C_{\delta}\log^2(l+2)}(l+1)^C < \infty.
    \]
    Also, the ending interval $[-1,0)$ for the event $B_{t,h^{\frown},\chi(t),1}$ in \cite[Definition 9.7]{Kriechbaum:time} should be $[-1,\infty)$, the exponent of $(\lfloor \chi(t)\rfloor-1)$ should be $\frac13$ instead of $\frac43$, and hence one may take $\chi(t)=\log^6(t)$ instead, and the ending interval $J_0$ in the RHS of \cite[Equation (9.25)]{Kriechbaum:time} should be $J_{-1}$ (hence our choice $\overline p_{r,B}^{1;1,1}(W)$ in the above claim).
    
    Nevertheless, the original arguments leading to \cite[Lemmas 12.6 and 12.9]{Kriechbaum:time} still work under these modifications, and they prove \Cref{claim:remove-deterministic}.

    It remains to prove \Cref{claim:short-time-uniform-LB}. We prove \eqref{eqn:simultaneous-short} with error probability $\delta/2$.
    The same argument applies to $\rev W$, whose increments have the law of $-W_1$ which satisfies our assumptions, and a union bound gives both estimates on one event of probability at least $1-\delta$.

    As in \Cref{lem:Gaussian-crude-LB}, we adapt \cite[Theorem 10.2]{Kriechbaum:time} to obtain a quenched (almost sure) lower bound for $\overline p_{r,B}^{1;1,1}(W)$ (similar to \eqref{eqn:Gaussian-crude-pathwise}) for every real $r\in(0,N]$, then we control simultaneously all the $r$-dependent exponents (similar to \eqref{eqn:Gaussian-crude-R}).
    
    However, we need to first modify the definition for the piecewise linear minorant function $g_r$ in \cite[Definition 10.1]{Kriechbaum:time}, since its definition as stated does not necessarily satisfy $g_r(s)\le W_s$ for $r-1<s\le r$, as needed in their Theorem 10.2, when $r$ is non-integer.
    Instead, we define
    \[
    \wt g_r(s) := 
    \begin{cases}
        W_s,\qquad r\in(0,6),\, s\in[0,r];\\
        g_{\lfloor r\rfloor}(s),\qquad r\ge 6,\, s\in[0,\lfloor r\rfloor],\\
        W_s,\qquad r\ge 6,\, s\in[\lfloor r\rfloor, r]
    \end{cases}
    \]
    where for any positive integer $k$, $g_k$ is defined as in \cite[Definition 10.1]{Kriechbaum:time}.
    Since $(W_s)_{s\ge 0}$ is linearly interpolated between integers, $\wt g_r$ is also piecewise linear and satisfies $\wt g_r(s)\le W_s$ for any $0\le s< r$, while $\wt g_r(r)=W_r$.

    When $r\ge 6$, applying monotonicity and then the argument leading to \cite[Theorem 10.2]{Kriechbaum:time} gives, almost surely,
    \begin{equation}
        \overline p_{r,B}^{1;1,1}(W) \ge \P\left(\sup_{s\in[0,r]}B_s-\wt g_r(s)\le 1,\, B_r-\wt g_r(r)\in[-1,0],\, B_{\lfloor r\rfloor}-\wt g_r(\lfloor r\rfloor) \ge -1\middle| W\right) \ge r^{-\gamma}e^{-cR_r},\label{eqn:adapted10.2}
    \end{equation}
    for some $\gamma,c>0$ and (recall $R_N(W,0)$ from \eqref{eqn:Gaussian-crude-R} with $f_*$ replaced by $0$)
    \begin{equation}\label{eqn:adapted-R}
        R_r := R_{\lfloor r\rfloor}(W,0)+|W_{\lceil r\rceil}-W_{\lfloor r\rfloor}|(2+|W_{\lceil r\rceil}-W_{\lfloor r\rfloor}|).
    \end{equation}
    Indeed, the only extra Girsanov change-of-measure cost (compared to removing $g_{\lfloor r\rfloor}$) is 
    \[
    \exp\left(-\int_{\lfloor r\rfloor}^{r} \wt g_r'(s)~d B_s+\frac12\int_{\lfloor r\rfloor}^{r} \wt g_r'(s)^2~d s\right) \ge \exp\left(-(W_{\lceil r\rceil}-W_{\lfloor r\rfloor})(B_r-B_{\lfloor r\rfloor})\right), 
    \]
    and on the event $\{B_r-\wt g_r(r)\in[-1,0],\, B_{\lfloor r\rfloor}-\wt g_r(\lfloor r\rfloor)\in[-1,1]\}$ we have (note $\wt g_r(s)=W_s$ on $[\lfloor r\rfloor, r]$)
    \[
    \exp\left(-(W_{\lceil r\rceil}-W_{\lfloor r\rfloor})(B_r-B_{\lfloor r\rfloor})\right) \ge \exp\left(-|W_{\lceil r\rceil}-W_{\lfloor r\rfloor}|(2+|W_{\lceil r\rceil}-W_{\lfloor r\rfloor}|)\right).
    \]
    When $0<r<6$, a similar Girsanov change-of-measure calculation gives, almost surely,
    \begin{align}
        \inf_{0<r<6}\overline p_{r,B}^{1;1,1}(W) &\ge e^{-\sum_{j=1}^{6}|W_{j}-W_{j-1}|(2+|W_{j}-W_{j-1}|)}\inf_{0<r<6}\P(\sup_{s\le r}|B_s|\le 1,\,B_r\in[-1,0]) =: c_* e^{-R_*},\label{eqn:short r}
    \end{align}
    where $c_*:=\inf_{0<r<6}\P(\sup_{s\le r}|B_s|\le 1,\,B_r\in[-1,0])>0$.

    Now we control $(R_r)_{6\le r\le N}$ and $R_*$ on the following event.
    Let $\cG_{N,\delta}$ be the event
    \[
    |W_j-W_i| \le K_{\delta}\left(\sqrt{(j-i)\log(j+2)}+\log(j+2)\right),\qquad \forall 0\le i<j\le N,
    \]
    where $K_{\delta}>0$ is a large constant to be determined later.
    Since $(W_j)_{j\ge 0}$ starts from zero, has iid centered and sub-exponential increments, there are constants $c_0,c_1>0$ such that for any pair $0\le i<j\le N$ and $v>0$,
    \[
    \P(|W_j-W_i|>v) \le 2\exp\left(-c_0\min\{\frac{v^2}{c_1(j-i)}, v\}\right).
    \]
    By enlarging $K_{\delta}$ according to $c_0, c_1$ if necessary, for any $0\le i<j\le N$ we have
    \[
    \P\left(|W_j-W_i|>K_{\delta}\left(\sqrt{(j-i)\log(j+2)}+\log(j+2)\right)\right) \le \exp\left(-\sqrt{K_{\delta}}\log(j+2)\right),
    \]
    and hence by the union bound we have
    \begin{equation}
    \P(\cG_{N,\delta}^c) \le \sum_{j=1}^{N}\frac{j}{(j+2)^{\sqrt{K_{\delta}}}} \le \sum_{j\ge 3} j^{-\sqrt{K_{\delta}}+1} \le \delta,\label{eqn:typical-W-path}
    \end{equation}
    for sufficiently large $K_{\delta}$. Note $K_{\delta}$ depends on the law of $W_1$ and $\delta$, but is independent of $N$.

    In the sequel, we derive estimates assuming $\cG_{N,\delta}$ holds. The absolute constant $C<\infty$ may change from line to line.
    Direct computation shows that, there is an absolute constant $C<\infty$ such that
    \begin{equation}
        R_* \le CK^2_{\delta}.\label{eqn:R*-ub}
    \end{equation}
    and also
    \begin{equation}
        |W_{\lceil r\rceil}-W_{\lfloor r\rfloor}|(2+|W_{\lceil r\rceil}-W_{\lfloor r\rfloor}|) \le CK_{\delta}^2\log^2(r+2).\label{eqn:Rr-extra-ub}
    \end{equation}
    Recall \eqref{eqn:defn-Xi} and \eqref{eqn:Gaussian-crude-R}. We have
    \begin{align}
        R_{\lfloor r\rfloor}(W,0) &= 1+\sum_{j=0}^{\lfloor\log_2(r/3)\rfloor-1}\left(\frac{|\Xi_{j+1}(W)-\Xi_j(W)|}{2^{j/2}}+\frac{|\Xi_{j+1}(W)-\Xi_j(W)|^2}{2^{j}}\right) \notag\\
        &+\Xi_0(W)(\Xi_0(W)-4)+\frac{W_{\lfloor r\rfloor}^2+O(\sqrt{r})|W_{\lfloor r\rfloor}|}{\lfloor r\rfloor-2^{\lfloor\log_2(r/3)\rfloor}+1}.
    \end{align}
    Enlarging $C$ if necessary, on the event $\cG_{N,\delta}$, we have
    \begin{align}
        |\Xi_0(W)(\Xi_0(W)-4)| \le CK_{\delta}^2\log^2(r+2),\qquad \frac{W_{\lfloor r\rfloor}^2+O(\sqrt{r})|W_{\lfloor r\rfloor}|}{\lfloor r\rfloor-2^{\lfloor\log_2(r/3)\rfloor}+1} \le CK^2_{\delta}\log(r+2).\label{eqn:Rr-easy-ub}
    \end{align}
    Moreover, for any $0\le j\le \lfloor \log_2(r/3)\rfloor-2$ we have
    \[
    |\Xi_{j+1}(W)-\Xi_j(W)|\le\max_{2^{j+1}\le i\le 2^{j+2}}|W_i-W_{2^{j+1}}|\,\vee \max_{\lfloor r\rfloor-2^{j+2}\le i\le \lfloor r\rfloor-2^{j+1}}|W_i-W_{\lfloor r\rfloor-2^{j+1}}|,
    \]
    so on the event $\cG_{N,\delta}$, enlarging $C$ if necessary,
    \begin{align}
        \frac{|\Xi_{j+1}(W)-\Xi_j(W)|}{2^{j/2}} &\le 2K_{\delta}\max\left\{\Big(\sqrt{\log(2^{j+2}+2)}+2^{-j/2}\log(2^{j+2}+2)\Big),\right.\notag\\
        &\left.\Big(\sqrt{\log(r-2^{j+1}+2)}+2^{-j/2}\log(r-2^{j+1}+2)\Big)\right\}\notag\\
        &\le CK_{\delta}\left(\sqrt{\log(r+2)}+2^{-j/2}\log(r+2)\right).\label{eqn:Rr-Deltaj-ub1}
    \end{align}
    For $j=\lfloor \log_2(r/3)\rfloor-1$, since $\Xi_{\lfloor \log_2(r/3)\rfloor}(W)\le \Xi_{\lfloor \log_2(r/3)\rfloor-1}(W)\le 0$, we have
    \begin{align}
        \frac{|\Xi_{\lfloor \log_2(r/3)\rfloor}(W)-\Xi_{\lfloor \log_2(r/3)\rfloor-1}(W)|}{2^{(\lfloor \log_2(r/3)\rfloor-1)/2}} \le \frac{2\max_{1\le j\le\lfloor r\rfloor}|W_j|}{2^{(\lfloor \log_2(r/3)\rfloor-1)/2}} \le CK_{\delta}\sqrt{\log(r+2)}.\label{eqn:Rr-Deltaj-ub2}
    \end{align}
    Summing \eqref{eqn:Rr-Deltaj-ub1} and \eqref{eqn:Rr-Deltaj-ub2} over $j$, we obtain
    \begin{align}
        &\sum_{j=0}^{\lfloor\log_2(r/3)\rfloor-1}\left(\frac{|\Xi_{j+1}(W)-\Xi_j(W)|}{2^{j/2}}+\frac{|\Xi_{j+1}(W)-\Xi_j(W)|^2}{2^{j}}\right) \le CK_{\delta}^2\log^2(r+2),\label{eqn:Rr-hard-ub}
    \end{align}
    and combining with \eqref{eqn:Rr-easy-ub} we have on the event $\cG_{N,\delta}$, there is an absolute constant $C<\infty$ such that
    \begin{equation}
        R_r \le CK_{\delta}^2\log^2(r+2),\qquad\forall r\in[6,N].\label{eqn:Rr-ub}
    \end{equation}

    Plugging \eqref{eqn:R*-ub} and \eqref{eqn:Rr-ub} into \eqref{eqn:adapted10.2} and \eqref{eqn:short r} respectively, recalling \eqref{eqn:typical-W-path}, proves \Cref{claim:short-time-uniform-LB}, which completes the proof of \Cref{claim:remove-deterministic} and hence of \Cref{lem:Kriech1}.
    \end{proof}

The next lemma compares the quenched ballot-type probabilities between a random walk and a Gaussian random walk. A similar comparison result in a different setting can be found in \cite[Corollary 5.23]{Kriechbaum:space}. 

\begin{lemma}\label{lem:Kriech2}
    Let $S,W$ be independent centered random walks with iid increments and $S_0=W_0=0$. 
    Let $(\mathcal F_n)_{n\ge0}$ be a filtration such that $W$ is adapted, $W_n-W_{n-1}$ is independent of $\mathcal F_{n-1}$ for every $n\ge1$, and $S$ is independent of $\mathcal F_\infty$.
    Assume $\sigma^2=\Var(S_1)>0$, $\Var(W_1)=1$, and $\E e^{\eta|S_1|}, \E e^{\eta |W_1|}<\infty$ for some $\eta>0$. 
    Let $(r_n)_{n\ge 1}$ satisfy that $\lim_{n\to\infty}\frac{r_n}{\log\log n}=\infty$ and $\lim_{n\to\infty}\frac{r_n}{\log n}=0$.
    Fix $\alpha\in(0,\frac16]$. For $A\in\R$,
    and $1\le j\le n$ let
    \[
    h_j^{(n,A)}:=\min(j, n-j)^\alpha - \frac jnA\,.
    \]
    For any fixed $\delta>0$ and $M<\infty$, there exist $n_{\delta,M}\in\N$ such that for all $n\ge n_{\delta,M}$, with probability at least $1-\delta$ the following event holds:
    \begin{align}
        &\qquad\qquad p_{n,S}^{x;y,z}(W-h^{(n,A)}) = e^{O_{\delta}(\log\log n)}p_{n,B}^{10;10,1}(W/\sigma),\notag\\
        &\; \forall \;
        |A|\le M\log n\,,\; x\in[r_n,2\log^2 n]\,, \;y\in[0,\log^2 n]\,, \;z\in[r_n, 3\log n].\label{event:Kriech2}
    \end{align}
\end{lemma}

\begin{proof}
    Throughout this proof, implicit constants may depend on the laws of $S_1$ and $W_1$ and on $M$, but are uniform in $n,x,y,z$, and $A$ over the stated ranges. Dependence on $\delta$ will be indicated explicitly. Implicit constants may vary from line to line.
    
    We first prove the upper bound in \eqref{event:Kriech2}.
    Fix $q>0$, to be chosen below. By the KMT coupling \cite{komlos1975approximation,komlos1976approximation}, there is a constant $C_q<\infty$ and a coupling of $S$ with a standard Brownian motion $(B_t)_{t\ge0}$ such that, for $G_j:=\sigma B_j$ and all sufficiently large $n$,
    \begin{equation}\label{eqn:KMT-upper-event}
        \ee_n:=C_q\log n,\qquad \cK_n:=\left\{\max_{j\le n}|S_j-G_j|\le\ee_n\right\}, \qquad \P(\cK_n)\ge1-n^{-q}.
    \end{equation}
    Let $\ell_j^{(n,A)}:=-\frac jnA$ be the linear part of $h_j^{(n,A)}$. On $\cK_n$ one has the deterministic inclusion
    \begin{align}
        &\left\{\max_{j\le n}S_j-W_j+h_j^{(n,A)}\le x,\, S_n-W_n+h_n^{(n,A)}-x\in[-y-z,-y]\right\}\notag\\
        &\subseteq
        \left\{\max_{j\le n}G_j-W_j+h_j^{(n,A)}\le x+\ee_n,\, G_n-W_n+h_n^{(n,A)}-(x+\ee_n)\in[-y-z-2\ee_n,-y]\right\}\notag\\
        &\subseteq \left\{\max_{j\le n}G_j-W_j+\ell_j^{(n,A)}\le x+\ee_n,\, G_n-W_n+\ell_n^{(n,A)}-(x+\ee_n)\in[-y-z-2\ee_n,-y]\right\}.
        \label{eqn:KMT-upper-inclusion}
    \end{align}
    In the second inclusion we used $\ell_j^{(n,A)}\le h_j^{(n,A)}$ and $\ell_n^{(n,A)}=h_n^{(n,A)}$. 
    Since the coupling is independent of $W$, \eqref{eqn:KMT-upper-event} and \eqref{eqn:KMT-upper-inclusion} give, pathwise in $W$,
    \begin{equation}\label{eqn:KMT-upper-first}
        p_{n,S}^{x;y,z}(W-h^{(n,A)}) \le p_{n,G}^{x+\ee_n;y,z+2\ee_n}(W-\ell^{(n,A)})+n^{-q}.
    \end{equation}

    We next remove $\ell^{(n,A)}$ using the Girsanov change of measure.
    Set
    \[
        x':=\frac{x+\ee_n}{\sigma},\qquad y':=\frac{y}{\sigma},\qquad z':=\frac{z+2\ee_n}{\sigma},
    \]
    and let $\P_\mu$ denote the law under which $G$ has independent $N(\mu,\sigma^2)$ increments. Then
    \begin{equation}\label{eqn:affine-drift-identity}
    \begin{split}
        p_{n,G}^{\sigma x';\sigma y',\sigma z'}(W-\ell^{(n,A)}) = \P_{-A/n}\Big(&\max_{j\le n}G_j-W_j\le\sigma x',\,G_n-W_n-\sigma x'\in[-\sigma(y'+z'),-\sigma y']\Big).
    \end{split}
    \end{equation}
    On $\sigma(G_1,\ldots,G_n)$,
    \begin{equation}\label{eqn:affine-RN-density}
        \frac{\d\P_{-A/n}}{\d\P_0} = \exp\left(-\frac{A}{\sigma^2n}G_n-\frac{A^2}{2\sigma^2n}\right).
    \end{equation}
    On the event in \eqref{eqn:affine-drift-identity}, one has
    \[
        G_n=W_n+\sigma(x'-u)
    \]
    for some $u\in[y',y'+z']$. Under $\P_0$, dividing the event by $\sigma$ gives
    \[
        \P_{0}\Big(\max_{j\le n}G_j-W_j\le\sigma x',\,G_n-W_n-\sigma x'\in[-\sigma(y'+z'),-\sigma y']\Big) = p_{n,B}^{x';y',z'}(W/\sigma).
    \]
    Consequently,
    \begin{align}
    &p_{n,G}^{x+\ee_n;y,z+2\ee_n}(W-\ell^{(n,A)}) \quad\le \exp\left\{\frac{|A|}{\sigma^2n}\bigl(|W_n|+\sigma(x'+y'+z')\bigr)+\frac{A^2}{2\sigma^2n}\right\}p_{n,B}^{x';y',z'}(W/\sigma).
    \end{align}
    On the event
    \begin{equation}\label{eqn:Wn-upper-good}
        |W_n|\le\sqrt n\log n,
    \end{equation}
    whose probability tends to one, we therefore have
    \begin{equation}\label{eqn:remove-affine-upper}
        p_{n,G}^{x+\ee_n;y,z+2\ee_n}(W-\ell^{(n,A)}) \le e^{o(1)}p_{n,B}^{x';y',z'}(W/\sigma)
    \end{equation}
    uniformly over all parameters in the lemma.

    It remains to compare the last probability with the target probability $p_{n,B}^{10;10,1}(W/\sigma)$. 
    By the union bound,
    \begin{equation}\label{eqn:long-endpoint-upper-sum}
        p_{n,B}^{x';y',z'}(W/\sigma) \le \sum_{k=0}^{\lceil z'\rceil}p_{n,B}^{x';y'+k,1}(W/\sigma).
    \end{equation}
    For every $k$, monotonicity in the barrier gives
    \[
        p_{n,B}^{x';y'+k,1}(W/\sigma) \le p_{n,B}^{x'+4;y'+k+4,1}(W/\sigma).
    \]
    For all sufficiently large $n$, uniformly over the parameter ranges in the statement and all $0\le k\le\lceil z'\rceil$,
    \[
        4\le x'+4\le C_{q}\log^2n<\log^5n,\quad 4\le y'+k+4\le C_{q}\log^2n<\log^5n,\quad \lceil z'\rceil = O_{q}(\log n).
    \]
    We apply the upper bound in \Cref{lem:Kriech1} to the last probability, summing over $k$ in \eqref{eqn:long-endpoint-upper-sum} and enlarging $C_\delta$ to absorb the factor
    \[
        1+\lceil z'\rceil=O_{q}(\log n),
    \]
    we obtain, on an event of probability at least $1-\delta/4$,
    \begin{equation}\label{eqn:Gaussian-upper-beta}
        p_{n,B}^{x';y',z'}(W/\sigma) \le (\log n)^{C_\delta}p_{n,B}^{10;10,1}(W/\sigma)
    \end{equation}
    uniformly over all the parameters under consideration.
    
    Finally, \Cref{prop:LB-wiggle}, used with $\wt W=W$ and $(x,y,z)=(10,10,1)$, implies that
    \begin{equation}\label{eqn:beta-polynomial-lower-for-KMT}
        p_{n,B}^{10;10,1}(W/\sigma)\ge n^{-C_*}\qquad\text{a.a.s.}
    \end{equation}
    for a deterministic $C_*<\infty$. Choose $q>C_*+1$.
    Combining \eqref{eqn:KMT-upper-first}, \eqref{eqn:remove-affine-upper}, \eqref{eqn:Gaussian-upper-beta}, and \eqref{eqn:beta-polynomial-lower-for-KMT} yields
    \begin{equation}\label{eqn:KMT-upper}
        p_{n,S}^{x;y,z}(W-h^{(n,A)}) \le (\log n)^{C_\delta}p_{n,B}^{10;10,1}(W/\sigma)
    \end{equation}
    uniformly in the asserted parameter ranges, on an event of probability at least $1-\delta/2$ for all sufficiently large $n$.

    We now turn to the lower bound. 
    By monotonicity, we may assume $\alpha=1/6$.
    Set
    \[
        m:=\left\lceil C_{\mathrm{cut}}(\log n)^4\right\rceil, \quad \cI:=[-3(\log n)^2,-(\log n)^2], 
    \]
    with $C_{\mathrm{cut}}>0$ being a large constant chosen later, depending only on $\sigma$.
    By the Markov property, we have
    \begin{align}
        p_{n,S}^{x;y,z}(W-h) \ge p_{m,S}^{x;\cI}(W-h)\inf_{u\in-\cI}p_{m,n-m;S}^{u;\cI}(W-h) \inf_{u\in-\cI} p_{n-m,n;S}^{u;[-y-z,-y]}(W-h).
        \label{eqn:lower-three-block-Markov}
    \end{align}
    We bound each of the three terms from below as follows.
    
    \begin{claim}\label{claim:two edges LB}
        For every $\delta>0$, there exist $C_\delta<\infty$ and $n_\delta\in\N$ such that, for every $n\ge n_\delta$, with probability at least $1-\delta$, the following holds uniformly over the stated parameter ranges:
        \[
        \min\{p_{m,S}^{x;\cI}(W-h),\, \inf_{u\in-\cI} p_{n-m,n;S}^{u;[-y-z,-y]}(W-h)\} \ge m^{-C_{\delta}}.
        \]
    \end{claim}
    
    \begin{claim}\label{claim:middle LB}
        For every $\delta>0$, there exist $C_\delta<\infty$ and $n_\delta\in\N$ such that, for every $n\ge n_\delta$, with probability at least $1-\delta$, the following estimate holds uniformly over $|A|\le M\log n$:
        \begin{equation}\label{eqn:middle comparison}
        \inf_{u\in-\cI}p_{m,n-m;S}^{u;\cI}(W-h) \ge (\log n)^{-C_{\delta}}
        \overline p_{n,B}^{10;10,1}((W-h)/\sigma).
        \end{equation}
    \end{claim}

    \begin{claim}\label{claim:middle compare to full}
        For every $\delta>0$, there exist $C_\delta<\infty$ and $n_\delta\in\N$ such that for every $n\ge n_\delta$, with probability at least $1-\delta$, uniformly over $|A|\le M\log n$:
        \[
        \overline p_{n,B}^{10;10,1}((W-h)/\sigma) \ge (\log n)^{-C_\delta}p_{n,B}^{10;10,1}(W/\sigma).
        \]
    \end{claim}
    
    Applying each of the three claims with $\delta/6$ in place of $\delta$, and using \eqref{eqn:lower-three-block-Markov} and $m=\Theta(\log^4 n)$, proves the required lower bound with probability at least $1-\delta/2$, after increasing $C_\delta$.
    Combining with \eqref{eqn:KMT-upper}, we conclude the proof of \Cref{lem:Kriech2}.
    \medskip
    
    It remains to prove \Cref{claim:two edges LB}, \Cref{claim:middle LB}, and \Cref{claim:middle compare to full}.
    Recall the KMT coupling \Cref{lem:KMT} for proofs of all three claims.
    \medskip

    \noindent\emph{Proof of \Cref{claim:two edges LB}.}
    Take $N=m$ and $\lam=C_{\mathrm{KMT}}\log m$ in \eqref{eqn:KMT}, where $C_{\mathrm{KMT}}>0$ is a large constant to be determined later.
    Let $d_{1}:=(c_0+C_{\mathrm{KMT}})\log m = \Theta(\log\log n) = o(r_n)$ be the total KMT coupling loss between $S$ and $G$.
    Then for sufficiently large $n$, we have
    \begin{equation}\label{eqn:estimate-I-KMT}
    p_{m,S}^{x;\cI}(W-h) \ge p_{m,G}^{x-d_{1};\log^2 n, \sigma}(W-h)-c_1m^{-c_2C_{\mathrm{KMT}}} \ge \overline p_{m,G}^{x-d_{1};\log^2 n, \sigma}(W-h)-c_1m^{-c_2C_{\mathrm{KMT}}}.
    \end{equation}
    We bound the last quenched probability from below by \cite[Theorem 10.2]{Kriechbaum:time}. To match their assumptions, we make the following reductions.
    Firstly, we rewrite
    \[
    \overline p_{m,G}^{x-d_{1};\log^2 n, \sigma}(W-h) = \overline p_{m,B}^{\frac{x-d_{1}}{\sigma};\frac{\log^2 n}{\sigma}, 1}((W-h)/\sigma)
    \]
    so that the Brownian motion is standard.
    
    Next, we introduce $\wt h_s:=(h_s-\frac{h_m}{m}s)/\sigma$ for $s\in[0,m]$ (so $\wt h_s = (s^{1/6}-m^{-5/6}s)/\sigma$ at integer times and linearly interpolated in between), as a tilted version of the deterministic barrier $h/\sigma$. By definition, we have $\wt h_0=\wt h_m=0$ satisfying the boundary convention at the beginning of \cite[Section 8]{Kriechbaum:time}, and $\wt h$ is ``nice" in the sense of \cite[Definition 11.5]{Kriechbaum:time} with any exponent $\eta\in(1/6,1/2)$ and for any sufficiently large $n$.
    
    Note the slope of the removed line is $h_m/m\sigma = \sigma^{-1}(m^{-5/6}-A/n)=O(m^{-5/6})$ uniformly for $|A|\le M\log n$.
    By Girsanov change-of-measure, with a computation similarly as in \eqref{eqn:remove-affine-upper}, the logarithm of the Radon--Nikodym factor on the relevant barrier event is of the size
    \[
    O_{\sigma}\Big(|\frac{h_m}{m}|\big(|W_m|+x+\log^2 n+|h_m|\big)\Big),
    \]
    which is uniformly $o(1)$ on the event $\{|W_m|\le \sqrt{m}\log m\}$ (whose probability goes to one) over the parameter range in the statement.
    Therefore, on the event $\{|W_m|\le \sqrt{m}\log m\}$, we have 
    \begin{align}
        \overline p_{m,G}^{x-d_{1};\log^2 n, \sigma}(W-h) &\ge e^{-o(1)}\overline p_{m,B}^{\frac{x-d_{1}}{\sigma};\frac{\log^2 n}{\sigma},1}(W/\sigma-\wt h) \label{eqn:lem6.4-I-reduction}
    \end{align}
    Fix $C_{\mathrm{cut}}>100\sigma^{-2}$ to ensure the parameters are uniformly in the range of \Cref{lem:Gaussian-crude-LB} for all sufficiently large $n$.
    Applying \Cref{lem:Gaussian-crude-LB} and then choosing $C_{\mathrm{KMT}}$ sufficiently large, we obtain the desired result: with probability at least $1-\delta/2$, we have
    \begin{equation}\label{eqn:estimate-I}
        \inf_{x\in[r_n,2\log^2 n]} p_{m,S}^{x;\cI}(W-h) \ge m^{-C_{\delta}}
    \end{equation}
    for some $C_{\delta}<\infty$.
    
    The argument leading to \eqref{eqn:estimate-I} also gives a uniform lower bound for $\inf_{u\in-\cI}p_{n-m,n;S}^{u;[-y-z,-y]}(W-h)$. We therefore indicate only the necessary modifications.
    
    Firstly, with the same KMT coupling parameters, for all $y, z$ in the stated ranges and all $n$ sufficiently large (recall $d_1=o(r_n)$ while $z\ge r_n$), the analogue of \eqref{eqn:estimate-I-KMT} is
    \begin{align}
        \inf_{u\in-\cI}p_{n-m,n;S}^{u;[-y-z,-y]}(W-h) &\ge \inf_{u\in-\cI}p_{n-m,n;G}^{u-d_1;[-y-\sigma,-y]}(W-h)-c_1m^{-c_2C_{\mathrm{KMT}}} \notag\\
        &\ge \inf_{u\in-\cI}\overline p_{n-m,n;G}^{u-d_1;[-y-\sigma,-y]}(W-h)-c_1m^{-c_2C_{\mathrm{KMT}}} \notag\\
        &= \inf_{u\in-\cI}\overline p_{n-m,n;B}^{\frac{u-d_1}{\sigma};[-\frac{y}{\sigma}-1,-\frac{y}{\sigma}]}((W-h)/\sigma)-c_1m^{-c_2C_{\mathrm{KMT}}}, \label{eqn:estimate-III-KMT}
    \end{align}
    where, similarly, when we pass to the continuous time barrier event, we linearly interpolate $W$ and $h$ between integer times. The rest follows again from \Cref{lem:Gaussian-crude-LB} and enlarging $C_{\mathrm{KMT}}$ when necessary. Together with \eqref{eqn:estimate-I} this proves \Cref{claim:two edges LB}.
    \medskip

    \noindent\emph{Proof of \Cref{claim:middle LB}.}
    Take $N=n-2m$ and $\lambda = C_{\mathrm{KMT}}\log N$ in \eqref{eqn:KMT}, and let $d_2 := (c_0+C_{\mathrm{KMT}})\log N = \Theta(\log n)$ be the KMT coupling loss. Then for sufficiently large $n$, by \eqref{eqn:KMT}, dividing by $\sigma$, restriction to a unit endpoint interval contained in $\sigma^{-1}[-3\log^2 n+2d_2, -\log^2 n]$, and passing to continuous time, we have
    \begin{equation}
        \inf_{u\in-\cI}p_{m,n-m;S}^{u;\cI}(W-h) \ge \inf_{u\in-\cI}\overline p_{m,n-m;B}^{\frac{u-d_2}{\sigma};\frac{\log^2 n}{\sigma},1}((W-h)/\sigma) - c_1 N^{-c_2C_{\mathrm{KMT}}}.\label{eqn:estimate-II-KMT}
    \end{equation}
    As before, in passing to the continuous-time barrier event, we linearly interpolate $W$ and $h$ between integer times.

    As in the argument leading to \eqref{eqn:estimate-I}, applying \Cref{lem:Gaussian-crude-LB} and enlarging $C_{\mathrm{KMT}}$ if necessary gives, for all sufficiently large $n$, on an event $\cG_{\delta}^{(0)}$ with probability at least $1-\delta/10$,
    \begin{equation}
        \inf_{u\in-\cI}p_{m,n-m;S}^{u;\cI}(W-h) \ge \frac12\inf_{u\in-\cI}\overline p_{m,n-m;B}^{\frac{u-d_2}{\sigma};\frac{\log^2 n}{\sigma},1}((W-h)/\sigma).\label{eqn:II-after-KMT}
    \end{equation}

    Next, we compare the RHS of \eqref{eqn:II-after-KMT} with the RHS of \eqref{eqn:middle comparison}.
    We start from shortening the time interval in the latter quenched probability to match that in the former.
    With $W$ and $h$ linearly interpolated between integer times, define
    \[
        Z_{a,b}(q):=q+G_b-G_a-(W_b-W_a)+(h_b-h_a),\qquad 0\le a\le b\le n.
    \]
    Note $G=\sigma B$ is a rescaled Brownian motion and is not linearly interpolated.
    Let $\cJ:=[-11\sigma,-10\sigma]$.
    By Markov property and monotonicity, we have
    \begin{align*}
        \overline p_{n,B}^{10;10,1}((W-h)/\sigma) = \overline p_{n,G}^{10\sigma;\cJ}(W-h) \le \P\left(Z_{0,m}(-10\sigma)\le -\log^3 n\middle| W\right) + \sup_{0\le q\le \log^3 n}\overline p_{m,n;G}^{q;\cJ}(W-h).
    \end{align*}
    Among the stated parameter ranges, we uniformly have $h_m = m^{1/6}-\frac{Am}{n} = O((\log n)^{2/3})$ while $h_0=W_0=0$. So, on the event 
    \[
    \cG_{1} := \{|W_m|\le \sqrt{m}\log m\},\qquad \lim_{n\to\infty}\P(\cG_1)=1,
    \]
    we have
    \[
    \bone(\cG_1)\P\left(Z_{0,m}(-10\sigma)\le -\log^3 n\middle| W\right) \le \P(G_m \le -\log^3 n/2) \le e^{-c\log^2 n}
    \]
    for some $c>0$ when $n$ is sufficiently large.
    Meanwhile, by \Cref{lem:Gaussian-crude-LB}, for any $\delta>0$ there is a constant $C_{\delta}<\infty$ and an event $\cG_{\delta}^{(1)}$ with $\P(\cG_{\delta}^{(1)})>1-\delta/10$, on which we have
    \[
    \bone(\cG_{\delta}^{(1)})\sup_{0\le q\le \log^3 n}\overline p_{m,n;G}^{q;\cJ}(W-h) \ge \bone(\cG_{\delta}^{(1)})\overline p_{m,n;G}^{\sigma;\cJ}(W-h) \ge n^{-C_{\delta}}\bone(\cG_{\delta}^{(1)})
    \]
    for sufficiently large $n$ (depending on $\delta$).
    Hence when $n$ is sufficiently large, we obtain
    \begin{align}
        \bone(\cG_1\cap\cG_{\delta}^{(1)})\overline p_{n,G}^{10\sigma;\cJ}(W-h) \le 2\sup_{0\le q\le \log^3 n}\overline p_{m,n;G}^{q;\cJ}(W-h).\label{eqn:remove first m}
    \end{align}

    To remove the last block $[n-m,n]$, we introduce the following time-reversed processes.
    Let $\rev{G}_{s}:=G_{n-s}-G_n$ for any $s\in[0,n-m]$, and similarly $\rev W, \rev h$. Define $\rev Z_{a,b}(q)$ similarly as $Z_{a,b}(q)$ with each process replaced by its reversed counterpart. 
    Then
    \begin{equation}
    Z_{m,s}(q) = \rev Z_{0,n-s}\bigl(Z_{m,n}(q)\bigr),\,s\in[m,n]\,;\quad \rev Z_{m,s}(q') = Z_{m,n-s}\bigl(\rev{Z}_{m,n-m}(q')\bigr),\, s\in[m,n-m]\label{eqn: Z and Z reverse}
    \end{equation}
    for any $q,q'\in\R$.
    Using the above equation and $Z_{m,n}(-q)\in \cJ$ in the barrier event defining $\overline p_{m,n;G}^{q;\cJ}(W-h)$, we have
    \begin{align}
    \sup_{0\le q\le \log^3 n}\overline p_{m,n;G}^{q;\cJ}(W-h) \le \sup_{0\le q\le \log^3 n}\overline p_{n-m,\rev G}^{11\sigma; q,\sigma}(\rev W-\rev h).\label{eqn:reverse time 1}
    \end{align}
    By a similar argument towards \eqref{eqn:remove first m}, noticing that $\rev G, \rev W, \rev h$ on $[0,n-m]$ also satisfies the assumptions in \Cref{lem:Gaussian-crude-LB}, with $\cG_{2}:=\{|\rev W_m|\le\sqrt{m}\log m\} = \{|W_{n}-W_{n-m}|\le\sqrt{m}\log m\}$, for any $\delta>0$ there is a constant $C_{\delta}<\infty$ and an event $\cG_{\delta}^{(2)}$ with $\liminf_{n\to\infty}\P(\cG_{\delta}^{(2)})\ge 1-\delta/10$, such that
    \begin{align}
        \bone(\cG_2\cap\cG_{\delta}^{(2)})\sup_{0\le q\le \log^3 n}\overline p_{n-m,\rev G}^{11\sigma; q,\sigma}(\rev W-\rev h) \le 2\sup_{0\le q,q'\le \log^3 n}\overline p_{m,n-m;\rev G}^{q'; q,\sigma}(\rev W-\rev h).\label{eqn:remove last m}
    \end{align}
    Using again \eqref{eqn: Z and Z reverse} and $\rev Z_{m,n-m}(-q')\in[-q-\sigma, -q]$ in the barrier event defining $\overline p_{m,n-m;\rev G}^{q'; q,\sigma}(\rev W-\rev h)$, we have
    \begin{align}
        \overline p_{m,n-m;\rev G}^{q'; q,\sigma}(\rev W-\rev h) \le \overline p_{m,n-m; G}^{q+\sigma; q',\sigma}(W-h) = \overline p_{m,n-m; B}^{\frac{q}{\sigma}+1;\frac{q'}{\sigma},1}((W-h)/\sigma).\label{eqn:reverse time 2}
    \end{align}

    Finally, by monotonicity we have
    \[
    \overline p_{m,n-m; B}^{\frac{q}{\sigma}+1;\frac{q'}{\sigma},1}((W-h)/\sigma) \le \overline p_{m,n-m; B}^{\frac{q}{\sigma}+5;\frac{q'}{\sigma}+4,1}((W-h)/\sigma).
    \]
    Applying \Cref{lem:Kriech1} to the RHS of the last expression gives, for any $\delta>0$ there is a $C_{\delta}<\infty$ and an event $\cG_{\delta}^{(3)}$ with $\liminf_{n\to\infty}\P(\cG_{\delta}^{(3)})\ge 1-\delta/10$, such that
    \begin{align}
        \sup_{0\le q,q'\le \log^3 n}\overline p_{m,n-m; B}^{\frac{q}{\sigma}+1;\frac{q'}{\sigma},1}((W-h)/\sigma) \le (\log n)^{C_{\delta}}\inf_{u\in-\cI}\overline p_{m,n-m;B}^{\frac{u-d_2}{\sigma};\frac{\log^2 n}{\sigma},1}((W-h)/\sigma).\label{eqn:comparison m to n-m}
    \end{align}
    Since $\liminf_{n\to\infty}\P\Big(\bigcap_{i=1,2}\cG_i\cap\bigcap_{i=0}^{3}\cG_{\delta}^{(i)}\Big)\ge 1-\delta$, combining \eqref{eqn:II-after-KMT}, \eqref{eqn:remove first m}, \eqref{eqn:reverse time 1}, \eqref{eqn:remove last m}, \eqref{eqn:reverse time 2}, and \eqref{eqn:comparison m to n-m} finishes the proof of \Cref{claim:middle LB}.
    \medskip

    \noindent\emph{Proof of \Cref{claim:middle compare to full}}
    This claim is an application of \Cref{lem:Kriech1}.
    Indeed, we apply \Cref{lem:Kriech1} on $[0,n]$ to the environmental walk $W/\sigma$ and the additional perturbative curve $h/\sigma$.
    The linear component of $h/\sigma$ has slope $-A/(\sigma n)$, so the total absolute change along the line is at most $(M/\sigma)\log n$.
    The remaining curve is the linear interpolation of $\sigma^{-1}\min(j,n-j)^{1/6}$, and when $n$ is sufficiently large it satisfies
    \[
        0\le\frac{h_s+As/n}{\sigma}\le 10\frac{((1+s)\wedge(1+n-s))^{1/6}-1}{\sigma},\qquad 0\le s\le n.
    \]
    The inequality holds at integer times and extends between them by concavity of the expression on the right.
    Thus \eqref{eqn:Gaussian-comparison-curve} holds with fixed constants, uniformly over $|A|\le M\log n$.
    Taking $(x,y,z)=(10,10,1)$ in the continuous-time part of \Cref{lem:Kriech1} gives, with probability at least $1-\delta$, there is a constant $C_{\delta}<\infty$ such that for all sufficiently large $n$,
    \[
        \overline p_{n,B}^{10;10,1}((W-h)/\sigma)\ge (\log n)^{-C_{\delta}}p_{n,B}^{10;10,1}(W/\sigma)
    \]
    uniformly over $|A|\le M\log n$.
    This concludes the proof of \Cref{claim:middle compare to full}.    
\end{proof}

The next several lemmas address the stability of quenched barrier probabilities under small perturbations.

\begin{lemma}\label{lem:Kriech3}
    Recall $S, W, (\mathcal F_j)_{j\ge 0}, \sigma, r_n, h^{(n,A)}$ as in \Cref{lem:Kriech2}.
    Let $\wt W$ be an $(\mathcal F_j)$-adapted process satisfying $\wt W_0=0$ and
    \[
    \max_{1\le j\le n}|W_j-\wt W_j|\le\frac1{100}\qquad\text{a.s.}
    \]
        
    For any fixed $\delta, M>0$ there exists $n_{\delta}\in\N$ such that for all $n\ge n_{\delta}$, with probability at least $1-\delta$ we have uniformly for $|A|\le M\log n, x\in[2r_n,\log^2 n], y\in[0, \log^2 n], z\in[2r_n, 2\log n]$:
    \begin{equation}\label{eqn:rounded comparison}
        p_{n,S}^{x;y,z}(\wt W-h^{(n,A)}) = e^{O_{\delta}(\log\log n)}p_{n,B}^{10;10,1}(\wt W/\sigma).
    \end{equation}
\end{lemma}
\begin{proof}
    Fix $\delta, M>0$ and let $\epsilon_0:=\frac{1}{100}$. Since $|\wt W_j-W_j|\le\epsilon_0$, monotonicity gives, pathwise,
    \begin{align}
        p_{n,S}^{x-\epsilon_0;y,z-2\epsilon_0}(W-h^{(n,A)}) &\le p_{n,S}^{x;y,z}(\wt W-h^{(n,A)}) \le p_{n,S}^{x+\epsilon_0;y,z+2\epsilon_0}(W-h^{(n,A)}).\label{eqn:perturb-curve}
    \end{align}
    For all sufficiently large $n$, uniformly over the stated ranges of $x,y,z$, the parameters $x\pm\epsilon_0,y,z\pm2\epsilon_0$ lie in the corresponding ranges of \Cref{lem:Kriech2}.
    Applying this lemma with error probability $\delta/2$, we obtain for every sufficiently large $n$, an event of probability at least $1-\delta/2$, on which
    \begin{equation}\label{eqn:rounded-comparison-W}
        p_{n,S}^{x;y,z}(\wt W-h^{(n,A)})=e^{O_\delta(\log\log n)}p_{n,B}^{10;10,1}(W/\sigma)
    \end{equation}
    uniformly over the stated ranges for $A, x, y, z$.

    It remains to replace $p_{n,B}^{10;10,1}(W/\sigma)$ by $p_{n,B}^{10;10,1}(\wt W/\sigma)$ in the above equation.
    Let $f_j:=(W_j-\wt W_j)/\sigma$ for $0\le j\le n$, and linearly interpolate in between integer times. Then
    \[
    p_{n,B}^{10;10,1}(\wt W/\sigma) = p_{n,B}^{10;10,1}(W/\sigma-f).
    \]
    Since $f_0=0$, and $|f_j|\le \epsilon_0/\sigma$, we have
    \[
    |\frac{f_n}{n}| \le \frac{\epsilon_0}{\sigma n},\qquad |f_{s}-f_0-\frac{f_n}{n}s| \le \frac{10\,\epsilon_0}{\sigma(2^{1/6}-1)}\left(((1+s)\wedge(1+n-s))^{1/6}-1\right), \; s\in[0,n],
    \]
    where the second inequality follows from boundedness of $f$ on integer values and concavity of the function on the right. 
    Thus $f$ satisfies the curve assumptions of \Cref{lem:Kriech1}. Since that lemma is uniform over all curves satisfying those assumptions, it applies also to this quenched random curve.
    By \Cref{lem:Kriech1}, when $n$ is sufficiently large there is an event of probability at least $1-\delta/2$, on which we have
    \begin{equation}\label{eqn:beta-rounded-comparison}
        p_{n,B}^{10;10,1}(\wt W/\sigma) = e^{O_\delta(\log\log n)}p_{n,B}^{10;10,1}(W/\sigma).
    \end{equation}
    Combining \eqref{eqn:rounded-comparison-W} and \eqref{eqn:beta-rounded-comparison} proves the claim.
\end{proof}

\begin{lemma}[Fine stability of quenched barrier probabilities]\label{lem:stab-gen}
    Recall $S, W, (\mathcal F_j)_{j\ge0}, r_n$ as in \Cref{lem:Kriech2}. 
    Further assume $S_1$ has the distribution of $f(\upsilon)$, where $\upsilon$ is uniform in $\R/\Z$ and $f$ is a smooth 1-periodic function satisfying \Cref{assu:f}.
    Fix $M<\infty$, $\alpha\in(0,1/6]$, and $0\le v<\esssup S_1$.
    Let $\wt W=\wt W^{(n)}$ be an $(\mathcal F_j)$-adapted process with $\wt W_0=0$ and $\max_{0\le j\le n}|W_j-\wt W_j|\le \eps:=e^{-r_n}$ a.s.
    For $h^{(n,A)}$ as in \Cref{lem:Kriech2}, set
    \[
    Z_j:=S_j-vj-\wt W_j+h_j^{(n,A)}.
    \]
    There is a constant $c>0$, depending only on $v$ and the laws of $S_1,W_1$, such that, uniformly for $|A|\le M\log n$ and $x, y\in[r_n,\log n]$,
    \[
    \P\Big(\max_{1\le j\le n}Z_j\le x+\eps,\ Z_n\ge x-y-\eps\mid\mathcal F_\infty\Big)
    \le(1+\eps^c)\P\Big(\max_{1\le j\le n}Z_j\le x,\ Z_n\ge x-y\mid\mathcal F_\infty\Big)\qquad a.a.s.
    \]
\end{lemma}

\begin{proof}
Let $\cE_{t}:=\{\max_{1\le j\le n}Z_j\le x+t,\,Z_n\ge x-y-t\}, \forall t\ge 0$. Then our goal is to show that there is a constant $c>0$ such that
\begin{equation}\label{eq:stab-gen}
    \P\left(\cE_{\eps}\mid\mathcal F_{\infty}\right) - \P\left(\cE_{0}\mid\mathcal F_{\infty}\right)\le \eps^c\P\left(\cE_{0}\mid\mathcal F_{\infty}\right)\qquad a.a.s.
\end{equation}
in the parameter ranges as stated in the lemma.

Let $z_j:=Z_j-S_1, \forall j\ge 1$ and $\mathcal H_n:=\sigma(S_j-S_1:2\le j\le n)$. Then we have
\[
    \cE_0 = \{L\le S_1 \le U\},\qquad \cE_\eps = \{L-\eps\le S_1\le U+\eps\},
\]
where an interval with its lower endpoint larger than its upper endpoint is considered empty, and
\[
    L:= x-y-z_n,\qquad U:= x-\max_{1\le j\le n}z_j,\qquad L,U\in \sigma(\mathcal F_{\infty}, \mathcal H_n).
\]
Let
\[
    \cG_{\eps}:=\left\{\P\left(\cE_\eps\middle| \mathcal F_{\infty}, \mathcal H_n\right)>0\right\}\in\sigma(\mathcal F_{\infty}, \mathcal H_n)
\]
denote the event of $(z_j)_{1\le j\le n}$ on which $\cE_{\eps}$ is possible. Since $f$ is smooth and 1-periodic, $S_1$ is bounded. Let $D:=\esssup S_1-\essinf S_1 = \diam(\supp(S_1))$. The boundedness of $S_1$ then implies
\[
    \cG_{\eps}\subset \cE_{D+\eps}\,,\qquad a.s.
\]
By the H\"older continuity of the distribution function of $S_1$ (\Cref{assu:f}) and independence between $S_1$ and $\sigma(\mathcal F_{\infty}, \mathcal H_n)$, there are constants $0<C_f,\kappa<\infty$ such that almost surely
\begin{align}
    \P\left(\cE_{\eps}\mid\mathcal F_{\infty}\right)-\P\left(\cE_{0}\mid\mathcal F_{\infty}\right) &= \E\left[\P\left(\cE_{\eps}\mid\mathcal F_{\infty},\mathcal{H}_n\right)-\P\left(\cE_{0}\mid\mathcal F_{\infty},\mathcal{H}_n\right)\middle| \mathcal{F}_\infty\right] \notag\\
    &\le \E\left[\P\left(S_1\in [L-\eps, L)\cup(U,U+\eps]\mid\mathcal F_{\infty},\mathcal{H}_n\right)\bone(\cG_\eps)\middle|\mathcal F_{\infty}\right] \notag\\
    &\le C_f\eps^{\kappa}\P\left(\cE_{D+\eps}\middle|\mathcal F_\infty\right).\label{eq:stability}
\end{align}

Note that under our assumptions $(S_j)_{j\ge 0}\overset{law}{=}(X_j(0))_{j\ge 0}$, so we recall $\lambda$ in \eqref{def:lambda}. 
Fix a $v\in[0,\esssup S_1)$. Let $\theta\ge 0$ be such that $\lambda'(\theta)=v$, and define the Girsanov change-of-measure
\[
    \frac{d \P_v}{d \P}:=\exp(\theta S_n-n\lambda(\theta)).
\]
Under $\P_v$, the process $T_j:=S_j-vj,\forall j\ge 1$ is a centered random walk with iid bounded increments, independent from $\mathcal{F}_\infty$ which is unchanged. Especially, $T$ also satisfies the assumptions of $S$ in \Cref{lem:Kriech2} under $\P_v$.

Now we tilt the probability so that $T$ is centered.
Denote $R_n:=S_n-Z_n+x-y\in\mathcal F_\infty$, we have $S_n\ge R_n-(D+\eps)$ on $\cE_{D+\eps}$ and therefore
\begin{align}
    \P\left(\cE_{D+\eps}\mid\mathcal F_{\infty}\right) = \E_v\left[e^{n\lambda(\theta)-\theta S_n}\bone(\cE_{D+\eps})\middle| \mathcal{F}_\infty\right] \le e^{n\lambda(\theta)-\theta R_n+\theta(D+\eps)}\P_v\left(\cE_{D+\eps}\middle|\mathcal{F}_\infty\right).\label{eq:D+eps}
\end{align}
Let $\ell_n:=\sqrt{r_n\log\log n}$, then $\log\log n\ll \ell_n\ll r_n$.
For all sufficiently large $n$ so that $y>2\ell_n$, define 
\[
    \cE_{\mathrm{short}}:=\left\{\max_{1\le j\le n}Z_j\le x,\, Z_n\in[x-y,x-y+\ell_n]\right\}.
\]
Then
\begin{equation}
    \P\left(\cE_{0}\mid\mathcal F_{\infty}\right) \ge \P\left(\cE_{\mathrm{short}}\mid\mathcal F_{\infty}\right) = \E_v\left[e^{n\lambda(\theta)-\theta S_n}\bone(\cE_{\mathrm{short}})\middle| \mathcal{F}_\infty\right] \ge e^{n\lambda(\theta)-\theta R_n-\theta\ell_n}\P_v\left(\cE_{\mathrm{short}}\mid\mathcal F_{\infty}\right).\label{eq:short}
\end{equation}

Fix a $\delta>0$. Applying \Cref{lem:Kriech3} with $(T, \ell_n/10)$ in place of $(S,r_n)$ inside it, for all sufficiently large $n$ (depending on $\delta$), there is an event $\cG'_{\delta}\in\mathcal{F}_\infty$ with probability at least $1-\delta$ and a constant $C_\delta<\infty$, such that uniformly in the parameter ranges stated in \Cref{lem:stab-gen},
\[
    \P_v\left(\cE_{D+\eps}\middle|\mathcal{F}_\infty\right) \le (\log n)^{C_{\delta}}\P_v\left(\cE_{\mathrm{short}}\mid\mathcal F_{\infty}\right).
\]
Joint with \eqref{eq:D+eps} and \eqref{eq:short}, we have (constants may change)
\begin{equation}\label{eq:bound D+eps by 0}
    \bone(\cG'_\delta)\P\left(\cE_{D+\eps}\mid\mathcal F_{\infty}\right) \le \bone(\cG'_\delta)e^{\theta\ell_n}(\log n)^{C_\delta}\P\left(\cE_{0}\mid\mathcal F_{\infty}\right).
\end{equation}
Since $\log\log n\ll\ell_n\ll r_n$ and $\eps=e^{-r_n}$, by \eqref{eq:stability} and \eqref{eq:bound D+eps by 0}, for all sufficiently large $n$ (depending on $\delta$), uniformly in the stated parameter ranges we have
\[
    \bone(\cG'_\delta)\left(\P\left(\cE_{\eps}\mid\mathcal F_{\infty}\right)-\P\left(\cE_{0}\mid\mathcal F_{\infty}\right)\right) \le \bone(\cG'_\delta)\eps^{\kappa/2}\P\left(\cE_{0}\mid\mathcal F_{\infty}\right),
\]
i.e.
\[
    \liminf_{n\to\infty}\P\left[\P\left(\cE_{\eps}\mid\mathcal F_{\infty}\right)-\P\left(\cE_{0}\mid\mathcal F_{\infty}\right) \le \eps^{\kappa/2}\P\left(\cE_{0}\mid\mathcal F_{\infty}\right)\right] \ge 1-\delta.
\]
Since $\kappa$ is independent of $\delta$ and $\delta$ is arbitrary, taking $\delta\downarrow 0$ proves \eqref{eq:stab-gen} with $c=\kappa/2$.
\end{proof}

\appendix

\section{Logarithmic correlation}
\label{sec:covariance}

In this appendix we establish the logarithmic correlation structure \eqref{Cov-sim0} for the field $X_n$ of \eqref{def:field}, in a quantitative form.
This fact is irrelevant to the proof of the main theorem, but we include it for completeness (the decorrelation estimate needed for our proof of \Cref{thm:main} is encapsulated in \Cref{prop:decouple}). 

Let $\widehat{f}(k):=\int_{0}^{1}f(x)e^{-2\pi\ii k x}\d x$, $k\in\Z$ be the Fourier coefficients of $f$ (recall the notation $e(x):= \exp( 2\pi \ii x)$). In this section we abbreviate $\|f\|_2:= \|f\|_{L^2([0,1))}$.

\begin{prop}
    For any $s,t\in\R$, 
    \begin{align}   \label{Cov-sim}
        \Cov(X_n(s),X_n(t))
        &\sim  \|f\|_2^2\,\min\Big\{ n\,,\,  \log_{+}\frac1{|s-t|}\Big\}
        +O(1)\,.
    \end{align}  
    In fact we have the following more accurate asymptotic. 
    Denoting the Fourier projection
    \[
    \cF_{K}f(t):= \sum_{|k|<K}\wh f(k)e(kt)\,.
    \]
    and defining
    \begin{equation}    \label{def:hn}
        h_n(r):= n \E \min \bigg\{ 1\,, \,\frac1{\arr_n} \log_+(1/r)\bigg\} 
    \end{equation}
    where $\log_+x:=\log( 1\vee x)$,
    we have that for any $c_0\in(0,\frac12)$, 
    \begin{equation}    \label{Cov-hn}
        \Cov(X_n(s),X_n(t)) = \|\cF_{c_0/|s-t|}f\|_2^2\, h_n(|s-t|) + O_{c_0}(1).
    \end{equation}
    Moreover, $h_n$ satisfies the bounds
    \begin{equation}    \label{hn-bounds}
    -O(1)\le 
        \min\big\{ n\,,\, \log_+(1/r)\big\}
        - h_n(r) 
        \le 
        \begin{cases}
            n^{-\omega(1)} & \log\frac1r \ge n+ \sqrt{n}\log n\\
            \sqrt{n}\log n & |n-\log\frac1r|< \sqrt{n}\log n\\
            O(1) & \log\frac1r\le n-\sqrt{n}\log n
        \end{cases}
    \end{equation}
    so 
    \begin{align*}
        & -O_{c_0}(\sqrt{n}\log n) 
    \le \Cov(X_n(s),X_n(t)) 
-\|\cF_{c_0/|s-t|}f\|_2^2\,\min\Big\{ n\,,\, \log_+\frac1{|s-t|}\Big\}
\le O_{c_0}(1) 
    \end{align*}   
    and the lower bound holds with $O_{c_0}(1)$ in place of $O_{c_0}(\sqrt{n}\log n)$ unless $|s-t|= e^{-n+O(\sqrt{n}\log n)}$.
\end{prop}

\begin{remark}
In particular, if $\wh f$ is finitely supported then 
\begin{equation}    \label{Cov-finspt}
    \Cov(X_n(s),X_n(t)) = \|f\|_2^2\, h_n(|s-t|) +  O(1).
\end{equation}
Indeed, this holds for $|s-t|\gs1$ since $h_n(|s-t|)=O(1)$ in this case by \eqref{hn-bounds}, so we may assume $|s-t|<c_0^2$, say, and take $c_0$ sufficiently small that $\supp(\wh f)\subset[-c_0^{-1},c_0^{-1}]$.
\end{remark}

\begin{proof}
By rescaling $f$ we may assume without loss of generality that $\|f\|_2=1$.
By stationarity of the field it suffices to establish the claim for the case that $s=0<t$. 
We have
\begin{align}
C_n(t):= \Cov(X_n(0),X_n(t))
    &= \sum_{j=1}^{n}\E\left[f(e^{\arr_j}t+\rot_j)f(\rot_j)\right] 
    = \sum_{j=1}^{n}\E\left[g(e^{\arr_j}t)\right]\label{Cov-g}
\end{align}
where $g(v):=\E[f(v+\rot)f(\rot)]=f*\tilde f(v)$ with $\tilde f(v):=f(-v)$. The function $g$ has absolutely convergent Fourier expansion with coefficients $\widehat g(k) = |\widehat f(k)|^2=|\widehat f(-k)|^2$ for $k\in\Z$, and since $\int_{0}^{1}f(v)dv=0$ we have $\widehat g(0)=\widehat f(0)=0$. Thus,
\begin{align}\label{g-xi-ell}
    \E\left[g(e^{\arr_j}t)\right] &= 2\sum_{k\geq 1}|\widehat{f}(k)|^2\E\left[e(e^{\arr_j}k t)\right].
\end{align}
Since $\{\arr_j\}$ are Poisson arrivals with rate-$1$, when conditioned on $\arr_n$, the arrivals $\{\arr_1,\cdots,\arr_{n-1}\}$ have the same distribution as the ordered statistics of $(n-1)$ i.i.d.~uniform random variables on $[0,\arr_n]$. For $L$ uniform in $[0,\arr_n]$ we have for any $v>0$ that
\begin{align}   \label{eee-tv}
    \E\bigg[\sum_{j=1}^{n}e(e^{\arr_j}v)\,\bigg|\, \arr_n\bigg] &= (n-1)\E\left[e(e^L v)\mid \arr_n\right]+O(1)
    = \frac{n}{\arr_n}\int_0^{\arr_n}e(e^x v)\d x+O(1).
\end{align}
A standard computation with integration by parts shows that for any $M>0$ and $0<v\le \frac12$,
\[
    \int_0^M e(e^xv)dx = \min\big\{ M, \log(1/v)\big\} + O(1)
\]
while for $v\gs 1$, 
\[
    \int_0^M e(e^xv)dx = O(1).
\]
Thus, for any $c_0\in(0,\frac12)$, 
\[
\int_0^{\arr_n} e(e^x v) dx  = O_{c_0}(1) + 1_{v<c_0} \min\{ \arr_n, \log(1/v)\}\qquad a.s.
\]
Combining with \eqref{Cov-g}, \eqref{g-xi-ell} and \eqref{eee-tv}, we obtain
\begin{align}
    C_n(t) 
    &= 2\sum_{k\ge1} |\wh f(k)|^2 \E \sum_{j=1}^n e(e^{\arr_j}kt)
    =   \wt h_n(t) + O_{c_0}(1)    \label{Cn-bd1}
\end{align}
where
\[
    \wt h_n(t):= 2\sum_{1\le k<c_0/t} |\wh f(k)|^2 h_n(kt)\,.
\]
Since $h_n$ is a monotone non-increasing function,
\begin{align*}
    \wt h_n(t) \le 2\sum_{1\le k<c_0/t} |\wh f(k)|^2 h_n(t) = h_n(t) \|\cF_{<c_0/t}f\|_2^2 \,.
\end{align*}
Together with \eqref{Cn-bd1} this yields the upper bound in \eqref{Cov-hn}.
For the lower bound, using that $\min(1,\cdot)$ is Lipschitz on $\R$, we have for $kt\in (0,c_0)$ that
\[
h_n(t)\ge
h_n(kt) = n\E \min\bigg\{ 1, \, \frac1{\arr_n} \log(1/t) - \frac{\log k}{\arr_n}\bigg\} \ge h_n(t) -n\E \frac{\log k}{\arr_n} = h_n(t) - O(\log k). 
\]
Hence,
\begin{align*}
    \wt h_n(t) = 2 h_n(t) \sum_{1\le k<c_0/t} |\wh f(k)|^2  + O\bigg(  \sum_{k\ge c_0/t} |\wh f(k)|^2\log k \bigg)
    =  2 h_n(t) \sum_{1\le k<c_0/t} |\wh f(k)|^2  + O(1) 
\end{align*}
which together with \eqref{Cn-bd1} yields the lower bound in \eqref{Cov-hn}.

Turning to the bounds \eqref{hn-bounds} for $h_n$, the first follows from Jensen's inequality:
\begin{equation}    \label{hn-bd}
    h_n(t) \le n \min\Big\{ 1, \E \frac1{\arr_n} \log_+(1/t) \Big\} =  \min\big\{ n,  \log_+(1/t) \big\} + O(1).
\end{equation}
Next,
since 
$\arr_n$ is a sum of iid exponential variables with rate 1, 
standard concentration gives
\begin{equation}    \label{hn-good}
    \P(\cG^c) = n^{-\omega(1)} \qquad \text{ for } \quad \cG:=\{ |\arr_n-n|\le \sqrt{n}\log n\}\,.
\end{equation}
Writing $T:=\log_+(1/t)$, we can thus bound
\begin{align}
    h_n(t) 
    &\ge n\E 1(\cG) \min\bigg\{ 1\,,\, \frac{T}{n+\sqrt{n}\log n} \bigg\}    \notag\\
    &\ge n\min\bigg\{ 1\,,\, \frac{T}{n+\sqrt{n}\log n} \bigg\} - n^{-\omega(1)} \notag\\
    &\ge \min \{ n, T\} - Tn^{-1/2}\log n - n^{-\omega(1)} \,.   \label{hn-LBmid}
\end{align}
For $T\le n-\sqrt{n}\log n$ we can do better:
\[
    h_n(t) 
    \ge \E T\frac{n}{\arr_n} \bone(\arr_n\ge T) 
    = Tn\E \arr_n^{-1} - T \E \arr_n^{-1} \bone(\arr_n< T).
\]
Since $\E \arr_n^{-1}=1/(n-1)$ and $\E \arr_n^{-2}= O(n^{-2})$, from Cauchy--Schwarz we conclude
\begin{equation}    \label{hn-LBlow}
    h_n(t) \ge T\frac{n}{n-1} - Tn(\E \arr_n^{-2})^{1/2} \P(\arr_n<T)^{1/2} =  T -O(1)
\end{equation}
since $\{\arr_n<T\}\subset\cG^c$.
For $T\ge n+\sqrt{n}\log n$ we similarly use \eqref{hn-good} to bound
\begin{equation}    \label{hn-LBhigh}
    h_n(t) \ge n\P( \arr_n< T) = n- n^{-\omega(1)}.
\end{equation}
This completes the proof.
\end{proof}

\begin{ack}
\phantomsection
\addcontentsline{toc}{section}{Acknowledgments}
We are grateful to Ofer Zeitouni and Xaver Kriechbaum for helpful discussions, and for feedback on a preliminary version of the paper.
\end{ack}

\phantomsection

\bibliographystyle{abbrv}
\bibliography{bibfile}

\end{document}